\documentclass[10pt]{article}

\usepackage{pdfsync}
\usepackage{enumitem}
\usepackage{array}
\usepackage{amsfonts}
\usepackage{amsmath}
\usepackage{amssymb}
\usepackage{graphicx}
\usepackage{url}
\usepackage{color}
\usepackage{hyperref}
\newtheorem{thm}{Theorem}[section]
\newtheorem{cor}[thm]{Corollary}
\newtheorem{lem}[thm]{Lemma}
\newtheorem{prop}[thm]{Proposition}

\newtheorem{rem}[thm]{Remark}
\numberwithin{equation}{section}

\newcommand{\dx}{\,{\rm d}x}
\newcommand{\dy}{\,{\rm d}y}

\newcommand{\dt}{\,{\rm d}t}
\newcommand{\dtau}{\,{\rm d}\tau}
\newcommand{\rd}{{\rm d}}

\def\LL{\mathrm{L}}

\newcommand{\ka}{\overline{\kappa}}

\newcommand{\kb}{\underline{\kappa}}

\newcommand{\RR}{\mathbb{R}}
\newcommand{\NN}{\mathbb{N}}

\def\op{\mathcal{L}_{\gamma, \beta}}

\def\qed{\,\unskip\kern 6pt \penalty 500
\raise -2pt\hbox{\vrule \vbox to8pt{\hrule width 6pt
\vfill\hrule}\vrule}\par}
\definecolor{darkblue}{rgb}{0.05, .05, .65}
\definecolor{darkgreen}{rgb}{0.1, .65, .1}
\definecolor{darkred}{rgb}{0.8,0,0}

\begin{document}

\title{\textbf{Fine properties of solutions to the Cauchy problem\\ for a
 Fast Diffusion Equation with \\ Caffarelli-Kohn-Nirenberg weights}\\[7mm]}

\author{\Large Matteo Bonforte$^{\,a}$
~and~ Nikita Simonov$^{\,b}$\\[3mm]
}
\date{}

\maketitle

\begin{abstract}
We investigate fine global properties of nonnegative, integrable solutions to the Cauchy problem for the Fast Diffusion Equation with weights (WFDE) $u_t=|x|^\gamma\mathrm{div}\left(|x|^{-\beta}\nabla u^m\right)$ posed on $(0,+\infty)\times\mathbb{R}^d$, with $d\ge 3$, in the so-called \emph{good fast diffusion range} $m_c<m<1$,  within the range of parameters $\gamma, \beta$ which is optimal for the validity of the so-called Caffarelli-Kohn-Nirenberg inequalities.

It is a natural question to ask in which sense such solutions behave like the Barenblatt $\mathfrak{B}$ (fundamental solution): for instance, asymptotic convergence, i.e.   $\|u(t)-\mathfrak{B}(t)\|_{{\rm L}^p(\mathbb{R}^d)}\xrightarrow[]{t\to\infty}0$, is well known for all $1\le p\le \infty$, while only few partial results tackle a finer analysis of the tail behaviour. \textit{We characterize the maximal set of data $\mathcal{X}\subset{\rm L}^1_+(\mathbb{R}^d)$ that produces solutions which are pointwise trapped between two Barenblatt }\textsl{(Global Harnack Principle)}, and \textsl{uniformly converge  in relative error} (UREC), i.e.
${\rm d}_\infty(u(t))=\|u(t)/\mathcal{B}(t)-1\|_{{\rm L}^\infty(\mathbb{R}^d)}\xrightarrow[]{t\to\infty}0$. Such characterization is in terms of an integral condition on $u(t=0)$.

To the best of our knowledge, analogous issues for the linear heat equation $m=1$, do not possess such clear answers, only partial results. Our  characterization is also new for the classical, non-weighted, FDE.
We are able to provide \textit{minimal rates of convergence }to $\mathcal{B}$ in different norms. Such rates are almost optimal in the non weighted case, and become optimal for radial solutions.
To complete the panorama, we show that solutions with data in  ${\rm L}^1_+(\mathbb{R}^d)\setminus\mathcal{X}$, preserve the same ``fat'' spatial tail for all times, hence UREC fails and ${\rm d}_\infty(u(t))\!=\!\infty$, even if $\|u(t)-\mathcal{B}(t)\|_{{\rm L}^1(\mathbb{R}^d)}\xrightarrow[]{t\to\infty}0$.

\vspace{3cm}
~\,.
\end{abstract}

\noindent {\sc Keywords. }Fast diffusion equation; Caffarelli-Kohn-Nirenberg weights; Global Harnack inequalities; Tail behaviour; Asymptotic behaviour.

\noindent{\sc Mathematics Subject Classification}. 35B40, 35B45, 35K55, 35K67, 35K65.

\vfill
\begin{itemize}[leftmargin=*]\itemsep2pt \parskip3pt \parsep0pt
\item[(a)] Departamento de Matem\'{a}ticas, Universidad Aut\'{o}noma de Madrid,\\
ICMAT - Instituto de Ciencias Matem\'{a}ticas, CSIC-UAM-UC3M-UCM, Calle Nicol\'{a}s Cabrera 13-15\\
Campus de Cantoblanco, 28049 Madrid, Spain.\\
E-mail:\texttt{~matteo.bonforte@uam.es }\;\;
Web-page:\texttt{~http://verso.mat.uam.es/\~\,matteo.bonforte/}
\item[(b)] Departamento de Matem\'{a}ticas, Universidad Aut\'{o}noma de Madrid, Campus de Cantoblanco, 28049 Madrid, Spain, and
CEREMADE (CNRS UMR no 7534), PSL University, Universit\'e Paris-Dauphine, Place de Lattre de
Tassigny, 75775 Paris 16, France  \\
E-mail:\texttt{~nikita.simonov@uam.es }\;\;
Web-page:\texttt{~https://sites.google.com/view/simonovnikita/}
\end{itemize}

\newpage
\small

\tableofcontents

\normalsize

\newpage

\section{Introduction and Main Results}\label{sec:intro}
The purpose of this paper is to investigate fine decay properties of solutions to the Cauchy problem for the Fast Diffusion Equation posed on $(0,+\infty)\times \RR^d$
\[
\partial_t u = \Delta u^m\qquad\mbox{with $m\in (0,1)$}\,,
\]
which is a classical prototype of singular nonlinear diffusion in homogeneous media, see~\cite{Vazquez2006,Vazquez2007}. Indeed, our results also cover a more general case, in which Caffarelli-Kohn-Nirenberg type weights are allowed, including also the case of inhomogeneous media, for more details, see~\cite{Rosenau1982,Kamin1982} and Subsection~\ref{ssec:setup}. Our goal is to understand the quantitative tail behaviour of nonnegative solutions, depending on the behaviour at infinity of the initial datum. In order to focus on the main questions and answers, we explain here a simplified version of our results.   Some of them, in the non-weighted case, were already known (mostly in a non-sharp form) and in this case we provide a new proof, see for instance \cite{Carrillo2001,Carrillo2002,DelPino2002,Carrillo2003,Vazquez2003,Vazquez2006,Carrillo2006,Bonforte2006,McCann2006,Kim2006} and also \cite{Friedman1980,Herrero1985}.

It is well known that nonnegative and integrable solutions tend to behave like the Barenblatt profile $\mathfrak{ B}$ with the same mass \cite{Vazquez2003,Vazquez2006,Vazquez2007}. However, the issue of making a precise and quantitative statement about such ``similar behaviour'' presents serious difficulties.  Nonnegative integrable solutions  (namely, solutions in $\LL^1_{+}(\RR^d)$) convergence to $\mathfrak{B}$ in different norms, $\|u-\mathfrak{B}\|\xrightarrow[]{t\to\infty}0$. However, none of these convergences provide enough information about the tails. Here we explore finer properties of such solutions to give an answer to the following question:
\begin{center}\textit{$Q_1$: Do  solutions in $\LL^1_{+}(\RR^d)$ have  the same tail as the Barenblatt (fundamental solution)?}
\end{center}

When the answer to this question is positive we can ask for a finer convergence, namely whether or not the quotient $u/\mathfrak{B}\xrightarrow[]{t\to\infty}1$  uniformly in $\RR^d$. We call this Uniform Convergence in Relative Error (UREC), and we can state the main question in this direction as follows:
\begin{center}\textit{$Q_2$: Do  solution in $\LL^1_{+}(\RR^d)$  behave asymptotically as the Barenblatt, uniformly in relative error? }
\end{center}

In this paper, we completely answer to both $Q_1$ and $Q_2$ in the so-called good fast diffusion range, when $m\in (\frac{d-2}{d},1)$. This is the natural range of parameters to deal with these questions, since the Barenblatt profile $\mathfrak{B}$ essentially represents the asymptotic behaviour of all integrable solutions and the mass is preserved along the flow. In order to provide an answer to the above questions, we split the cone of nonnegative integrable initial data $\LL^1_+(\RR^d)$ in two disjoint subspaces $\mathcal{X}\bigsqcup\mathcal{X}^c=\LL^1_+(\RR^d)$,  where $\mathcal{X}$ is a set of functions satisfying a suitable tail condition specified below.

The main results of this paper can be roughly explained as a quantitative version of the following fact. On one hand, the answer to both $Q_1$ and $Q_2$ is affirmative  \textit{if and only if }$u_0\in \mathcal{X}$, (Thm. \ref{thm.ghp}, \ref{equivalent.condition} and Section \ref{section:GHP}). On the other hand, the answer to both  $Q_1$ and $Q_2$ is No \textit{if and only if }$u_0\in \mathcal{X}^c$ (Thm. \ref{thm.GGHP} and Section \ref{Sec4}).

It is remarkable that such a complete answer can be given for a nonlinear equation, while, in the - a priori simpler - linear heat equation $m=1$, things are not so clear.  Partial - non sharp - answers to $Q_1$ can be deduced from the representation formula, an extremely useful tool that we do not have at our disposal in the nonlinear case. As for $Q_2$, the question seems to be completely open: to the best of our knowledge, there is no characterization of the class of initial data for which the corresponding solution converge to the Gaussian (with the same mass) uniformly in relative error. Some examples, in the negative direction, are shown in \cite{Vazquez2017}.

Our results are sharp and turn into an explicit characterization of the ``Tail Condition'' that the initial datum has to satisfy to be in $\mathcal{X}$ (hence answering yes to $Q_1$ and $Q_2$), which amounts to requiring that
\[
\sup_{R > 0} R^{\frac{2}{1-m}-d}\int_{B^{c}_R(0)}|f(x)|\dx  <\infty\,,\qquad\mbox{or equivalently}\qquad
\int_{B_{|x|/2}(x)}|f(y)|\dy = \mathrm{O}\left(|x|^{d-\frac{2}{1-m}}\right)\,.\vspace{-2mm}
\]
The proof of the equivalence of the above two conditions is not trivial, indeed it requires one of our main results, Theorem \ref{thm.ghp}; see Section \ref{ssect5.1}.
The latter condition was introduced by V\'azquez in 2003~\cite{Vazquez2003} to give a positive answer to $Q_2$, we show here that \emph{a posteriori} it was the sharp one. Notice that this condition allows for a wider class of data than the (non sharp) pointwise condition used in \cite{Bonforte2006, Vazquez2006}, namely $u_0(x)\lesssim |x|^{-\frac{2}{1-m}}$, see Section \ref{ssec.BadGuy}. We also show that the above condition, when fulfilled by the data, is enough to prove  polynomial rates of convergence in several norms. In the radial case, we deduce sharp rates of convergence in uniform relative error and we provide an answer to a question left open by Carrillo and V\'azquez in \cite{Carrillo2003}, see Remark \ref{Carrillo-Vazquez}.

Concerning initial data in $\mathcal{X}^c$: we show the existence of a class of solutions which exhibits, for all times, a \emph{fat tail} (bigger then the Barenblatt's). This is done by constructing explicit sub and super solutions. Such class provides the negative answer to both $Q_1$ and $Q_2$. Furthermore, we show that in $\mathcal{X}^c$ no (power-like) rate of convergence to the Barenblatt profile is possible, see Theorem~\ref{no.rates}, Subsection~\ref{ssec:no.rates}.

In the rest of this section, we set the problem in its whole generality, also including equations with Caffarelli-Kohn-Nirenberg type weights, and give precise statements of our results.

\subsection{The Setup of the problem and precise statement of the Main Results}\label{ssec:setup}
In this paper we study the following Cauchy-Problem for the Weighted Fast Diffusion Equation (WFDE)
\begin{equation}\label{cauchy.problem}\tag{CP}
\begin{cases}
\begin{aligned}
&\partial_t u  = |x|^{\gamma}\nabla \cdot \left(|x|^{-\beta}\nabla u^m \right)  \qquad\, &\text{in} \left(0, \infty \right)\times \RR^d, \\
&u\left(0,x\right) =u_0(x) \qquad\, &\text{in $\RR^d$}. \\
\end{aligned}
\end{cases}
\end{equation}
where the parameters $d, \gamma, \beta$ are as follows
\[
d\ge3\,,\qquad \gamma< d\,,\qquad\mbox{and}\qquad\gamma -2 < \beta \leq \gamma(d-2)/d\,.
\]
This is a natural restriction since it represents the optimal domain of validity of the so-called Caffarelli-Kohn-Nirenberg inequalities, see~\cite{Caffarelli1984, Bonforte2019a}.
The exponent $m$ is in the so-called \textit{good fast diffusive range}, namely
\[
m \in \left(m_c, 1\right)\qquad\mbox{where}\qquad m_c:= \frac{d-2-\beta}{d-\gamma}\,.
\]
From now on we will fix the parameters $d,m,\gamma,\beta$ as above (unless explicitly stated). \\

\noindent\textit{Modelling and related results.} The problem \eqref{cauchy.problem} was introduced in the 80s by Kamin and Roseau  to model singular/degenerate diffusion in inhomogeneous media, see~\cite{Kamin1981,Rosenau1982,Kamin1982}. Since then, there has been a systematic study of similar equations, mostly in the case $m\ge 1$ and/or with only one weight, see~\cite{Abdellaoui2004,Bonafede1999,DallAglio2004,Eidus1990,Eidus1994,Grillo2013,Grillo2015,Iagar2014,Kamin1998,Kamin2010,Nieto2013,Reyes2009,Reyes2008,Reyes2006,Reyes2006a}.
Recently,~\eqref{cauchy.problem} has proven to be an essential tool in the study of symmetry/symmetry breaking phenomena in Caffarelli-Kohn-Nirenberg inequalities, see~\cite{Bonforte2017a,Bonforte2017,Catrina2001,Dolbeault2015,Dolbeault2016,Dolbeault2017,Dolbeault2011,Dolbeault2017a,Felli2003}. Several intriguing connections between nonlinear diffusions on Riemannian manifolds and weighted parabolic equations were explored in~\cite{Bonforte2008,Bonforte2010,Grillo2014,Ishige1999,Ishige2001,Ishige1998,Vazquez2015}. \\

\noindent\textbf{Existence, uniqueness, comparison and mass conservation. }The basic theory is well established: existence, uniqueness and comparison for nonnegative and bounded integrable data is well known, see Section 2.2 of \cite{Bonforte2017a}, where it can also be found a suitable definition of weak solutions (cf. also Definition 1.1 of \cite{Bonforte2019a}). In view of the smoothing effects of \cite{Bonforte2019a}, it is straightforward to extend those results to weak solutions corresponding to merely integrable (and possibly unbounded) data
\begin{equation}\label{definition.L1}
u_0 \in \LL^1_{\gamma, +}(\RR^d)=\left\{u_0:\RR^d\rightarrow\RR : u_0\ge0\,, \int_{\RR^d}u_0\,|x|^{-\gamma}\dx <\infty\right\}\,.
\end{equation}
We refrain from giving further details that would involve weighted Sobolev spaces which we never use in this paper and we choose not to define here. What we want to emphasize, is that  data in $\LL^1_{\gamma, +}$ produce solutions that turn out to be bounded, positive and regular, (at least H\"older continuous, see also Appendix \ref{app.holder.cont} and \cite{Bonforte2019a}) and that solutions considered in this paper possess enough regularity to guarantee the validity of all the calculations performed here.
We also recall that in the good fast diffusive range, nonnegative integrable solutions conserve mass along the flow,
\[
M(t):=\int_{\RR^d}u(t,x)|x|^{-\gamma}\dx = \int_{\RR^d}u_0(x)|x|^{-\gamma}\dx=:M\qquad\mbox{for any $t>0$.}
\]
For a proof, see Section 2.2 of \cite{Bonforte2017} and Proposition 2.4 of \cite{Bonforte2019a} and the Remark thereafter.

\noindent\textbf{The fundamental solution, }is of \emph{self-similar type} and it is often called Barenblatt solution:
\begin{equation}\label{barenblatt}\tag{B}
\mathfrak{B}(t+T, x; M)= \frac{\zeta^{d-\gamma}}{R_\star(t+ T)^{d-\gamma}}\,\mathcal{B}_M\left(\frac{\zeta\,x}{R_\star(t+ T)}\right)= \frac{(t+T)^{\frac{1}{1-m}}}{\left[b_0\frac{\left(t+T\right)^{\sigma\vartheta}}{M^{\sigma\vartheta(1-m)}}+b_1 |x|^{\sigma}\right]^\frac{1}{1-m}}\,,
\end{equation}
where
\begin{equation}\label{parameters.sigma.theta}
\sigma:= 2+\beta-\gamma\,, \qquad \frac{1}{ \vartheta}=(d-\gamma)(m-m_c)\,, \qquad\zeta^\frac{1}{\vartheta}=\frac{1-m}{\sigma\,m}\,,\qquad\mbox{and}\qquad R_\star(t) = \left(\frac{t}{\vartheta}\right)^\vartheta\,,
\end{equation}
the parameter $M$ is the \emph{mass} of the solution, $T$ is a free parameter and $b_0, b_1$ are constants which depends on $m, d, \gamma, \beta$. The profile $\mathcal{B}_M$ is given by
\begin{equation}\label{barenblatt.profile.stationary}
\mathcal{B}_M= (C(M)+|x|^\sigma)^\frac{1}{m-1}\,
\end{equation}
where $C(M)$ depends on $M, d ,m, \gamma, \beta$, and has an explicit expression, see Appendix~\ref{appendix:mass}.
In what follows we shall frequently use the solution~\eqref{barenblatt} with the parameter $T=0$. Recall that by the very definition of fundamental solutions we have $\mathfrak{B}(0,x; M)=M\delta_0$, in the sense of measures: the computation goes as for the standard FDE, and the extra weight $|x|^{-\gamma}$ does not cause any problem. Also, we will sometimes drop the dependence on the $x$ variable and write $\mathfrak{B}(t;M)$ or $\mathfrak{B}(t, \cdot; M)$ when no confusion arises.

\noindent\textbf{The Tail Condition. }We say that $f \in \LL^1_\gamma(\RR^d)$ satisfies the \emph{tail-condition}  -or equivalently that $f \in \mathcal{X}$- if
\begin{equation}\label{tail.condition}\tag{TC}
|f|_{\mathcal{X}}:=\sup_{R > 0} R^{\frac{2+\beta-\gamma}{1-m}-(d-\gamma)}\int_{B^{c}_R(0)}|f(x)||x|^{-\gamma}\dx  <\infty\,.
\end{equation}
Recall that since $m \in (m_c, 1)$ we have $\frac{2+\beta-\gamma}{1-m}-(d-\gamma) > 0$. It is easily seen that $|\cdot|_\mathcal{X}$ is a norm. Intuitively the quantity $|f|_\mathcal{X}$ measures \emph{how fast} the function $f$ decays at $\infty$ relatively to the decay of the Barenblatt profile $\mathfrak{B}_{M}$.  We now introduce a subspace of $\LL^1_\gamma(\RR^d)$ of functions that satisfy the tail condition~\eqref{tail.condition}, that will play a key role in the rest of the paper:
\begin{equation}\label{dacay.space}
\mathcal{X} := \{ u \in \LL^1_{\gamma}(\RR^d) :  |u|_{\mathcal{X}} < +\infty\}.
\end{equation}

We adapt to our setting an alternative tail condition proposed by Vazquez~\cite{Vazquez2003}:  we say that $f\in \LL^1_\gamma(\RR^d)$ satisfies~\eqref{tail.condition.1} if
\begin{equation}\label{tail.condition.1}\tag{TC'}
\int_{B_{\frac{|x|}{2}}(x)}|f(y)||y|^{-\gamma}\dy = \mathrm{O}\left(|x|^{d-\gamma-\frac{2+\beta-\gamma}{1-m}}\right)\,.
\end{equation}
We will show in section \ref{ssect5.1} that~\eqref{tail.condition} and~\eqref{tail.condition.1} are indeed equivalent.

We will provide now a precise and sharp answer to questions $Q_1$ and $Q_2$, in the form of our main results.       \\

\noindent\textbf{The space $\mathcal{X}$: affirmative answer to $Q_1$ and $Q_2$, and a characterization. }

As we already explained in the Introduction, the answer to both $Q_1$ and $Q_2$ are affirmative \emph{if and only if} the initial data is in $\mathcal{X}$.  The main tool in providing such answers is the so called \emph{Global Harnack Principle} (GHP): a lower and upper bound in terms of Barenblatt profiles, see Theorem~\ref{thm.ghp} below. The GHP provides a complete answer both to $Q_1$ and, surprisingly, also  to $Q_2$, as we shall see later. In the non-weighted case, the GHP in the form of Theorem \ref{thm.ghp} was introduced in~\cite{Bonforte2006} (under the stronger pointwise assumption $u_0(x)\lesssim |x|^{-\frac{2}{1-m}}$) and was inspired by the pioneering results of~\cite{Vazquez2003}, in which condition \ref{tail.condition.1} was introduced.\\
Our main contribution in this case consists in the  \emph{characterization} of the maximal set $\mathcal{X}$ of initial data that generate solutions satisfying the GHP.

We shall see in Section \ref{X-flow} that the space $\mathcal{X}$ is \textit{invariant} under the WFDE-flow: indeed $u(t)\in \mathcal{X}$ if and only if $u(0)\in \mathcal{X}$, and the same holds for $\mathcal{X}^c$, see Proposition \ref{curve.X.prop} and Theorems \ref{thm.ghp} and \ref{thm.GGHP}.

\begin{thm}[Characterization of the GHP]\label{thm.ghp}Let $m\in \left(m_c, 1\right)$ and let $u$ be a solution to \eqref{cauchy.problem} with $u_0 \in \LL^1_{\gamma,+}(\RR^d)$.  Then, for any $t_0>0$ there exists $\overline{\tau}, \underline{\tau} > 0$  and $\overline{M}, \underline{M}>0$ such that
\begin{equation}\label{ghp.inequality}
\mathfrak{B}(t-\underline{\tau},x; \underline{M}) \leq u(t, x) \leq \mathfrak{B}(t+\overline{\tau},x; \overline{M})\,,\qquad\mbox{for any $x \in \RR^d$ and $t \geq t_0$ },
\end{equation}
if and only if
\[
\mbox{$ u_0 \in \mathcal{X}\setminus\{0\}\,.$}
\]
\end{thm}

\begin{rem}\rm
The proof of the above result is quantitative and provides explicit expressions for $\underline{\tau}, \overline{\tau}, \underline{M}, \overline{M}$. It follows by combining the upper bound of Theorem \ref{thm.upper} with the lower bound of Theorem \ref{lower.estimate.theorem}. For the upper bound the hypothesis $0\le u_0 \in \mathcal{X}$ is strictly necessary. Indeed, for data $u_0 \notin \mathcal{X}$ we are able to construct explicit (sub)solutions that provide precise lower bounds that clearly contradict the upper bound of formula \eqref{ghp.inequality}. More precisely, for any $t>0$ and for any $x \in \RR^d$ we have that
 \begin{equation*}
u(t,x) \geq \frac{1}{\left(D(t) + |x|^\sigma\right)^{\frac{1}{1-m}-\varepsilon}} \gvertneqq \mathfrak{B}(t, x; M)\,,
 \end{equation*}
where $\varepsilon>0$ is small, and $D(t) \sim t^{\frac{2}{\varepsilon (1-m)}}$.\\
On the other hand, such  hypothesis is not necessary for the lower bound of formula \eqref{ghp.inequality}: indeed, lower bounds hold for any data $0\leq u_0 \in \LL^1_{\rm \gamma, loc}(\RR^d)$, and produce a \emph{minimal lower tail}, see  Theorem \ref{lower.estimate.theorem} and Corollary \ref{herrero.pierre.limit.corollary}. This provides a partial answer to $Q_1$.
\end{rem}
Let us turn our attention to question $Q_2$. The convergence of solution  to the Barenblatt profile has been studied by many researchers and under different sets of assumptions, especially in the non-weighted case $\gamma=\beta=0$, see for instance  \cite{Friedman1980,Vazquez2003,Carrillo2003,Carrillo2006,Kim2006,McCann2006,Blanchet2009,Bonforte2010} and references therein. We will discuss now some of the existing results which are strictly related to ours, but all in the non-weighted case. To the best of our knowledge, no results about the weighted case are present in the literature, except some partial results of \cite{Bonforte2017a,Bonforte2017,Bonforte2019a}. In \cite{Friedman1980} the authors proved uniform convergence on expanding sets of the form $|x|\le C t^\vartheta$, namely
\begin{equation}\label{local.convergence.relative.error.noweights}
\lim_{t\rightarrow 0}\sup_{x \in \{|x| \le C\, t^\vartheta\}} \left|\frac{u(t,x)-\mathfrak{B}(t, x; M)}{\mathfrak{B}(t, x; M)} \right|=0\,,
\end{equation}
under the condition $u_0 \in \LL^1(\RR^d)\cap \LL^2(\RR^d)$. We will prove an analogous result in the weighted case, see Theorem~\ref{local.convergence.relative.error.thm}. Lately V\'azquez  in \cite{Vazquez2003} has completed the proof of the previous result for the whole class of positive initial data which belongs to $\LL^1(\RR^d)$, he also shows uniform convergence in $\LL^\infty(\RR^d)$ and in  $\LL^1(\RR^d)$. In~\cite{Vazquez2003} V\'{a}zquez proved that UREC takes place for all data which satisfies the pointwise condition $u_0(x)\lesssim |x|^{-\frac{2}{1-m}}$, and he also introduces~\eqref{tail.condition.1}. In 2003, Carrillo and V\'{a}zquez in~\cite{Carrillo2003} obtain the estimates
\begin{equation}\label{open.problem.carrillo.vazquez}
\sup_{x\in\RR^d}\Big\|\frac{u(t,x)}{\mathfrak B (t, x; M)}-1\Big\|_{\LL^\infty(\RR^d)} \le \frac{C(u_0)}{t}\,,
\end{equation}
for radial initial data which satisfies the pointwise condition $u_0(x)\lesssim |x|^{-\frac{2}{1-m}}$. An intriguing open question was left in \cite[pag. 1027]{Carrillo2003}: to extend the validity of~\eqref{open.problem.carrillo.vazquez} to a larger class of initial data. The question was partially answered in~\cite{Kim2006,Blanchet2009,Bonforte2017} in some non-optimal classes of data, possibly non radial.

Our main contribution in this paper is to \emph{characterize} the maximal set $\mathcal{X}$ of initial data whose solution converge to the Barenblatt profile \emph{uniformly in relative error}.

\begin{thm}[Characterization of the UREC]\label{equivalent.condition}
Let $m\in \left(m_c, 1\right)$ and let $u$ be a solution to \eqref{cauchy.problem} with initial data $ u_0 \in \LL^1_{\gamma,+}(\RR^d)$ and  $M=\|u_0\|_{\LL^1_{\gamma}(\RR^d)}$.  Then,
\begin{equation}\label{global.convergence.relative.error.equivalence}
\lim_{t\rightarrow\infty}\Big\|\frac{u(t,x)}{\mathfrak{B}(t, x; M)}-1\Big\|_{\LL^\infty(\RR^d)} =0
\end{equation}
if and only if
\[
u_0 \in \mathcal{X}\setminus\{ 0 \}
\]
\end{thm}
\begin{rem}\label{Carrillo-Vazquez}Sharp convergence rates for radial solutions. \rm
We notice here that if $f\in \mathcal{X}_{rad}$, the class of radial functions in $\mathcal{X}$, then it does not necessarily satisfy the pointwise condition $f\lesssim |x|^{-\frac{2}{1-m}}$, see Section \ref{ssec.BadGuy}. Hence, our Theorem \ref{rates.radial.data} which shows the validity of \eqref{open.problem.carrillo.vazquez} for any $u_0\in \mathcal{X}_{rad}$, provides a sharp answer to the question raised by Carrillo and Vazquez. The maximality of $\mathcal{X}$ is guaranteed by Theorem \ref{equivalent.condition}, indeed, if $u_0\not\in \mathcal{X}$ then the limit \eqref{global.convergence.relative.error.equivalence} is infinite, see Proposition \ref{curve.X.prop}.

In Theorem \ref{almost.universal.rates.noweights}, we provide \textit{almost optimal rates }of convergence for all data in $\mathcal{X}$ in the non weighted case, valid also for non-radial solutions. Analogously,  Theorem \ref{convergence.rates.weights} shows  minimal rates in the weighted case. Sharp rates of convergence can be obtained under more restrictive assumptions (but for the whole range $m<1$): this happens if the initial datum is trapped between two Barenblatt solutions with exactly the same tail (which is stronger than the GHP of Theorem \ref{thm.ghp}). We refer to \cite{Blanchet2009,Bonforte2010,Bonforte2010c,Bonforte2017a} and references therein for an overview of previous results; see also \cite{Bonforte2010c} for a brief historical overview.
\end{rem}

\begin{figure}[h]
\begin{minipage}[t]{0.5\textwidth}
   \includegraphics[width=1\textwidth]{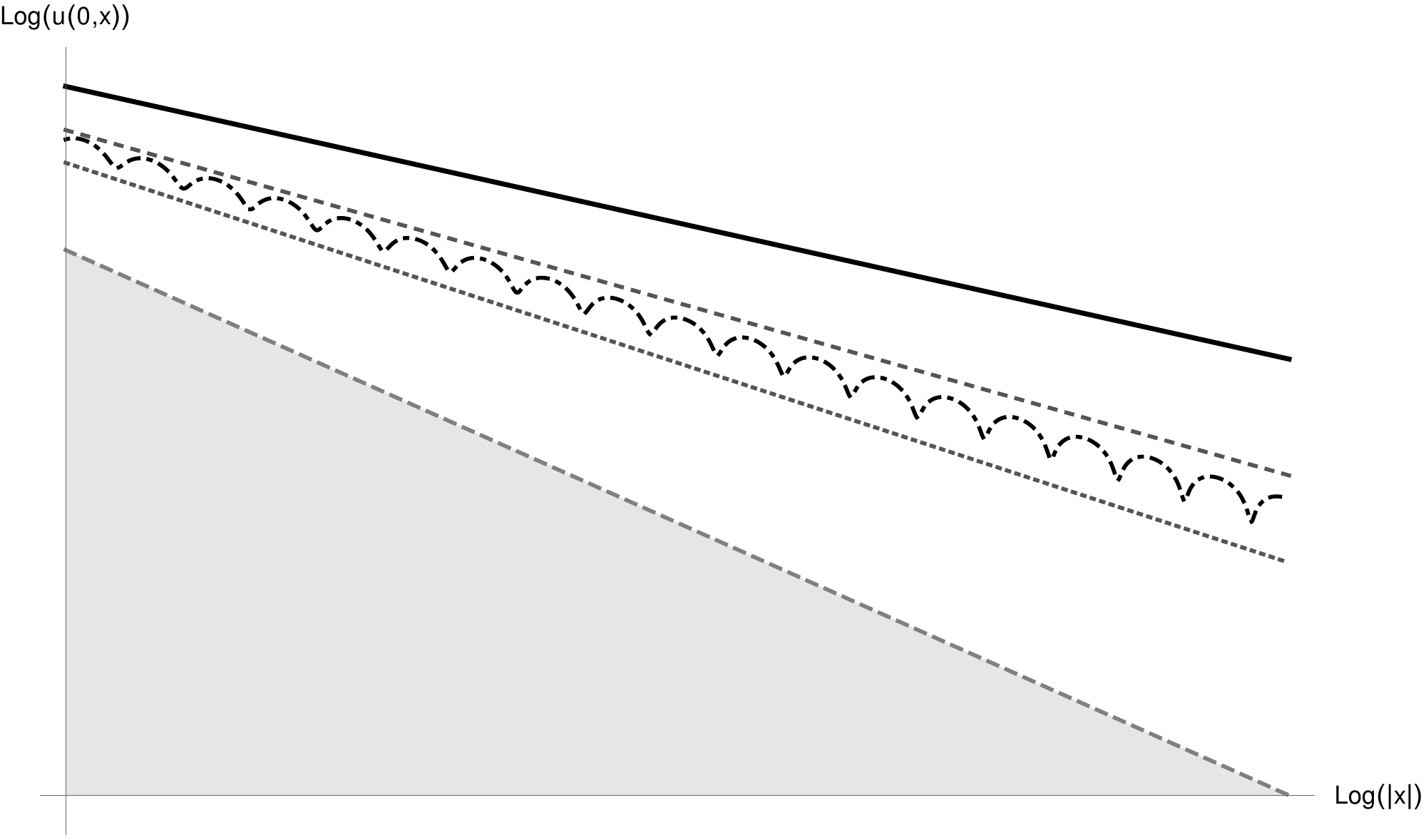}
   \caption{The picture represents different possible power-like tail behaviour of an initial data in dimension $d=3$, $\gamma=\beta=0$ and $m=2/3$. Every line represents a power-like behaviour. The dashed ``sinusoidal'' curve represents a function whose behaviour is trapped between two different tail powers. The thick line represents the power $|x|^{-3}$, the dashed line $|x|^{-3.8}$, the dotted line $|x|^{-4.4}$ and the grey region represents any decay below the Barenblatt's decay ($|x|^{-6}$).  In both plots $|x|\in \big[10^3,10^{8}\big]$.}\label{fig:picture.global.initial.fig}
\end{minipage}
   \hspace{0.2cm}
\begin{minipage}[t]{0.5\textwidth}
   \includegraphics[width=1\textwidth]{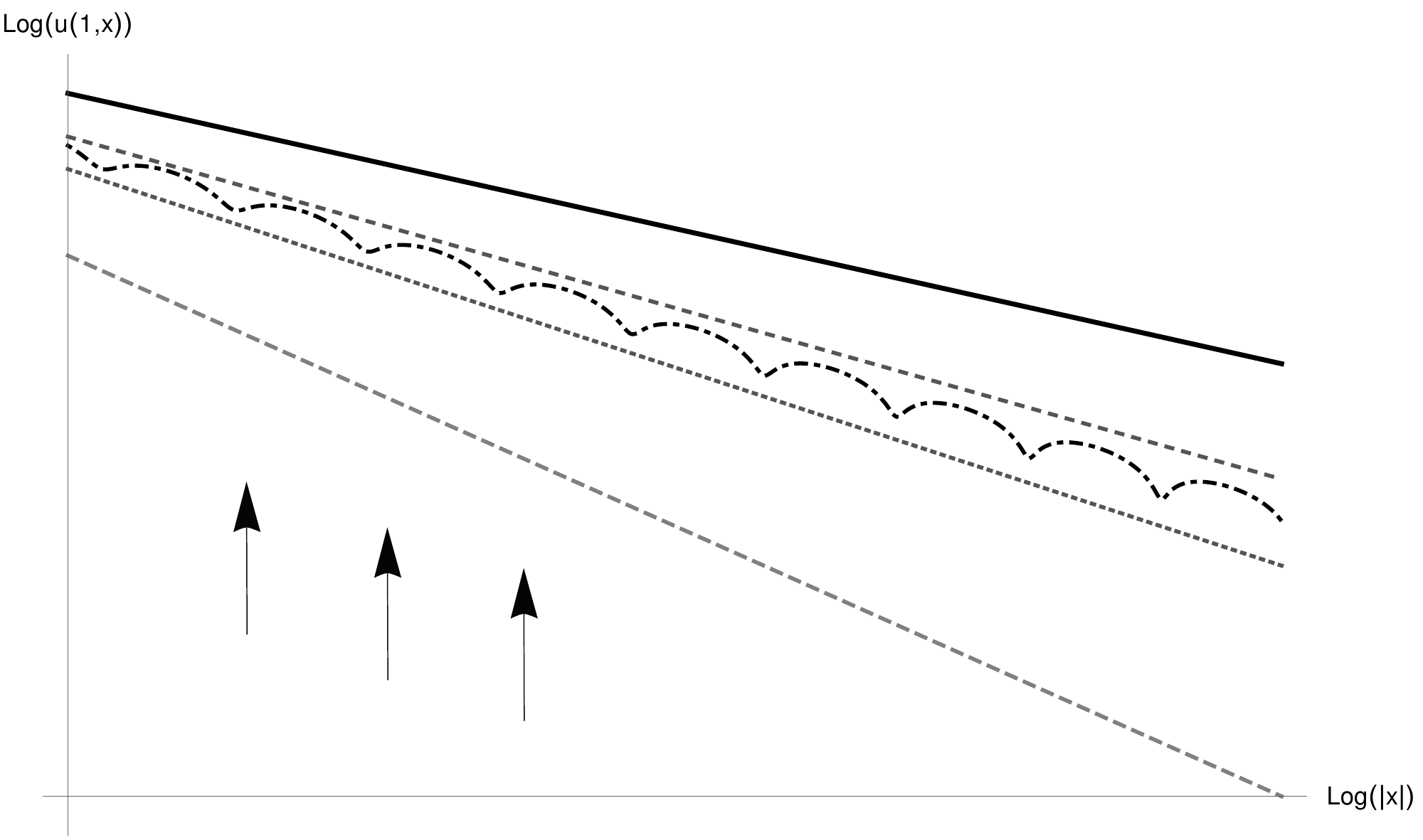}
   \caption{General Panorama for $\LL^1_{\gamma, +}$ solutions. The picture represents different possible power-like tail behaviours in dimension $d=3$, $\gamma=\beta=0$ and $m=2/3$ for a solution $u(t,x)$ at time $t=1$. We can appreciate that: (i) there are no solutions below the lowest line $|x|^{-6}$, the one corresponding to the Barenblatt behaviour; (ii) if the tail of the initial data is a line above the Barenblatt then the solution preserves the same tail behaviour for all times; (iii) the dashed curve remains trapped between the same initial power-like behaviors. }\label{fig:picture.global.time1.fig}
\end{minipage}
\end{figure}

\noindent\textbf{The space $\mathcal{X}^c$: negative answer to $Q_1$ and $Q_2$. Counterexamples. }In order to complete the panorama, we still have to answer the next natural question: \textit{what happens to the solutions with data in $\mathcal{X}^c$? }On one hand, the space $\mathcal{X}^c$ is also \textit{invariant} under the WFDE-flow: indeed, $u(t)\in \mathcal{X}^c$ if and only if $u_0\in \mathcal{X}^c$; moreover the uniform relative error \eqref{RE.Xc} is always infinite, see Proposition \ref{curve.X.prop}. As a consequence, answer to $Q_1$ and $Q_2$ is definitively negative in $\mathcal{X}^c$.
On the other hand, we will show that -somehow stable- \textit{anomalous tail behaviour }can happen in this case.
Let us begin with an illuminating example in the simplest possible case, when $\gamma=\beta=0$.
Let $m> \frac{d}{d+2}$,
consider the solution $w(t,x)$ with initial data
\[
w_0 = \frac{1}{\left(1+|x|^2\right)^\frac{m}{1-m}}\,.
\]
It is clear that for $w_0$ does not satisfy the assumption of Theorem~\ref{thm.upper} and, for $|x|$ large enough, we have that $w_0(x)> \mathfrak B(t_0, x;M)$ for any $t_0, M>0 $. However, $w_0 \in \LL^1(\RR^d)$ whenever $m > \frac{d}{d+2}$. The tail-behaviour of $w(t,x)$ is strongly different from the Barenblatt profiles, this can be better appreciated in logarithmic scale see for instance Figure~\ref{fig:picture.global.time1.fig}, indeed for all $t>0$
\[
\frac{1}{\left(\,(t+1)^\frac{1}{1-m}+|x|^2\right)^\frac{m}{1-m}}\lesssim w(t,x) \lesssim \frac{(1+t)^\frac{m}{1-m}}{\left(1+t+|x|^2\right)^\frac{m}{1-m}}\,,\qquad\mbox{for all }x\in \RR^d\,.
\]
The above inequality  gives us remarkable insights about the long time behaviour of the solution $w(t,x)$. First, for any time $t>0$, $w(t,x)$ has a power-like behaviour at infinity different from the Barenblatt one, namely as $|x|\rightarrow \infty$ we have that $w(t,x) \sim |x|^\frac{-2m}{1-m}$ versus  $\mathfrak B(t, x;M)\sim |x|^\frac{-2}{1-m}$. The upper part of GHP fails outside a space-time region that we explicitly identify, as a consequence of this anomalous tail behaviour, indeed
\begin{equation}\label{RE.Xc}
\sup_{x\in \RR^d}\left|\frac{w(t,x)}{\mathfrak B (t, x; M)}-1\right|  = \infty\,,
\end{equation}
where $\mathfrak B (t, x; M)$ has the same mass of $w(t,x)$. The same considerations apply by replacing $\mathfrak B (t, x; M)$ with any other Barenblatt solution. Obviously, \emph{uniform converge in relative-error (UREC) fails}.

The anomalous behaviour found in this particular example is indeed shared by an entire class of solutions. We prove here a \emph{generalized} version of the GHP, valid for initial data decaying slower than the Barenblatt profile. The proof is based on the construction of two  families of sub and super solutions.
We will cover all admissible $\gamma, \beta$ and $m \in (m_c, 1)$, extending  the above considerations to the weighted case, as in the following
\begin{thm}[Generalized Global Harnack Principle]\label{thm.GGHP}
Let $m \in (m_c, 1)$,  $\varepsilon \in (0, \frac{2}{1-m}-\frac{2}{\sigma}(d-\gamma))$ and $\alpha=\frac{1}{1-m} - \frac{\varepsilon}{2}>0$. Assume that the initial data $u_0$ satisfies
\[
\frac{A}{\left(\underline{t}^\frac{1}{1-\alpha(1-m)}+B\,|x|^\sigma\right)^\alpha}\le u_0(x) \le \frac{E\,\overline{t}^\sigma}{\left(\overline{t} + F\,|x|^\sigma\right)^\alpha}\,,
\]
for some  $A, B, E, F ,\underline{t}, \overline{t}>0$. Then for any $t>0$ we have that
\[
\underline{V}(t,x):= \frac{A}{(D(t)+B |x|^\sigma)^{\alpha}}  \le u(t,x) \le \frac{E\,G(t)^\alpha}{(G(t)+F |x|^\sigma)^{\alpha}}=:\overline{V}(t,x)\,
\]
where
\[
D(t):= \left(\sigma\,A^{m-1}\,m\,B\,(d-\gamma)\, (1-\alpha(1-m)) \,t + \underline{t}\right)^{\frac{1}{1-\alpha(1-m)}}\quad\mbox{and}\quad G(t):=\overline{t} + H\,t,\,
\]
where $H \ge m\,\sigma\,F^2\,E^{m-1}\,\left(2+\beta-d+\sigma\,\alpha\,m\right)$.
\end{thm}
The proof of the above Theorem  is just the combination of the results of Propositions~\ref{subsolution.prop} and \ref{supersolution.prop}.

\noindent\textbf{Remark. }The above Theorem shows that an initial ``fat-tail'' is preserved for all times. This marks a clear difference between the ``good'' range $(m_c, 1)$ and the very fast diffusive range $(0, m_c)$: in the latter case there can be
solutions with a power-like tails which change with time, see for instance \cite{Daskalopoulos2008}.


\noindent\textbf{The space $\LL^1_{\gamma, +}$. General Picture. }We provide here a general picture for solutions in $\LL^1_{\gamma}$, depending on the spatial decay of its initial data. This is better understood in the following $\log\left(u(x)\right)-\log\left(x\right)$ plot, where different kind of possible tail behaviours are represented. For instance, the Barenblatt profile $\mathcal{B}$, marked in dashed-grey  below, corresponds to the curve $\log\left(\mathcal{B}(x)\right)=-\frac{2+\beta-\gamma}{1-m}\log{|x|}+o(\log(|x|))$. The different lines represent other possible power-like tail behaviours; the thick line is the natural barrier for solutions to $\LL^1_{\gamma, +}(\RR^d)$ since it corresponds to the case $|x|^{-(d-\gamma)}$.

Let us begin our analysis. As we have already explained, every nonnegative solution to~\eqref{cauchy.problem} develops a minimal power-like tail, at least $|x|^{-\frac{\sigma}{1-m}}$, therefore in Figure~\ref{fig:picture.global.time1.fig}, there are no solutions below the dashed-grey line. Initial data in $\mathcal{X}$ develop exactly the Barenblatt's tail, see Theorem~\ref{thm.ghp}, hence, roughly speaking ``they live on the dashed-grey line'', the $\mathcal{X}$-curve. Things are different for initial data in $\mathcal{X}^c$. We only analyze power like behaviours at infinity. Roughly speaking, the Generalized GHP, Theorem~\ref{thm.GGHP}, says that all solutions must live between the   dotted and dashed line (recall that the red like correspond to the case $|x|^{-(d-\gamma)}$). More precisely, the Generalized GHP tells us that solutions with data decaying like $|x|^{-\alpha}$, with $(d-\gamma)<\alpha<\frac{\sigma}{1-m}$, will have the same decay $|x|^{-\alpha}$. Indeed, Theorem~\ref{thm.GGHP} tells us more: any initial datum  in $\mathcal{X}^c$ with a tail behaviour trapped between  two different lines  (maybe oscillating between two power-tails at infinity) produces a solution trapped between the same lines. For instance, if the datum is between the dotted and dashed lines, then the solution is trapped among those barriers for all times.

\subsection{A dynamical system interpretation}\label{sec.dyn.sys.intro}

The aim of this Section is to describe a global picture of the fine behaviour of the solutions to WFDE, in terms of convergence to equilibrium states of a (infinite dimensional) dynamical systems. It is convenient to pass to selfsimilar variables in order to make stationary the ``asymptotic'' Barenblatt solution.

\noindent\textbf{Self similar variables. Nonlinear Fokker-Plank equation. }
Let $u(t,x)$ be a solution to~\eqref{cauchy.problem} with initial data $u_0$, and consider $R(t)=R_\star(t+1)$. The \emph{self-similar change of variables}
\begin{equation}\label{change.fokker.planck}
v(\tau, y)= \frac{R(t)^{d-\gamma}}{\zeta^{d-\gamma}}\,u(t, x)\quad\mbox{where}\quad\tau=\frac{1}{\sigma}\log\frac{R(t)}{R(0)}\,,\,\,y=\frac{\zeta\,x}{R(t)}\,,
\end{equation}
transforms $u(t,x)$ into a solution to the following nonlinear \emph{Fokker-Planck} type equation
\begin{equation}\label{fokker.planck.equation}\tag{NLWFP}
\frac{\partial v}{\partial \tau} + |x|^\gamma\mathrm{div}\left(|x|^{-\beta}\,v\,\nabla v^{m-1}\right)=|x|^\gamma\mathrm{div}\left(|x|^{-\beta}\,v\,\nabla |x|^\sigma\right)\,,
\end{equation}
with initial data $v_0(y)=\frac{\zeta^{d-\gamma}}{R(0)^{d-\gamma}}u_0(\frac{\zeta\, x}{R(0)})$, with the same mass.
Also notice that among all the Barenblatt profiles $\mathfrak{B}(t+\tau,x; M)$, only the one with $\tau=1$ becomes stationary, and we call it Barenblatt  profile $\mathcal{B}_M(y)$: this is the unique attractor or the unique equilibrium (asymptotically stable).

In what follows we shall assume that solutions to~\eqref{cauchy.problem} with initial data $u_0 \in \LL^1_{\gamma, +}(\RR^d)$ will converge to the Barenblatt solution in the following sense (recall that the mass is preserved along the flow)
\[
\|u(t)-\mathfrak{B}(t;M)\|_{\LL^1_\gamma(\RR^d)}\rightarrow 0\quad\mbox{as}\quad t \rightarrow \infty\,.
\]
In \emph{self-similar} variables the previous result can be restated as
\[
\|v(\tau)-\mathcal{B}_M\|_{\LL^1_\gamma(\RR^d)}\rightarrow 0\quad\mbox{as}\quad \tau \rightarrow \infty\,.
\]

\noindent\textbf{Dynamical system approach: ``infinite dimensional phase plane analysis'' and the space $\mathcal{X}$. }It is well-known that the~\eqref{fokker.planck.equation} can be seen as the gradient flow of an Entropy functional, cf \cite{McCann1997,Otto2001}. It can be shown that solutions corresponding to nonnegative initial data will converge to a stationary solution with the same mass. To be more precise, let us define the $\omega$-limit of the~\eqref{fokker.planck.equation}  as the one dimensional manifold of the so-called Barenblatt solutions:
\begin{equation}\label{Barenblatt.Manifold}
\mathcal{M}_{\mathcal{B}}:=\left\{\mathcal{B}_M\;:\; M>0\right\}\,,
\end{equation}
and the distance
\[
\rd_1(f):=\inf_{\mathcal{B}_M\in \mathcal{M}_{\mathcal{B}}}{\|f-\mathcal{B}_M\|_{\LL^1_\gamma(\RR^d)}}.
\]
It has been proven in~\cite{Friedman1980,Vazquez2003} (for the case $\gamma=\beta=0$) and in~\cite{Bonforte2017a,Bonforte2017} (for the weighted case), that for any $u_0\in \LL^1_{\gamma,+}(\RR^d)$ there exists a unique $\mathcal{B}_{M_0}$ ($M_0$ being the mass of $u_0$) such that
\[
\rd_1(v(t))\le \|v(\tau)-\mathcal{B}_{M_0}\|_{\LL^1_\gamma(\RR^d)}\xrightarrow[t\to \infty]{}0\,.
\]
Hence solutions of~\eqref{fokker.planck.equation} can be seen as a continuous path with respect to the $\LL^1_{\gamma, +}$ topology,  that will eventually converge to a point of the manifold $\mathcal{M}_{\mathcal{B}}$. This fact can be rephrased as follows:  \textit{the basin of attraction of $\mathcal{M}_{\mathcal{B}}$ in the $\LL^1_{\gamma, +}$-topology is the whole space $\LL^1_{\gamma,+}$}.
\begin{figure}[h]
        \centering
        \includegraphics[trim= 1cm 0cm 0cm 0cm, clip=true,width=15cm]{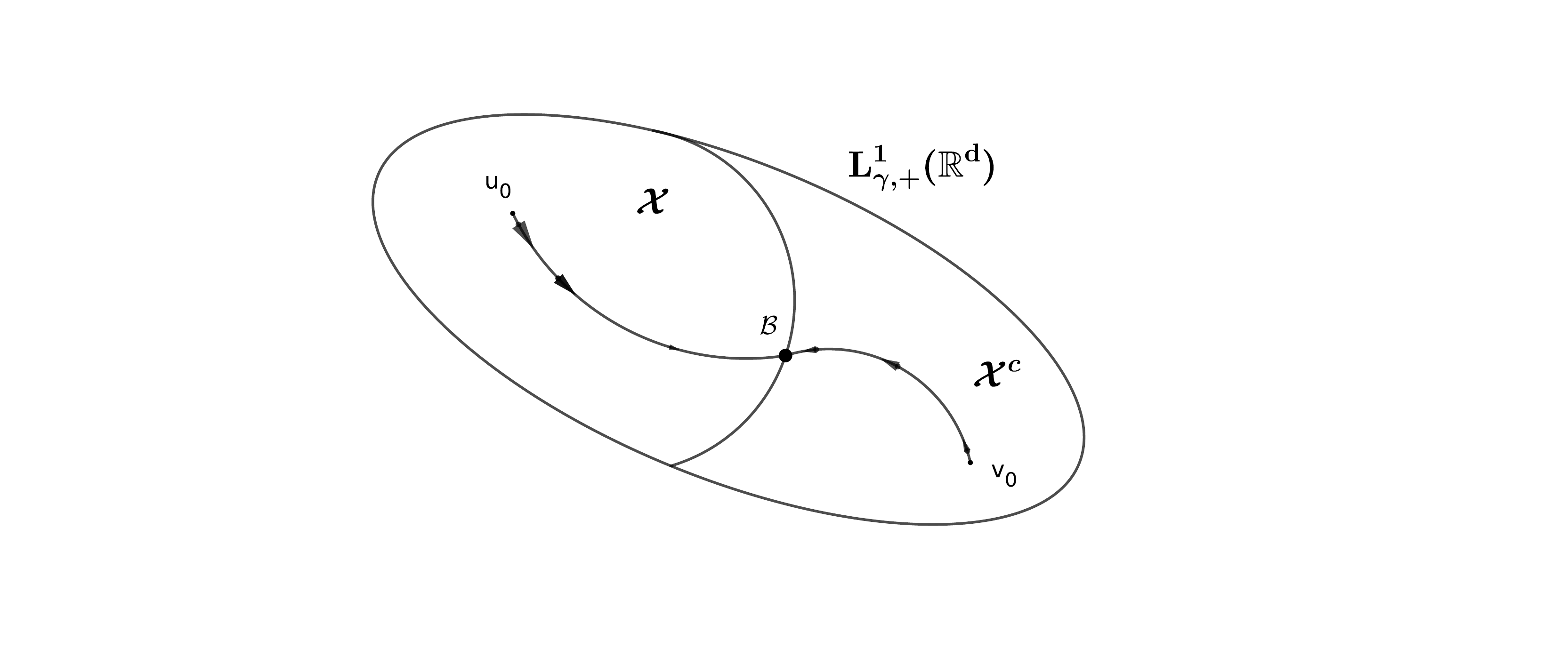}\label{fig:picture.L1.fig}
        \caption{ We represent two possible paths in $\LL^1_{\gamma, +}$. Since $\mathcal{X}$  and $\mathcal{X}^c$ are invariant sets for the flow there are no crossing lines between them. We notice that the manifold $\mathcal{M}_{\mathcal{B}}$ is contained in the topological boundary (with respect of the $\LL^1_{\gamma, +}$ topology) of $\mathcal{X}$, $\mathcal{M}_\mathcal{B}\subset \partial_{\LL^1_{\gamma, +}} \mathcal{X}$.}
\end{figure}\\

We can ask a similar question for a stronger convergence that allows to have a better asymptotic knowledge of the tails, \textit{the uniform converge in relative error} (UREC), properly measured by the following distance from $\mathcal{M}_{\mathcal{B}}$:
\[
\rd_\infty(f):=\inf_{\mathcal{B}_M\in \mathcal{M}_{\mathcal{B}}}\left\|\frac{f}{\mathcal{B}_M}-1\right\|_{\LL^\infty(\RR^d)}\,.
\]
\begin{figure}[h]
        \centering
        \includegraphics[trim= 1cm 0cm 2cm 0cm, clip=true,width=14cm]{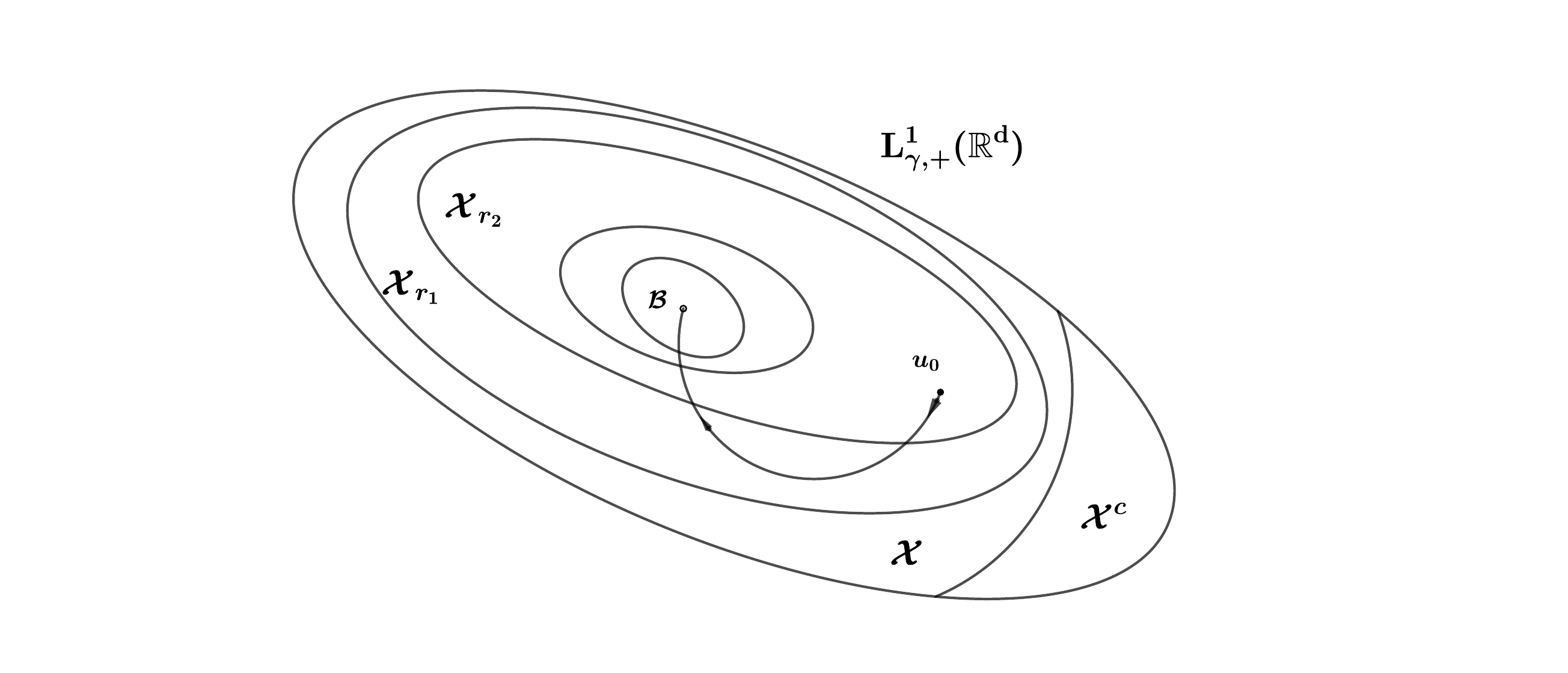}\label{fig:picture.Linfty.fig}
        \caption{Illustration of the stability of the sets $\mathcal{X}_r$: if the solution starts from one of those sets, it will forever stay in one of those sets. Indeed, if $u_0\in \mathcal{X}_r$ there exist is a maximal  $\mathcal{X}_{\overline{r}}$ which is invariant under the flow, i.e. s.t $u(t)\in \mathcal{X}_{\overline{r}}$ for all $t\ge 0$.}\vspace{-2mm}
\end{figure}
The above distance is induced by the  norm $\|f\|_{m, \gamma, \beta}:=\sup\limits_{x\in\RR^d}\left|f(1+|x|^{2+\beta-\gamma})^\frac{1}{1-m}\right|$, hence the topology is equivalent and we will refer to it as relative error topology.
As a consequence of $\LL^1_\gamma$ convergence, $v(\tau,x)\to\mathcal{B}_M(x)$ as $\tau\to \infty$ a.e., hence also  pointwise, i.e. $v(\tau,x)\mathcal{B}^{-1}_M(x)\to 1$ as $\tau\to \infty$.  However the \emph{uniform} convergence in relative error (UREC)  may fail. The main results of this paper solve this issue: \textit{$\mathcal{X}$ is the basin of attraction of the manifold $\mathcal{M}_{\mathcal{B}}$ in the relative error topology. }Notice that $\mathcal{X}$ is defined in terms of a practical -easy to check- condition on the initial datum, which a priori does not have any relation with the asymptotic behaviour. In what follow we explain our main results in terms of assumptions on the initial datum. Our main contribution in this direction is that we show that \emph{only} three things can happen: (i) $v_0\in\mathcal{X}^c$, (ii) $v_0\in\mathcal{X}$ and $\rd_\infty(v_0)<\infty$, (iii) $v_0\in\mathcal{X}$ and $\rd_\infty(v_0)=\infty$.  We analyze each case separately.\vspace{-2mm}
\begin{itemize}[leftmargin=*]\itemsep3pt \parskip3pt \parsep0pt
\item[(i)]\textit{If $v_0\in \mathcal{X}^c$. }Roughly speaking, in this case we show that if initial datum $v_0$ has a tail strictly above the Barenblatt one, then that ``fat tail'' is preserved in time. More precisely, Proposition~\ref{curve.X.prop} implies that
    \[
    v_0\in \mathcal{X}^c\Longrightarrow \left\|\frac{v(\tau)}{\mathcal{B}_M}-1\right\|_{\LL^\infty(\RR^d)}=+\infty\Longrightarrow v(\tau)\in \mathcal{X}^c\;\mbox{and}\; d_\infty(v(\tau))=\infty,
    \]
    for all $\tau>0$. In particular, since $d_\infty(v_0)=\infty$, it reveals that it is impossible to have bounds of the form $\rd_\infty(v(\tau))\le \rd_1(v(\tau))$, if we do not have it already (at least) for the initial datum. We can appreciate here a strong  difference along the flow between the $\LL^1_{\gamma, +}(\RR^d)$ and the $\|\cdot\|_{m, \gamma, \beta}$ topologies. On one hand, any $v_0\in \LL^1_{\gamma,+}$ is sent by the flow to a unique element of $\mathcal{M}_{\mathcal{B}}$ in the $d_1$ distance. On the other hand, this is not true in the $\rd_\infty$ distance, in which case, the flow stays always at infinite $\rd_\infty$-distance from $\mathcal{M}_{\mathcal{B}}$.

\item[(ii)]\textit{If $v_0\in \mathcal{X}$ and $\rd_\infty(v_0)<\infty$. }This is the stable case: if the initial datum is  close to the manifold $\mathcal{M}_{\mathcal{B}}$, then the flow will stay close to it and eventually $\rd_\infty$-converge to a unique element of $\mathcal{M}_{\mathcal{B}}$. More precisely, recall that $\mathcal{M}_{\mathcal{B}}\subset \mathcal{X}$ is the $\omega$-limit set, made of stationary solutions or equilibria of our dynamical system. The GHP of Theorem~\ref{thm.ghp} together with the UREC of  Theorem \ref{equivalent.condition}, imply
    \[
    v_0\in \mathcal{X} \Longrightarrow \left\|\frac{v(\tau)}{\mathcal{B}_M}-1\right\|_{\LL^\infty(\RR^d)}<\infty\qquad\mbox{for all $\tau>0$}\Longrightarrow v(\tau)\in \mathcal{X}\Longrightarrow d_\infty(v(\tau))\xrightarrow[\tau\to\infty]{}0\,.
    \]
    Indeed the GHP tells us a global stability result for the flow, since it can be rewritten as $\rd_\infty(v(t))\le F(\|v_0\|_{\mathcal{X}})$, for some locally bounded function $F$. A finer analysis is performed below.
\item[(iii)]\textit{If $v_0\in  \mathcal{X}$ and $\rd_\infty(v_0)=\infty$. }We show that even if the initial datum is at infinite distance  from the manifold $\mathcal{M}_{\mathcal{B}}$, but still in $\mathcal{X}$, the solution will eventually $\rd_\infty$-converge to $\mathcal{M}_{\mathcal{B}}$. Indeed, the GHP of Theorem~\ref{thm.ghp} only needs the assumption $v_0\in \mathcal{X}$, regardless of $\rd_\infty(v_0)=\infty$. Hence the same argument as case (ii) applies.
\end{itemize}
\noindent\textbf{Finer analysis in $\mathcal{X}$. }It is possible to show that
\[
\mathcal{X}=\{d_\infty=\infty\}\cup\bigcup_{r>0} \mathcal{X}^\infty_r=\{d_\infty=\infty\}\cup\bigcup_{r>0} \left\{f\in \mathcal{X}\;:\; d_\infty(f)<r\right\}
\]
The GHP of Theorem~\ref{thm.ghp}, reveals an important stability property of the fast diffusion flows: for any $v_0\in \mathcal{X}$ there exists $r_0,\tau_0>0$ s.t. $d_\infty(v(\tau))<r_0$ for all $\tau\ge \tau_0$, hence the flow never exit from a certain $\mathcal{X}^\infty_{r_0}$. Indeed, we show more: $d_\infty(v(\tau))\to 0$ as $\tau\to \infty$, which means that the flow always leaves the manifolds $d_\infty(v(\tau))=r$ (level sets of distance from $\mathcal{M}_{\mathcal{B}}$) to enter one at a lower level, say $d_\infty(v(\tau))=r-\varepsilon$.

This can be summarized as follows: we show that the solution map sends immediately $\mathcal{X}$ in a more regular subspace $\mathcal{X}^\infty_{r_0}$
\[
T_\tau: \mathcal{X}\to \bigcup_{r>0} \mathcal{X}^\infty_r\qquad\mbox{indeed there exists }r_0>0\,:\qquad
T_\tau: \mathcal{X}\to  \mathcal{X}^\infty_{r_0}\,.
\]
On one hand, in the relative error topology we have a dichotomy: $\lim\limits_{\tau\to\infty}T_\tau(\mathcal{X}\setminus\{0\})=\mathcal{M}_{\mathcal{B}}$ and $d_\infty(T_\tau(\mathcal{X}^c))=\infty$ for all $\tau>0$. On the other hand, in the $d_1$-topology we always have $\lim\limits_{\tau\to\infty}T_\tau(\LL^1_{\gamma,+}(\RR^d)\setminus\{0\})=\mathcal{M}_{\mathcal{B}}$\\

\subsection{Others ranges of $ m $ and generalizations} In the study of~\eqref{cauchy.problem} a dramatic change occurs when we consider the exponent $m\le m_c$ or $m\ge1$. When $m=1$, it is known that solutions to the Cauchy problem for the classical heat equation develop eventually Gaussian tails whenever the initial data is compactly supported. Surprisingly enough, no sharp results are known it this contest, to the best of our knowledge. Uniform convergence in relative error is (in general) false, see~\cite{Vazquez2017}, and several tail behaviours are indeed possible, see~\cite{Herraiz1999}.   The slow diffusion case $m>1$ has been widely investigated, see the monograph~\cite{Vazquez2007}, but due the finite speed of propagation, fine results about convergence in relative error are still missing. Let us now comment on the very fast diffusion case $m\le m_c$, where several new difficulties arise such as, for instance, loss of mass and extinction in finite time. In the case $m=m_c$ conservation of mass still holds, nevertheless the asymptotic behaviour becomes quite involved, see~\cite{Kalashnikov1987,Galaktionov2000}. Below the critical exponent $m_c$  fewer results are known and some considerations are in order.  The fundamental solution does not exist anymore,~\cite{Vazquez2006}. A large class of solutions vanish in finite time with different possible behaviour near the extinction time. The vanishing profile of a suitable class of initial data is represented by the so called \emph{Pseudo}-Barenblatt solutions, see~\cite{Blanchet2009}. The extinction behaviour of bounded and integrable solutions for $0<m<m_c$ is only known in the radial case, \cite{Galaktionov1997,Vazquez2006}. In the Yamabe flow case, $m=\frac{d-2}{d+2}$, finer results are known, see~\cite{Pino2001,Daskalopoulos2008,Daskalopoulos2018,Daskalopoulos2013}.
The situation is completely different the so-called ultra fast diffusion case, $m\le0$. Indeed, for the Cauchy problem nonnegative $\LL^1$ data do not produce solutions~\cite{Rodriguez1990,Vazquez1992}. As a consequence, the are no solutions for the homogeneous Dirichlet problem, while a special class of solutions can be found for the Neumann problem, \cite{Schimperna2012,Schimperna2016}. Sharp existence and non existence conditions, for the Cauchy problem, has been given in \cite{Daskalopoulos1994,Daskalopoulos1995,Daskalopoulos1997} and \cite{Rodriguez1997}.

\noindent\textit{Possible generalizations. }
The Global Harnack Principle (Theorem~\ref{thm.ghp}) can be generalized to solutions to equations of the form
\begin{equation}\label{eq.with.coefficients}
u_t = \mathrm{div}\left(A(t,x, u, \nabla u^m)\right)\,,
\end{equation}
where $A(t,x, u, \eta)$ satisfies suitable structural conditions, as those in~\cite{Ragnedda2013,Recalde2016}.  Even if the fundamental solution exists, see~\cite{Ragnedda2013,Recalde2016}, it is not clear whether or not (and in which sense) it represents the large time behaviour of nonnegative, integrable solutions to~\eqref{eq.with.coefficients}.

\medskip
\noindent\textbf{Organization of the paper. }In Section~\ref{section:l1} we collect some results that hold for all nonnegative integrable solutions: we find universal lower bounds in terms of Barenblatt profiles, which allow to identify the minimal tail of all nonnegative solutions, Corollary \ref{herrero.pierre.limit.corollary}. We prove that the answer to $Q_2$ can be yes for all $u_0\in\LL^1_+(\RR^d)$, but only on suitable parabolic cones, which represent the optimal domain of validity for such results.
In Section~\ref{section:GHP} we analyze the behaviour of solutions whose initial data are in $\mathcal{X}$ and provide a positive answer to both $Q_1$ and $Q_2$. We prove the upper part of inequality~\eqref{ghp.inequality} (Theorem~\ref{thm.ghp}) and Theorem~\ref{equivalent.condition}.  As a consequence, we obtain rates of convergence to the Barenblatt profile in several norms.
In Section~\ref{Sec4} we construct sub/super solutions with the anomalous tail behaviour analyzed above. We also show, by means of counterexamples, that the power-like rates obtained in Section~\ref{section:GHP} are not possible for data outside $\mathcal{X}$.
In Section~\ref{sec:on.X} we show the equivalence between~\eqref{tail.condition} and~\eqref{tail.condition.1} and we give example of function in $\mathcal{X}$ which do not satisfies the pointwise condition $u_0(x)\lesssim |x|^{-\frac{\sigma}{1-m}}$. We also give more details about the natural topology of $\mathcal{X}$ and analyze stability properties of the WFDE-flow as curve in $\mathcal{X}$.
\medskip

\noindent\textbf{Notations. }%
We will systematically use $\infty$ to indicate $+\infty$. We will use the following notations throughout the paper:  $a\wedge b=:\min\{a,b\}$, $a\vee b:=\max\{a,b\}$ and $a\asymp b$ means that there exist constants $c_1, c_2>0$ such that $c_1 a \le b \le c_2 a$; similarly we write $a\lesssim b$ whenever there exists $c>0$ such that $a\le c\,b$. Also, give $B \subset \RR^d$ we define $\chi_{B}$ as the characteristic function of $B$, namely $\chi_{B}(x)=1$ if $x\in B$, while $\chi_{B}(x)=0$ if $x\not\in B$.
\section{Initial Data in $\LL^1_{\gamma, +}(\RR^d)$}\label{section:l1}

In this Section we show the results that hold for all data in $\LL^1_{\gamma, +}(\RR^d)$, namely, we show that the lower-part of the GHP estimates hold true (Theorem~\ref{lower.estimate.theorem}) for (just) locally integrable data: this allows to measure the (infinite) speed of propagation as ``fatness of the tails''. On the other hand, on the whole space it is not possible to match the lower bounds with similar upper bounds for all initial data in $\LL^1_{\gamma, +}(\RR^d)$: we will provide explicit counterexamples and improved lower bounds in Section~\ref{Sec4}. This latter phenomenon, an anomalous tail behaviour, can only happen if we miss a control the tail of the initial datum: we will show that the sharp tail-condition is encoded in the space $\mathcal{X}$ thoroughly analyzed in Section~\ref{section:GHP}. As a consequence of the estimates of this section, we show also uniform convergence in relative error towards equilibrium on compact sets and even on parabolic cones, see Theorem~\ref{local.convergence.relative.error.thm}. All of these results are sharp, as shown in Section~\ref{Sec4} by means of suitable counterexamples.

\subsection{A universal global lower bound: measuring the speed of propagation}
We now state the main result of this section, which holds for nonnegative initial data which are merely locally integrable. We recall here a useful quantity, $t_*$ that will appear frequently throughout this section:
\begin{equation}\label{Minimal.life.time}
     t_{\ast}= t_{\ast}(u_0, R)= \kappa_*\, \|u_0\|_{\LL^1_{\gamma}(B_{R}(0))}^{1-m}\,R^{\frac{1}{\vartheta}} \,.
\end{equation}
where $\kappa_*>0$ depends on $d, m, \gamma, \beta$.   We are now ready to state the main result of this section.

\begin{thm}\label{lower.estimate.theorem}
Let $u$ be a solution to \eqref{cauchy.problem} with initial data $0 \leq u_0 \in \LL^1_{\gamma, \mathrm{loc}}(\RR^d)$ and let $t_0, R_0>0$ be such that $ \|u_0\|_{\LL^1_{\gamma}(B_{R_0}(0))} > 0$. Then there exists $\underline{\tau} > 0$  and $\underline{M}>0$ such that
\begin{equation}\label{lower.estimate.inequality}
u(t, x) \geq \mathfrak{B}(t-\underline{\tau},x; \underline{M}),\qquad\mbox{for all $x \in \RR^d$ and $t \geq t_0$}.
\end{equation}
where
\begin{equation}\label{lower.estimate.inequality.constants}
\underline{\tau}=  \frac{1}{2} \, \left(t_* \wedge t_0\right)\qquad\mbox{and}\qquad \underline{M}=b \, \|u_0\|_{\LL^1_\gamma(B_{R_0}(0))} \, \left(1 \wedge \frac{t_0}{t_*}\right)^\frac{1}{1-m}.
\end{equation}
The constant $ b>0$ depends only on $d, m, \gamma, \beta$ and has an explicit expression given in the proofs, while $t_*$ is as in \eqref{Minimal.life.time}.
\end{thm}
\noindent\textbf{Measuring the speed of propagation. }The above Theorem partially answers $Q_1$ and reveals a remarkable property of solutions to WFDE: the positivity spreads immediately for every nonnegative initial datum, showing infinite speed of propagation.
A delicate issue is how to discriminate in a quantitative way among two infinite speed of propagation. Our Theorem shows that we can put a (delayed) fundamental solution as a lower barrier for any data: this is how the WFDE immediately creates a fat tail (inverse power), which is clearly bigger than the ``standard Gaussian tail'' (exponentially decaying) created by the linear heat equation. This can be expressed as follows:
\begin{cor}[Minimal tails]\label{herrero.pierre.limit.corollary}
Under the assumption of Theorem~\ref{lower.estimate.theorem} we have that for any $t>0$
\begin{equation}\label{Herrero.Pierre.Limit}
\liminf_{|x|\rightarrow \infty}\, u(t,x)\,|x|^\frac{\sigma}{1-m} \ge b_1\, t^\frac{1}{1-m}
\end{equation}
The constant $b_1$ depends only on $m, d, \gamma, \beta$ and is achieved by the Barenblatt solutions.
\end{cor}
We will often call $|x|^{-\sigma/(1-m)}$ \textit{a minimal tail or a Barenblatt tail }. Finding matching upper bounds is simply not possible in such generality, we will need to ask the tail condition~\eqref{tail.condition} on $u_0$.

\medskip
\noindent {\bf Proof of Theorem \ref{lower.estimate.theorem}: }Let us first state an inequality proven in \cite[Theorem 1.4]{Bonforte2019a}, a sharp local lower bound (half-Harnack inequality), essential to this proof. We do not use here Aleksandrov Principle, as in \cite{Bonforte2006}, nor other moving planes argument. Under the running assumption we have that
\begin{equation}\label{local.lower.estimate.minimal}
\inf_{x \in B_{2R}(0)} u(t_\star,x) \ge \kb_1 \frac{\|u_0\|_{\LL^1_{\gamma}(B_{R}(0))}}{R^{d-\gamma}},
\end{equation}
where $\kb_1$ depends only on $d, m, \gamma, \beta$, and has an explicit expression given in \cite{Bonforte2019a}.
Let us now explain the strategy of the proof. The quantities $\underline{\tau}$ and $\underline{M}$ take different forms depending wether or not $t_*\le t_0$. We will assume first  that $t_0 \geq t_*$, then  we will discuss the case $0< t_0 < t_*$ at the end of the proof.

Let  $M_{R_0}=\|u_0\|_{\LL^1_{\gamma}(B_{R_0}(0))}$, $\underline{\tau}=at_*$ and $\underline{M}=bM_{R_0}$ where $a \in (0,1)$ and $b>0$ will be explicitly chosen later. Without loss of generality, we prove inequality \eqref{lower.estimate.inequality} only at $t=t_*$, namely
\begin{equation}\label{lower.ghp.1}
u(t_*, x) \geq \mathfrak{B}((1-a)t_*,x; \underline{M})\,.
\end{equation}
Once proven at $t=t_*$, the case $t \geq t_*$ will follow by comparison. To prove \eqref{lower.ghp.1} we need to determine the values of $a, b$. We need to separate two cases, namely inside a ball  and outside a ball, obtaining different conditions on $a,b$, respectively conditions \eqref{condition.lower.bounds.1} and \eqref{condition.lower.bounds.2}. Finally we check the compatibility of such conditions and choose $a,b$ explicitly as in \eqref{a.b.chosen}.

\noindent\textit{Condition on $a,b$ inside a ball. }We want to find condition on $a, b$ such that the following  inequality holds:
\begin{equation}\label{inequality.inside.ball}
\kb_1 \frac{M_{R_0}}{{R_0}^{d-\gamma}} \geq \frac{b^{\sigma \vartheta}M_{R_0}}{b_0^\frac{1}{1-m}(1-a)^{(d-\gamma)\vartheta}\kappa_*^{(d-\gamma)\vartheta}{R_0}^{d-\gamma}}= \sup_{x \in B_{2{R_0}}(0)} \mathfrak{B}(t_*-\underline{\tau},x; \underline{M})\,,
\end{equation}
where $\kb_1$ is as in \eqref{local.lower.estimate.minimal}. It is easily seen that the former is implied by the following condition on $a$ and $b$:
\begin{equation}\label{condition.lower.bounds.1}
b^{\sigma \vartheta} \leq \kappa_*^{(d-\gamma)\vartheta} \kappa_1 b_0^{\frac{1}{1-m}} (1-a)^{(d-\gamma)\vartheta}\,.
\end{equation}
Note that by inequality \eqref{local.lower.estimate.minimal} the first term in \eqref{inequality.inside.ball} is bounded above by $\inf_{x \in B_{2R_0}}u(t_*, x)$, therefore we obtain that
\begin{equation*}
\inf_{x \in B_{2R_0}}u(t_*, x) \geq     \sup_{x \in B_{2{R_0}}(0)} \mathfrak{B}(t_*-\underline{\tau},x; \underline{M})\,,
\end{equation*}
inequality \eqref{lower.ghp.1} is then proved for any $|x| \leq 2 R_0$.

\noindent\textit{Condition on $a,b$ outside a ball.} We want to find suitable conditions on $a, b$ such that \eqref{lower.ghp.1} holds in the outer region $|x|> R_0$. Such an inequality will be deduced by applying the comparison on the parabolic boundary of $Q=\left(\underline{\tau}, t_*\right)\times B_{R_0}^c(0)$, namely $\partial_p Q=\{\left\{\underline{\tau}\right\}\times B_{R_0}^c(0)\} \bigcup \{\left(\underline{\tau}, t_*\right)\times \left\{x \in \RR^d: |x|=R_0\right\}\}$, see for instance \cite[Lemma 3.4]{Herrero1985}.

It is clear that $u(\underline{\tau}, x) \geq \mathfrak{B}(0,x; \underline{M})=\delta_0(x)$, for any $|x| \geq R_0$, hence we just need to prove that
\begin{equation}\label{inequality.prove.2}
u(t, x) \geq \mathfrak{B}(t-\underline{\tau},x; \underline{M})\quad\mbox{for any}\quad|x|=R_0,\, t \in \left(\underline{\tau}, t_*\right).
\end{equation}
A straightforward computation shows that, under the running assumption,  for $|x|=R_0$ we have that
\begin{equation}
\sup_{t\ge\tau}\mathfrak{B}(t-\underline{\tau},x; \underline{M}) = \left(\frac{b_1}{b_0\, \vartheta}\right)^\frac{1}{\sigma\vartheta(1-m)}\,\frac{\left[(d-\gamma)(1-m)\right]^\frac{d-\gamma}{\sigma}}{\left[\kappa_\star\,\sigma\right]^\frac{1}{1-m}}\,\left(\frac{t_*}{b_0\,R_0^\sigma}\right)^\frac{1}{1-m}\, b\,.
\end{equation}
The following inequality
\begin{equation}\label{inequality.impose}
\kb \left(\frac{a t_*}{{R_0}^\sigma} \right)^{\frac{1}{1-m}}\geq \left(\frac{b_1}{b_0\, \vartheta}\right)^\frac{1}{\sigma\vartheta(1-m)}\,\frac{\left[(d-\gamma)(1-m)\right]^\frac{d-\gamma}{\sigma}}{\left[\kappa_\star\,\sigma\right]^\frac{1}{1-m}}\,\left(\frac{t_*}{b_0\,R_0^\sigma}\right)^\frac{1}{1-m}\, b\,,
\end{equation}
implies that inequality \eqref{inequality.prove.2} holds, indeed for any $|x|=R_0$ and $t \in \left(\underline{\tau}, t_*\right)$ we have that
\begin{equation*}
u(t, x) \geq \inf_{\substack{t \in \left(at_*, t_*\right),\\ x \in B_{2R_0}(0)}} u(t, x) \geq \kb \left(\frac{a t_*}{{R_0}^\sigma} \right)^{\frac{1}{1-m}}\geq\sup_{t\ge\tau}\mathfrak{B}(t-\underline{\tau},x; \underline{M}) \geq \mathfrak{B}(t-\underline{\tau},x; \underline{M})\,.
\end{equation*}
It is easy to show that inequality \eqref{inequality.impose} is equivalent to the following one
\begin{equation}\label{condition.lower.bounds.2}
b^{\sigma\vartheta} \leq b_0^{\frac{1}{1-m}}\, a^\frac{\sigma\vartheta}{1-m}\,\kb^{\sigma\vartheta} \left(\frac{\vartheta\,b_0^{\sigma\vartheta}}{b_1}\right)^\frac{1}{1-m}\,\frac{\left[\kappa_\star\,\sigma\right]^\frac{\sigma\vartheta}{1-m}}{\left[(d-\gamma)(1-m)\right]^{(d-\gamma)\vartheta}}
\end{equation}
which is the condition we were looking for.
\noindent\textit{Compatibility of condition \eqref{condition.lower.bounds.1} and \eqref{condition.lower.bounds.2}.} Both the conditions are satisfied by the following choice
\begin{equation}\label{a.b.chosen}
a=\frac{1}{2}\qquad\mbox{and}\qquad b^{\sigma\vartheta} = b_0^\frac{1}{1-m}\, \left[ \left(\frac{\vartheta\,b_0^{\sigma\vartheta}}{2^{\sigma\vartheta}\,b_1}\right)^\frac{1}{1-m}\,\frac{\kb^{\sigma\vartheta} \,\left[\kappa_\star\,\sigma\right]^\frac{\sigma\vartheta}{1-m}}{\left[(d-\gamma)(1-m)\right]^{(d-\gamma)\vartheta}} \wedge \frac{\kappa_*^{(d-\gamma)\vartheta} \kappa_1}{2^{(d-\gamma)\vartheta}}\right]\,.
\end{equation}
This concludes the proof of \eqref{lower.ghp.1} in the case when $t_0\geq t_*$. It only remains to analyze the case when $t_0 < t_*$.

\noindent\textit{Case $0<t_0 < t_*$. }Without loss of generality, we only need to prove inequality \eqref{lower.estimate.inequality} at time $t=t_0$, the full result will then follow by comparison. Recall the Benilan-Crandall-type estimate, \cite{Benilan1981},
\begin{equation}\label{benilan.crandall}
u(t_0, x) \geq u(t_*, x) \left(\frac{t_0}{t_*}\right)^\frac{1}{1-m}\,, \qquad\mbox{for all $0 < t_0 < t_* $.}
\end{equation}
Now we recall that inequality \eqref{lower.ghp.1}  holds under the choices of $a,b$ as in \eqref{a.b.chosen}.
Using inequality \eqref{lower.ghp.1} and inequality \eqref{benilan.crandall}  we get
\begin{equation*}\begin{split}
u(t_0, x) \geq u(t_*, x) \left(\frac{t_0}{t_*}\right)^\frac{1}{1-m} &\geq \frac{2^{-\frac{1}{1-m}}t_*^\frac{1}{1-m}}{\left[b_0 \frac{2^{-\sigma\vartheta} t_*^{\sigma\vartheta}}{\underline{M}^{\sigma\vartheta(1-m)}}+b_1 |x|^\sigma\right]^\frac{1}{1-m}} \, \left(\frac{t_0}{t_*}\right)^\frac{1}{1-m} \\
&=\frac{2^{-\frac{1}{1-m}}t_0^\frac{1}{1-m}}{\left[b_0 \frac{2^{-\sigma\vartheta} t_0^{\sigma\vartheta}}{\underline{M}^{\sigma\vartheta(1-m)} \left[\left(\frac{t_0}{t_*}\right)^\frac{1}{1-m}\right]^{\sigma\vartheta}}+b_1 |x|^\sigma\right]^\frac{1}{1-m}} = \mathfrak{B}\left(t_0-\frac{t_0}{2},x; \left(\frac{t_0}{t_*}\right)^\frac{1}{1-m}\underline{M}\right)\,.
\end{split}\end{equation*}
Recalling that in this case $\underline{\tau}=t_0/2$, the proof is concluded.\qed

We can now give  the proof of Corollary~\ref{herrero.pierre.limit.corollary}.

\noindent {\bf Proof of Corollary~\ref{herrero.pierre.limit.corollary}.} Let $R_0$ be such that $\|u_0\|_{\LL^1_\gamma(B_{R_0}(0))}>0$,  $t>0$, and $0<\varepsilon <t$. By applying Theorem~\ref{lower.estimate.theorem} at time $t_0=t-\varepsilon$ and radius $R_0$ we get the following inequality
\[
u(t, x) \ge \mathfrak{B}(t-\underline{\tau}, x; \underline{M})\,.
\]
As a consequence we obtain
\[
\liminf_{|x|\rightarrow \infty}\, u(t,x)\,|x|^\frac{\sigma}{1-m} \ge b_1\, \left[t-\frac{1}{2}\left(t_*\wedge t_0\right)\right]^\frac{1}{1-m},
\]
from which~\eqref{Herrero.Pierre.Limit} follows just by taking the limit for $\varepsilon\rightarrow t$. Notice that in such a limit $t_0\rightarrow 0$.\qed
\subsection{Harnack Inequality in Parabolic Cones}
We have shown in~\cite{Bonforte2019a} that nonnegative local solutions to WFDE satisfy Harnack inequalities of various kind:
an \emph{elliptic} form (in which the supremum and the infimum are taken at the same time), a \emph{forward} in time (the supremum is taken at a smaller time than the infimum) and a \emph{backward} in time (the supremum is taken at a bigger time than the infimum). We remark that for solutions to the \emph{heat equations} in general only the \emph{forward} Harnack inequality holds. Here we shall prove an \emph{elliptic} form of a Harnack inequality on conical space-time domains of the form
\begin{equation}\label{KKK}
K(t)=K_M(t) =\{|x| \le\,  t^\vartheta\,M^{(m-1)\vartheta}\}\,.
\end{equation}
for some fixed $M >0$. We will call these sets ``Parabolic Cones'', with a slight abuse of language, indeed for $\vartheta=1$, $K(t)$ are really cones in space time domains of the form $\RR_+\times \RR^N$.  A similar inequality on balls has been proven in~\cite[Theorem 1.4]{Bonforte2006}.
\begin{thm}[Harnack inequality in Parabolic Cones]\label{thm.Harnack.Parabolic.Cones}
Let $u$ be a solution to \eqref{cauchy.problem} with initial data $0 \leq u_0 \in \LL^1_\gamma(\RR^d)$. Let $M=\|u_0\|_{\LL^1_\gamma(\RR^d)}$ and $R_0>0$ be such that $\|u_0\|_{\LL^1_{\gamma}(B_{R_0}(0))} = M/2$, and let $t_*=\kappa_*\,R_0^\frac{1}{\vartheta}\,\left(M/2\right)^\frac{1}{1-m}$.  Then, there exists a positive constant $\mathcal H$ such that
\begin{equation}\label{elliptic.harnack.inequality}
\sup_{x \in K(t)} \frac{u(t,x)}{\mathfrak{B}(t,x; M)} \le \mathcal H \,\inf_{x \in K(t)} \frac{u(t,x)}{\mathfrak{B}(t,x; M)}\,,
\qquad\mbox{for any $t\ge 3\,t_*.$}
\end{equation}
where the constant $\mathcal H$ depends only on $m, d, \gamma, \beta$ and $K(t)$ depends on $M$ as in \eqref{KKK}
\end{thm}
\noindent\textbf{Proof. }In the proof we will make use the Smoothing Effect for solutions to~\eqref{cauchy.problem}, namely the following inequality which hold for any $t>0$
\begin{equation}\label{global.upper.estimate}
\|u(t)\|_{\LL^{\infty}(\RR^d)}\leq \frac{\ka_1}{t^{(d-\gamma)\vartheta}} \left[\int_{\RR^d}{|u_0(y)| \ |y|^{-\gamma} \ \dy} \right]^{\sigma  \vartheta},
\end{equation}
where $\vartheta$ and $\sigma$ are defined in \eqref{parameters.sigma.theta}. The constant $\ka_1$ has an explicit form, cf.~\cite{Bonforte2019a}, it depends only on $d, m,\gamma$ and $\beta$ and it is the same constant $\ka_1$ appearing in formula~\eqref{local.upper.estimte}. Inequality~\eqref{global.upper.estimate} has been obtained in~\cite{Bonforte2019a} and can be easily deduced by taking $x_0=0$ and letting $R_0\rightarrow \infty$ in estimate~\eqref{local.upper.estimte}.

Let us now begin the proof. By applying Theorem~\ref{lower.estimate.theorem} we deduce that $u(t,x)\ge \mathfrak{B}(t-\underline{\tau},x; \underline{M})$ with $\underline{\tau}=\frac{t_*}{2}=   \frac{\kappa_*}{2}\,\,R_0^\frac{1}{\vartheta}\,\left(\frac{M}{2}\right)^\frac{1}{1-m}$ and $\underline{M}=b\,M/2$ where $b$ is as in \eqref{condition.lower.bounds.2}. In view of the Smoothing Effects~\eqref{global.upper.estimate} and of inequality~\eqref{lower.estimate.inequality}, it is enough to prove that there exists $\mathcal H$ such that
\[
\ka_1 \,\left(b_0+b_1\right)^\frac{1}{1-m}\le \mathcal H \,b_0^\frac{1}{1-m} \frac{t^{(d-\gamma)\vartheta}}{M^{\sigma\vartheta}}\, \inf_{x \in K(t)} \mathfrak{B}(t-\underline{\tau},x; \underline{M})\,
\]
This amounts to prove that the following quotient is uniformly bounded by $\mathcal{H}$ for $t \ge 3\,t_*$:
\begin{equation*}
\quad\ka_1\left(1+\frac{b_1}{b_0}\right)^\frac{1}{1-m} \frac{M^{\sigma\vartheta}}{t^{(d-\gamma)\vartheta}} \frac{\left[\frac{b_0(t-\underline{\tau})^{\sigma\vartheta}}{\underline{M}^{2\vartheta(1-m)}}+ \frac{b_1 t^{\sigma\vartheta}}{M^{\sigma\vartheta(1-m)}}\right]^\frac{1}{1-m}}{\left(t-\underline{\tau}\right)^{\frac{1}{1-m}}}\le \mathcal{H}\,.
\end{equation*}
Since $\underline{\tau}=t_*/2$ we easily conclude that $\mathcal H$ can be taken as
\[
\mathcal H = \ka_1\, \left(1+\frac{b_1}{b_0}\right)^\frac{1}{1-m}\,5^\frac{1}{1-m}\, \left[b_0 \, \left(\frac{2}{b}\right)^{\sigma\vartheta} + b_1\right]^\frac{1}{1-m}\mbox{\,.\qed}
\]

\subsection{Uniform Convergence in  relative error in Parabolic Cones}\label{ssec:ucre.parabolic.cones}
In this section we will prove that solutions to~\eqref{cauchy.problem} with initial data $u_0 \in \LL^1_{\gamma, +}(\RR^d)$ converge to the Barenblatt profile $\mathfrak{B}(t, x; M)$ in \emph{relative error} uniformly in parabolic cones, and as a consequence uniformly on compact subsets of $\RR^d$. To obtain such a result we will use the convergence to the Barenblatt profile in $\LL^1_\gamma(\RR^d)$, namely
\begin{equation}\label{convergence.bareblatt.l1.fde}
\|u(t)-\mathfrak{B}(t;M)\|_{\LL^1_\gamma(\RR^d)}\rightarrow 0\quad\mbox{as}\quad t \rightarrow \infty\,,
\end{equation}
or equivalently, in \emph{self-similar variables}
\begin{equation}\label{convergence.bareblatt.l1.nfp}
\|v(\tau)-\mathcal{B}_M\|_{\LL^1_\gamma(\RR^d)}\rightarrow 0\quad\mbox{as}\quad \tau \rightarrow \infty\,,
\end{equation}
where  $v(\tau, y)$ is defined in~\eqref{change.fokker.planck} and it is a solution to~\eqref{fokker.planck.equation}. The proof of~\eqref{convergence.bareblatt.l1.fde} can be done by a straightforward adaptation to our setting of the so called ``4 step method'', carefully explained in~\cite[Theorem 1.1]{Vazquez2003}. We leave the details to the interested reader, just noticing that the proof contained in~\cite{Vazquez2003} deals with the case $m>1$, and uses compactly supported initial data, hence compactly supported solutions (when $m>1$ there is finite speed of propagation). In the present setting, the very same proof works, just by replacing the compactly supported solutions by the ones which satisfy the GHP.

\begin{thm}[Uniform Convergence in \emph{relative error} on parabolic cones]\label{local.convergence.relative.error.thm}
Let $m \in \left(m_c, 1\right)$ and let $u$ be a solution to \eqref{cauchy.problem} with initial data $0 \leq u_0 \in \LL^1_\gamma(\RR^d)$ and let $M=\|u_0\|_{\LL^1_\gamma(\RR^d)}$. Then for any  $\Upsilon>0$ we have that
\begin{equation}\label{local.convergence.relative.error}
\lim_{t\to \infty}\sup_{x \in \{|x| \le \Upsilon\, t^\vartheta\}} \left|\frac{u(t,x)-\mathfrak{B}(t, x; M)}{\mathfrak{B}(t, x; M)} \right|=0\,.
\end{equation}
\end{thm}
\noindent\textbf{Remark. }As an easy corollary of the previous Theorem, we obtain that
\[
\left\|\frac{u(t,x)-\mathfrak{B}(t, x; M)}{\mathfrak{B}(t, x; M)}\right\|_{\LL^\infty(K)}\xrightarrow[t\to\infty]{} 0\qquad\mbox{for any compact set $K\subset \RR^d$.}
\]
This follows from inequality \eqref{local.convergence.relative.error} just by observing that $K\subset \{|x| \le \Upsilon\, t^\vartheta\}$ for some $t_0>0$.

\noindent\textbf{Proof.} We split the proof into several steps. First we prove an uniform pointwise estimate on the solution $u(t,x)$ in domains of the form $\{|x| \le C R(t)\}$, where $R(t)$ is as in~\eqref{change.fokker.planck} and $C>0$. We remark that for any $t>0$ we have that $\{|x| \le C t^\vartheta\} \subset \{|x| \le C R(t)\}$. As a second step we will rescale $u(t,x)$ to self-similar variables (we recall that domains of type $\{|x| \le C R(t)\}$ are transformed into $B_\rho(0)$, where $\rho=\zeta\,C$) and, using the estimates obtained before, we  estimate $\lfloor v(\tau, \cdot) - \mathfrak B_M\rfloor_{C^\nu(B_{3r})}$ uniformly in time. Finally, by applying a clever interpolation, Lemma \ref{interpolation.lemma} we  prove that $\|v(\tau, \cdot)-\mathfrak B_M\|_{\LL^\infty(B_{r})}\rightarrow 0$ as $\tau \rightarrow \infty$, and finally~\eqref{local.convergence.relative.error} follows.\\

\noindent\textit{Uniform estimate on $u(t,x)$ in $\{|x| \le 3\,\Upsilon\, R(t)\}$.} Let $\rho>0$ be such that $\int_{B_\rho}u_0(x)|x|^{-\gamma}\dx=\frac{M}{2}$ and define $t_\star=\kappa_\star\,\rho^\frac{1}{\vartheta}\,\left(\frac{M}{2}\right)^{1-m}$ where $\kappa_\star$ is as in~\eqref{Minimal.life.time}. By applying Theorem~\ref{lower.estimate.theorem} and the global smoothing effect, inequality~\eqref{global.upper.estimate}, we obtain that for any $t\ge t_\star$
\[
\mathfrak{B}(t-\underline{t},x; \underline{M}) \le u(t,x) \le \ka_1 \frac{M^{\sigma\vartheta}}{t^{(d-\gamma)\vartheta}}\,,
\]
where $\underline{t}=\frac{t_\star}{2}$ and $\underline{M}=\frac{b}{2}\,M$. By the above inequality, we can deduce the following matching lower bound, by means of straightforward estimates relying on the explicit expression of $\mathfrak{B}$: there exists a constant $\kb_1>0$ which depends on $d, m, \gamma, \beta, \Upsilon$ and $\underline{M}$ such that
\begin{equation}\label{uniform.estimate.parabolic.cylinders}
\kb_1 \frac{M^{\sigma\vartheta}}{t^{(d-\gamma)\vartheta}} \le u(t,x) \le \ka_1  \frac{M^{\sigma\vartheta}}{t^{(d-\gamma)\vartheta}}\quad\mbox{for any}\quad t\ge t_\star\quad\mbox{and any}\quad x \in \left\{|x|\le 3\,\Upsilon R(t)\right\}\,.
\end{equation}
\noindent\textit{Uniform and H\"older estimates in self-similar variables. }We first rescale $u$ in selfsimilar variables, according to~\eqref{change.fokker.planck}, and get $v(\tau, y)$. Analogously, the domain $\left\{|x|\le 3\,\Upsilon R(t)\right\}$ is transformed into $B_{3r}(0)$ where $r=\Upsilon \zeta$. Inequality~\eqref{uniform.estimate.parabolic.cylinders} reads in rescaled variables:
\begin{equation}\label{uniform.estimates.balls.fp}
\frac{\kb_1}{\zeta^{d-\gamma}}\,\vartheta^\vartheta M^{\sigma\vartheta} \le v(\tau, y) \le 2 \frac{\ka_1}{\zeta^{d-\gamma}}\,\vartheta^\vartheta M^{\sigma\vartheta} \quad\mbox{for any}\quad \tau \ge \frac{1}{\sigma}\log\frac{R(t_\star \vee 1 )}{R(0)}\quad\mbox{and any}\quad y \in B_{3r}(0)\,.
\end{equation}
By applying Lemma~\ref{Technical.Lemma} we deduce that there exist $\nu>0, \ka>0$ such that for any $\tau \ge \frac{1}{\sigma}\log\frac{R(t_\star \vee 1 )}{R(0)}+1$ we have that
\[
\lfloor v(\tau, \cdot) \rfloor_{C^\nu\left(B_{\frac{3}{2}r}(0)\right)}\leq \ka\,\,\  2 \frac{\ka_1}{\zeta^{d-\gamma}}\,\vartheta^\vartheta M^{\sigma\vartheta}\,.
\]
Using the subadditivity of $\lfloor \cdot \rfloor_{C^\nu\left(B_{r}(0)\right)}$ and the fact that the above estimates can also be applied to the Barenblatt profile $\mathcal{B}_M(y)$, we conclude that
\begin{equation}\label{uniform.cnu.estimate.fp}
\lfloor v(\tau, \cdot)-\mathcal{B}_M \rfloor_{C^\nu\left(B_{\frac{3}{2}r}(0)\right)}\leq 4\,\ka\, \frac{\ka_1}{\zeta^{d-\gamma}}\,\vartheta^\vartheta M^{\sigma\vartheta}\quad\mbox{for any}\quad \tau \ge \frac{1}{\sigma}\log\frac{R(t_\star \vee 1 )}{R(0)}+1\,.
\end{equation}

\noindent\textit{Convergence in $\LL^\infty$ norm.} We only prove the case $0<\gamma<d$, which is the most delicate, the case $\gamma\le0$ being simpler. In what follows it is convenient to assume that $r\ge2$, namely that $\Upsilon\ge \frac{2}{\zeta}$, we will overcome this technical assumption at the end of the proof. From the convergence in $\LL^1_\gamma$, formula \eqref{convergence.bareblatt.l1.nfp},  we deduce that there exists $\tau_\star$ such that for any $\tau\ge \tau_\star$ we have that $\|v(\tau, \cdot)-\mathcal{B}_M\|_{\LL^1_{\gamma}(B_{\frac{3}{2}r}(0))}\le \frac{|\gamma|}{d}$. We are in the position to apply inequality~\eqref{interpolation.inequality.gamma.positive} of Lemma~\ref{interpolation.lemma} to $v(\tau, \cdot)-\mathcal{B}_M$ and get that for any $\tau \ge \tau_\star \vee \frac{1}{\sigma}\log\frac{R(t_\star \vee 1 )}{R(0)}+1$
\begin{equation}\label{L1-Linfty.interp}
\|v(\tau, \cdot)-\mathcal{B}_M\|_{\LL^\infty(B_r(0))}\le C_{d,\gamma, \nu, p}\, \left(1+r\right)^\gamma\, \left(1+4\,\ka\, \frac{\ka_1}{\zeta^{d-\gamma}}\,\vartheta^\vartheta M^{\sigma\vartheta}\right)^\frac{d}{d+p\nu}\,\|v(\tau, \cdot)-\mathcal{B}_M\|_{\LL^1_\gamma(B_{\frac{3}{2}r}(0))}^\frac{\nu}{d+\nu}\,
\end{equation}
where we have used \eqref{uniform.cnu.estimate.fp}.
Since $\mathfrak{B}\ge \left(C(M)+r^\sigma\right)^\frac{-1}{1-m}$ on $B_r(0)$, it follows that
\[
\sup_{y \in B_r(0)}\left|\frac{v(\tau, y)-\mathcal{B}_M(y)}{\mathcal{B}_M(y)}\right| \le \left(C(M)+r^\sigma\right)^\frac{1}{1-m} \|v(\tau, \cdot)-\mathcal{B}_M\|_{\LL^\infty(B_r(0))}\,,
\]
which, combined with \eqref{L1-Linfty.interp} and the convergence in $\LL^1_\gamma$, formula \eqref{convergence.bareblatt.l1.nfp}, shows that the \emph{relative error} approaches zero as $\tau\rightarrow\infty$. Rescaling back, we finally obtain~\eqref{local.convergence.relative.error}, recalling that $\{|x| \le \Upsilon t^\vartheta\} \subset \{|x| \le \Upsilon R(t)\}$.

It only remains to overcome the technical assumption $\Upsilon\ge\frac{2}{\zeta}$. If $\Upsilon\le\frac{2}{\zeta}$ we can repeat the same argument for $\Upsilon=\frac{2}{\zeta}$. Next, we conclude that~\eqref{local.convergence.relative.error} takes place for any $\Upsilon'\le\frac{2}{\zeta} $  using that $\left\{|x|\le \Upsilon' R(t)\right\}\subset \left\{|x|\le \Upsilon R(t)\right\}$ whenever $\Upsilon'<\Upsilon$. The proof is now concluded.\qed

\section{Initial data in $\mathcal{X}$. Global Harnack Principle and uniform convergence in relative error}\label{section:GHP}
The space $\mathcal{X}$ is naturally invariant under the fast diffusion flow as explained in the introduction, see also Proposition \ref{curve.X.prop}. As a consequence, solutions belonging to this space possess some extra properties, that we summarize here:\vspace{-2mm}\begin{itemize}[leftmargin=*]\itemsep1pt \parskip1pt \parsep1pt
 \item The tail is essentially the same as the Barenblatt solution, the GHP holds, see Subsection \ref{thm.ghp}.
 \item The Uniform convergence in relative error (UREC) takes place, see Subsection \ref{Sect.UCRE}. Moreover, we also provide \textit{Almost Optimal Rates of convergence} in Subsection \ref{SSec.ratesX}, which turn out to be  sharp in some cases.

 \item  Boundary Harnack type inequalities hold true, and the behaviour at infinity of solutions does not depend on the mass, see Subsection \ref{ssection:behaviour-infty}.\vspace{-2mm}
 \end{itemize}
 In Section \ref{Sec4} we will show that the above properties are false if $u_0\not\in \mathcal{X}$.

\subsection{Upper Bound and proof of Theorem~\ref{thm.ghp}}
As already observed in the Introduction Theorem~\ref{thm.ghp} is divided in two parts: the upper bound and the lower bound of inequality~\eqref{ghp.inequality}. In this section we are going to discuss the upper bound,  the main result of this section is the following Theorem.
\begin{thm}\label{thm.upper}
Let $u$ be  the solution to \eqref{cauchy.problem} corresponding to the initial data $0\le u_0\in \LL^1_{\gamma,+}(\RR^d)$. Then, for any $t_0 > 0$ there exist $\overline{\tau},\overline{M}>0$,   explicitly given in \eqref{parameters.upper.bound},  such that
\begin{equation}\label{ghp.upper.bound}
u(t, x) \leq \mathfrak{B}(t+\overline{\tau},x; \overline{M})\qquad\mbox{for any $x \in \RR^d$ and any $t > t_0 $}\,,
\end{equation}
if and only if
\[
\mbox{$u_0$ satisfies \eqref{tail.condition}, i.e. $u_0 \in \mathcal{X}\,.$}
\]
\end{thm}

The proof of inequality~\eqref{ghp.upper.bound} is constructive and we are able to give values of $\overline{\tau}$ and $\overline{M}$, see formulae \eqref{parameters.upper.bound} at the end of the proof. Here we just point out that they depend on   $M, A, d,m,\gamma,\beta$ and $t_0$.

\noindent\textbf{Proof of Theorem \ref{thm.ghp}}. The proof is a simple combination of Theorem \ref{thm.upper} and Theorem \ref{lower.estimate.theorem}. \qed

\begin{rem}\label{limit.herrero.pierre.upper.rem}\rm We easily deduce from the above upper bound that
\begin{equation}\label{herrero.pierre.limit.upper.inq}
\limsup\limits_{|x|\rightarrow \infty}\,u(t,x)\,|x|^\frac{\sigma}{1-m}\le b_1 \left(t +\overline{\tau} \right)^\frac{1}{1-m}\,.
\end{equation}
Equality is achieved by the Barenblatt solution translated in time by $\overline{\tau}$. Notice that this maximal tail behaviour only holds for $u_0\in \mathcal{X}$, in which case it matches the optimal minimal behaviour given in Corollary \ref{herrero.pierre.limit.corollary}. These two pieces of information combine well and allow to deduce the sharp behaviour at infinity, see Section \ref{ssection:behaviour-infty}, Corollary \ref{liminf-limsup}.
\end{rem}
\noindent {\bf Proof of Theorem \ref{thm.upper}: } Let us first explain the strategy of the proof. We will prove inequality \eqref{ghp.upper.bound} only at time $t_0$, then, by comparison (see for instance \cite[Corollary 9]{Bonforte2017})  it will hold for any $t\geq t_0$.  The proof is divided in several steps: first, we estimate the solution $u(t_0, x)$ in two different regions (on $B_{R_1}(0)$ and on $B_{R_1}(0)^c$, with $R_1$ to be chosen later), then we find conditions on $\overline{\tau}$ and $\overline{M}$ necessary for inequality \eqref{ghp.upper.bound} to hold. Finally, we show that such conditions can be fulfilled providing an explicit expression of $\overline{\tau}$ and $\overline{M}$ in terms of $t_0, M$ and $A$.

In this proof we will make use of the following estimate, \cite[Theorem 1.2]{Bonforte2019a}: there exists $ \ka_1,\ka_2>0 $  such that for any $ t >0 $,  $x_0 \in \RR^d$ and any $R_0 \in \left[|x_0|/16, |x_0|/32\right]$ (any $R_0>0$ if $x_0=0$) we have that
 \begin{equation}\label{local.upper.estimte}
\sup_{y\in B_{R_0}(x_0) }{u\left(t,y\right)} \leq \frac{\ka_1}{t^{(d-\gamma)\vartheta}} \left[\int_{B_{2R_0}(x_0)}{|u_0(y)| \ |y|^{-\gamma} \ \dy} \right]^{\sigma  \vartheta} + \ka_2 \left[\frac{t}{R_0^{\sigma}} \right]^{\frac{1}{1-m}},
 \end{equation}
where $\vartheta$ and $\sigma$ are defined in \eqref{parameters.sigma.theta}. The constants $\ka_1, \ka_2$ are explicit and depend only on $d, m,\gamma$ and $\beta$.  The constant $\ka_1$ is the same one which appears in the smoothing effect given in inequality~\eqref{global.upper.estimate}.

\noindent\textsc{Estimate inside a ball.} We want to find suitable conditions on $\overline{M}, \overline{\tau}$ and $R_1$ such that
\begin{equation}\label{upper.ghp.t0}
u(t_0, x) \leq \mathfrak{B}(t_0+\overline{\tau},x; \overline{M})=\frac{(t_0+\overline{\tau})^{\frac{1}{1-m}}}{\left[b_0\frac{(t_0+\overline{\tau})^{\sigma\vartheta}}{M^{\sigma\vartheta(1-m)}}+b_1 |x|^{\sigma}\right]^\frac{1}{1-m}}\,,\qquad\mbox{holds for all $|x|\leq R_1$.}
\end{equation}
Recall that $M= \int_{\RR^d}u_0 |x|^{-\gamma} \dx$. Inequality \eqref{local.upper.estimte} implies that
\begin{equation*}
u(t_0, x) \leq \ka_1 t_0^{-(d-\gamma)\vartheta} M^{\sigma\vartheta}\qquad\mbox{for any $x \in \RR^d$ and $t_0>0$}.
\end{equation*}
To deduce the above from~\eqref{local.upper.estimte} it suffices to take $x_0=0$ and let $R_0\rightarrow \infty$.
In view of the above inequality, to prove \eqref{upper.ghp.t0} it is enough to find suitable $\overline{M}, \overline{\tau}$ and $ R_1$ such that
\begin{equation*}
\ka_1 \frac{M^{\sigma\vartheta}}{t_0^{(d-\gamma)\vartheta} } \leq \frac{\left(t_0 + \overline{\tau}\right)^{\frac{1}{1-m}}}{\left[b_0\frac{\left(t_0 + \overline{\tau} \right)^{\sigma\vartheta}}{\overline{M}^{\sigma\vartheta(1-m)}}+b_1 |x|^{\sigma}\right]^{\frac{1}{1-m}}}
\qquad\mbox{for any $|x| \leq R_1$.}
\end{equation*}
Since the righthand side is decreasing in $|x|$ it suffices to have the previous inequality at $|x|=R_1$, i.e.
\begin{equation}\label{condition.upper.bounds.2}
b_0\frac{\left(t_0 + \overline{\tau} \right)^{\sigma\vartheta}}{\overline{M}^{\sigma\vartheta(1-m)}}+b_1 R_1^{\sigma} \leq \frac{(t_0 + \overline{\tau}) t_0^{(d-\gamma)\theta(1-m)}}{\ka_1^{1-m}M^{\sigma\vartheta(1-m)}}.
\end{equation}
Inequality \eqref{condition.upper.bounds.2} is nothing but a first condition on $\overline{M}, \overline{\tau}$ and $ R_1$ in order to guarantee the validity of \eqref{upper.ghp.t0}.

\medskip

\noindent\textsc{Estimate outside a ball. }The goal of this step is to extend inequality \eqref{upper.ghp.t0} outside a ball, namely for all $|x|\geq R_1$. This will end up to conditions on $\overline{M}, \overline{\tau}$ and $R_1$ different from \eqref{condition.upper.bounds.2}. In the next step we will take care of checking the compatibility of the two conditions.

We first prove that for any fixed $ t_0 > 0 $ there exists $C_1=C_1(t_0, A)>0$ such that
\begin{equation}\label{upper.decay.1}
u(t_0, x) \leq \frac{ C_{1}}{|x|^{\frac{\sigma}{(1-m)}}}\qquad\mbox{for any $ |x|> R_1$}.
\end{equation}
Let $x \in \RR^d$, $|x|\geq R_1$ and let $R$ be such that $B_{2R}(x) \subset B_{2R}(0)^c$, for instance $R=|x|/16$. Applying inequality \eqref{local.upper.estimte} to $u(t_0, x)$ in the ball $B_R(x)$, we get
\begin{equation*}\begin{split}
u(t_0, x)
& \leq \frac{\ka_1}{t_0^{(d-\gamma)\vartheta}} \left[\int_{B_{2R}^c(0)}u_0(y)\,|y|^{-\gamma}\dy\right]^{\sigma  \vartheta} + \ka_2(16)^{-\frac{\sigma}{1-m}} \left(\frac{t_0}{|x|^{\sigma}}\right)^{\frac{1}{1-m}} \\
& \leq \frac{\ka_1 8^\frac{\sigma}{1-m}}{t_0^{(d-\gamma)\vartheta}} \frac{A^{\sigma \vartheta}}{|x|^{\frac{\sigma}{1-m}}} +\frac{\ka_2}{16^{\frac{\sigma}{1-m}}}\left(\frac{t_0}{|x|^{\sigma}}\right)^{\frac{1}{1-m}}  \quad \leq \frac{C_1}{|x|^{\frac{\sigma}{1-m}} }\,,
\end{split}\end{equation*}
where in the third line we have used that $\int_{B^{c}_R(0)}u_0|x|^{-\gamma}\dx \leq A R^{(d-\gamma)-\frac{2+\beta-\gamma}{1-m}}$ with $R=|x|/16$ and that  $C_1=C_1(t_0, A)$ is given by
\begin{equation*}
C_1 = 8^{\frac{\sigma}{1-m}}\frac{\ka_1}{t_0^{(d-\gamma)\vartheta}}\,A^{\sigma\vartheta} + \frac{\ka_2}{16^{\frac{\sigma}{1-m}}} t_0^{\frac{1}{1-m}}\,.
\end{equation*}
Hence inequality \eqref{upper.decay.1} holds. It only remains to show that %
\begin{equation}\label{condition.upper.bounds.0}
\frac{C_1}{|x|^{\sigma/(1-m)}} \leq \frac{\left(t_0 + \overline{\tau}\right)^{\frac{1}{1-m}}}{\left[b_0\frac{\left(t_0 + \overline{\tau} \right)^{\sigma\vartheta}}{\overline{M}^{\sigma\vartheta(1-m)}}+b_1 |x|^{\sigma}\right]^{\frac{1}{1-m}}}\qquad\mbox{for any $|x| \geq R_1$}.
\end{equation}
This will give a condition on $\overline{\tau}, \overline{M} $ and $R_1$, as we explain next. Indeed, the above inequality is equivalent to
\begin{equation*}
 b_1 C_1^{1-m} + b_0\frac{\left(t_0 + \overline{\tau} \right)^{\sigma\vartheta}}{|x|^{\sigma}\overline{M}^{\sigma\vartheta(1-m)}}  \leq t_0+\overline{\tau}.
\end{equation*}
It is indeed enough to choose $R_1>0$ such that
\begin{equation}\label{condition.upper.bounds.1}
b_1 C_1^{1-m} + b_0\frac{\left(t_0 + \overline{\tau} \right)^{\sigma\vartheta}}{R_1^{\sigma}\overline{M}^{\sigma\vartheta(1-m)}} \leq t_0+\overline{\tau},
\end{equation}
since the second term in left-hand side is decreasing in $|x|$.
We conclude that inequality \eqref{upper.ghp.t0} holds for any $|x|\geq R_1$ whenever $\overline{\tau}, \overline{M} $ and $R_1$ satisfy condition \eqref{condition.upper.bounds.1}.

\medskip

\noindent\textsc{Compatibility among the conditions \eqref{condition.upper.bounds.2} and \eqref{condition.upper.bounds.1}. }We only need to show the compatibility of the conditions that imply the main estimates of the previous steps, i.e. that inequality \eqref{upper.ghp.t0} holds for all $x\in \RR^d$. The two conditions \eqref{condition.upper.bounds.2} and \eqref{condition.upper.bounds.1}correspond to the following system of inequalities
\begin{equation*}
\textbf{(A)}=\begin{cases}
\begin{aligned}
b_1 C_1^{1-m} R_1^\sigma + b_0\frac{\left(t_0 + \overline{\tau} \right)^{\sigma\vartheta}}{\overline{M}^{\sigma\vartheta(1-m)}} &\leq R_1^{\sigma}(t_0+\overline{\tau})\,, \\
b_1 R_1^{\sigma} + b_0\frac{\left(t_0 + \overline{\tau} \right)^{\sigma\vartheta}}{\overline{M}^{\sigma\vartheta(1-m)}} &\leq \frac{(t_0 + \overline{\tau}) t_0^{(d-\gamma)\theta(1-m)}}{\ka_1^{1-m}M^{\sigma\vartheta(1-m)}}\,. \\
\end{aligned}
\end{cases}
\end{equation*}
It is convenient to simplify the above system in order to be able to make explicit choices of $\overline{\tau}, \overline{M} $ and $R_1$.
The first simplification is the following:
\begin{equation}
\textbf{(B)}=\begin{cases}
\begin{aligned}
b_1 \left(1 \vee C_1\right)^{1-m} R_1^\sigma &\leq \frac{t_0+\overline{\tau}}{2} \left[R_1^\sigma \wedge\frac{t_0^{(d-\gamma)\theta(1-m)}}{\ka_1^{1-m}M^{\sigma\vartheta(1-m)}}\right]    \,,  \\
b_0\frac{\left(t_0 + \overline{\tau} \right)^{\sigma\vartheta}}{\overline{M}^{\sigma\vartheta(1-m)}} &\leq \frac{t_0+\overline{\tau}}{2} \left[R_1^\sigma \wedge \frac{t_0^{(d-\gamma)\theta(1-m)}}{\ka_1^{1-m}M^{\sigma\vartheta(1-m)}}\right] \,.\\
\end{aligned}
\end{cases}
\end{equation}
It is clear that any choice of $\overline{\tau}, \overline{M} $ and $R_1$ that satisfies $\textbf{(B)}$  also satisfies $\textbf{(A)}$.
We need a further simplification, but this time we will choose $R_1=R_1(R_0, t_0, M)$ in a particular way, as follows
\begin{equation*}
R_1 := \left(\frac{t_0^{(d-\gamma)\theta}}{\ka_1M^{\sigma\vartheta}}\right)^{\frac{1-m}{\sigma}}\quad\mbox{so that}\quad R_1^\sigma=\frac{t_0^{(d-\gamma)\theta(1-m)}}{\ka_1^{1-m}M^{\sigma\vartheta(1-m)}}\,,
\end{equation*}
and system $\textbf{(B)}$ simplifies to
\begin{equation*}
\textbf{(B')}=\begin{cases}
\begin{aligned}
b_1 \left(1 \vee C_1\right)^{1-m}  \leq \frac{\left(t_0+\overline{\tau}\right)}{2}\,,  \\
b_0\frac{\left(t_0 + \overline{\tau} \right)^{\sigma\vartheta}}{\overline{M}^{\sigma\vartheta(1-m)}} \leq \frac{\left(t_0+\overline{\tau}\right)}{2} \frac{t_0^{(d-\gamma)\theta(1-m)}}{\ka_1^{1-m}M^{\sigma\vartheta(1-m)}} \,,\\
\end{aligned}
\end{cases}
\Longrightarrow\begin{cases}
\begin{aligned}
\overline{\tau} &\geq \left[ 2 b_1 \left(1 \vee C_1\right)^{1-m} -t_0 \right]\,,  \\
\overline{M} &\geq \left(2 b_0 \ka_1^{1-m}\right)^{\frac{1}{\sigma\vartheta(1-m)}}  \, \left(\frac{t_0 + \overline{\tau}}{t_0 }\right)^{\frac{d-\gamma}{\sigma}}\, M\,,\\
\end{aligned}
\end{cases}
\end{equation*}
It is now clear that choosing $\overline{\tau}=\overline{\tau}(t_0, M_\infty, C_1, R_1)$ and $M=M(\overline{\tau}, t_0, M)$ of the form
\begin{equation*}
\overline{\tau} := 0 \vee \left[ 2 b_1 \left(1 \vee C_1\right)^{1-m} -t_0 \right] \quad \mbox{and}\quad
\overline{M}:= \left(2 b_0 \ka_1^{1-m}\right)^{\frac{1}{\sigma\vartheta(1-m)}}  \, \left(\frac{t_0 + \overline{\tau}}{t_0 }\right)^{\frac{d-\gamma}{\sigma}}\, M\,,
\end{equation*}
implies the validity of  the two inequalities of system $\textbf{(B')}$, hence of system $\textbf{(B)}$, and finally of $\textbf{(A)}$.

\noindent\textsc{Values of the constants. }Letting  $A:=|u_0|_\mathcal{X}$ and $M:=\|u_0\|_{\LL^1_\gamma(\RR^d)}$\,, we have
\begin{equation}\label{parameters.upper.bound}\begin{split}
\overline{\tau} := 0 \vee \left\{ 2 b_1 \left[1 \vee \left(8^{\frac{\sigma}{1-m}}\frac{\ka_1A^{\sigma\vartheta}}{t_0^{(d-\gamma)\vartheta}} + \frac{\ka_2t_0^{\frac{1}{1-m}}}{16^{\frac{\sigma}{1-m}}}  \right)\right]^{1-m} -t_0 \right\}\,, \\
\overline{M}:= \left(2 b_0 \ka_1^{1-m}\right)^{\frac{1}{\sigma\vartheta(1-m)}} \left(\frac{t_0 + \overline{\tau}}{t_0 }\right)^{\frac{d-\gamma}{\sigma}} M\,,
\end{split}\end{equation}
where $\ka_1, \ka_2>0$ depend on $d, m, \gamma, \beta$, and they have an explicit expression given at the end of the proof of Theorem 1.2 in \cite{Bonforte2019a}.  The proof is concluded. \qed

\subsection{Proof of Theorem \ref{equivalent.condition}. \emph{Uniform Relative Error Convergence.}}\label{Sect.UCRE}
In this subsection we prove the sufficiency part of Theorem \ref{equivalent.condition}. The converse implication will be proven in the next Section. We just recall that partial results in the non-weighted case, $\gamma=\beta=0$, have been proven in \cite{Carrillo2003,Kim2006,Vazquez2003, Vazquez2006,Bonforte2006,Blanchet2009,Bonforte2010,Bonforte2010c}. For the weighted case, see \cite{Bonforte2017a,Bonforte2017}.
\begin{thm}[UREC]\label{global.convergence.relative.error.thm}
Let $m \in \left(m_c, 1\right)$ and let $u$ be a solution to \eqref{cauchy.problem} with initial data $0 \leq u_0 \in \mathcal{X}\setminus\{0\}$. Then we have that
\begin{equation}\label{global.convergence.relative.error}
\lim_{t\rightarrow 0} \Big\|\frac{u(t,x)-\mathfrak{B}(t, x; M)}{\mathfrak{B}(t, x; M)} \Big\|_{\LL^\infty(\RR^d)}=0\,,\qquad\mbox{where}\qquad M=\|u_0\|_{\LL^1_\gamma(\RR^d)}.
\end{equation}
\end{thm}
\noindent\textbf{Proof.} It is convenient to work in self-similar variables: we transform $u(t,x)$ into $v(\tau, y)$ accordingly to formula~\eqref{change.fokker.planck}. We will prove that for any $\varepsilon>0$ there exits $\tau_{\varepsilon}>0$ such that
\begin{equation}\label{global.convergence.relative.error.inq}
\Big\|\frac{v(\tau,y)-\mathfrak{B}_{M}(y)}{\mathfrak{B}_{M}(y)} \Big\|_{\LL^\infty(\RR^d)} < 2\varepsilon\qquad\mbox{for any }\quad \tau\ge \tau_\varepsilon\,.
\end{equation}
We argue that we only need to prove the following claim.

\noindent\textit{Claim.} For any $1>\varepsilon>0$ there exists $\rho_\varepsilon>$ and $\overline{\tau}_\varepsilon>0$ such that
\begin{equation}\label{claim.inequality}
\sup_{|y|\ge \rho_\varepsilon}\left|\frac{v(\tau, y)-\mathcal{B}_M(y)}{\mathcal{B}_M(y)}\right|<\,\varepsilon\quad\mbox{for any}\quad \tau \ge \overline{\tau}_\varepsilon\,.
\end{equation}
Indeed, once the Claim is proven, we just combine it with the convergence inside parabolic cones, i.e. the main result of Theorem~\ref{local.convergence.relative.error.thm} and obtain inequality~\eqref{global.convergence.relative.error.inq} as follows:

\begin{equation*}\begin{split}
\Big\|\frac{v(\tau,y)}{\mathfrak{B}_{M}(y)} -1 \Big\|_{\LL^\infty(\RR^d)} \le \Big\|\frac{v(\tau,y)}{\mathfrak{B}_{M}(y)} -1 \Big\|_{\LL^\infty(\{|y|\le \Upsilon\})} + \Big\|\frac{v(\tau,y)}{\mathfrak{B}_{M}(y)} -1 \Big\|_{\LL^\infty(\{|y|\ge \Upsilon\})}\le 2\,\varepsilon\,.
\end{split}\end{equation*}
Recall that the change of variables~\eqref{change.fokker.planck} transforms the parabolic cones $\{|x| \le \Upsilon R(t)\}$ into balls  $\{|y| \le \Upsilon \}$.

\noindent\textit{Proof of the Claim.}  Let $t_0, R_0>0$ be such that $ \|u_0\|_{\LL^1_{\gamma}(B_{R_0}(0))} > 0$. We know by  Theorem~\ref{thm.ghp} that
\[
\mathfrak{B}(t-\underline{t},x; \underline{M}) \leq u(t, x) \leq \mathfrak{B}(t+\overline{t},x; \overline{M})\,.
\]
for suitable  $\overline{t}, \underline{t}>0$ and $\overline{M}, \underline{M}>0$. As a consequence, recalling the  change of variables~\eqref{change.fokker.planck}, we get
\begin{equation}\label{ghp.NFP}
a(t)^{d-\gamma}\,\left(1  \wedge a(t)\right)^{\frac{-\sigma}{1-m}}\,\, \mathfrak{B}_{\underline{M}}\left(y\right)\,\,\le v(\tau, y)\le \,\, b(t)^{d-\gamma}\,\left(b(t)\vee1\right)^\frac{-\sigma}{1-m}\,\,\mathfrak{B}_{\overline{M}}\left(y\right)\,,
\end{equation}
where $R(t)=R_\star(t+1)$ and
\[
\tau=\frac{1}{\sigma}\log\frac{R(t)}{R(0)}\qquad a(t)=\frac{R_\star(t+1)}{R_\star(t+\underline{t})}\,\quad\mbox{and}\quad\,b(t)= \frac{R_\star(t+1)}{R_\star(t+\overline{t})}\,.
\]
Since $a(t), b(t) \rightarrow 1$ as $t\rightarrow \infty$ we deduce that there exists $\tau_\varepsilon>0$ such that
\begin{equation}\label{time.convergence.ghp.NFP}
\left(1-\frac{\varepsilon}{3}\right)\,\mathcal{B}_{\underline{M}}\left(y\right)\le \,v(\tau, y)\, \le\left(1+\frac{\varepsilon}{3}\right)\,\mathcal{B}_{\overline{M}}\left(y\right)\qquad\mbox{for every $\tau>\tau_\varepsilon$.}
\end{equation}
Recall that all the Barenblatt solutions $\mathcal{B}_{M}$ have the same behaviour at infinity, which is independent of the mass $M$, namely $\lim\limits_{|y|\rightarrow \infty}\mathcal{B}_{M_1}(y)/\mathcal{B}_{M_2}(y)=1$for any $M_1, M_2 >0$. Hence,  there exists $\rho_\varepsilon=\rho_\varepsilon(\underline{M}, \overline{M})>0$ such that
\[
 1-\frac{\varepsilon}{3}\le \frac{\mathcal{B}_{\underline{M}}(y)}{\mathcal{B}_M(y)}\quad\mbox{and}\quad \frac{\mathcal{B}_{\overline{M}}(y)}{\mathcal{B}_M(y)}\le 1+\frac{\varepsilon}{3}\,,\qquad\mbox{for any}\quad |y|\ge \rho_\varepsilon\,.
\]
Combining the above inequality with~\eqref{time.convergence.ghp.NFP} we obtain the proof of the \textit{Claim}.  The proof is concluded. \qed
\subsection{Proof of Theorem~\ref{equivalent.condition}. The necessary part}\label{sec:proof.thm.equivalent.condition}
We have already shown that the tail condition \eqref{tail.condition} implies the Uniform Convergence in Relative Error: this result is contained in Theorem~\ref{global.convergence.relative.error.thm}. As a consequence, the sufficiency part of Theorem~\ref{equivalent.condition} has already been proven, and we only need to prove the converse implication: if a solution converges uniformly in relative error, then the initial datum $u_0$ satisfies the tail condition \eqref{tail.condition} or equivalently $u_0\in \mathcal{X}$.   The proof of the necessary part  is based on the following Lemma, which bears some similarity with a result by Herrero and Pierre \cite[Lemma 3.1]{Herrero1985}, valid in the non-weighted case, and then generalized by the authors in \cite{Bonforte2019a} to the present weighted case. Notice that we need to integrate on the complementary set of balls, instead of balls as in \cite{Herrero1985,Bonforte2019a}. The proof is quite similar hence we only sketch it.
\begin{lem}\label{lem.herrero.pierre.inverse}
Let $m \in \left(m_c, 1\right)$ and let $u$ be a solution to \eqref{cauchy.problem} with initial data $0 \leq u_0 \in \LL^1_\gamma(\RR^d)$, then for any $R>0$ and for any $t, s\ge0$ there exist constants $C_1, C_2>0$ which depend on $m, d, \gamma, \beta$ such that
\begin{equation}\label{Herrero.Pierre.inverse.2}
\int_{B_{2R}^c(0)}u(t,x)|x|^{-\gamma}\dx  \le C_1\,\int_{B_{R}^c(0)}u(s,x)|x|^{-\gamma}\dx   + C_2\,|t-s|^\frac{1}{1-m}\, R^{(d-\gamma)-\frac{\sigma}{1-m}}\,.
\end{equation}
\end{lem}
\noindent\textbf{Proof.} Let $R>0$ and define $A(R)=B_{2R}(0)\setminus B_{R}(0)$. Let $0\le\psi\in C^\infty(\RR^d)$ be such that $\psi=1$ in $B_{2R}(0)^c$ and $\psi=0$ in $B_{R}(0)$. Let us sketch the proof. For any $t>0$, let us formally compute
\begin{equation}\label{proof.1}\begin{split}
\left|\frac{\mathrm{d}}{\dt}\int_{\RR^d}u(t, x)\psi(x)\right. &\left.\frac{\dx}{|x|^\gamma} \right|= \left|\int_{A(R)}u^m(t, x)\,|x|^\gamma\mathrm{div}\left(|x|^{-\beta}\nabla \psi\right)\frac{\dx}{|x|^\gamma}\right| \\
& \le \int_{A(R)} u^m(t,x)\psi^m(x)\,\psi^{-m}(x) \Big||x|^\gamma\mathrm{div}\left(|x|^{-\beta}\nabla \psi\right)\Big|\frac{\dx}{|x|^\gamma} \\
& \le \left(\int_{A(R)} u(t,x)\psi(x)\frac{\dx}{|x|^\gamma}\right)^m\,\left(\int_{A(R)}\psi^\frac{-m}{1-m}\Big||x|^\gamma\mathrm{div}\left(|x|^{-\beta}\nabla \psi\right)\Big|^\frac{1}{1-m}\frac{\dx}{|x|^\gamma}\right)^{1-m} \\
& \le C\,  \left(\int_{\RR^d} u(t,x)\psi(x)\frac{\dx}{|x|^\gamma}\right)^m\,R^{(d-\gamma)(1-m)-\left(2+\beta-\gamma\right)}\,,
\end{split}\end{equation}
where we have used H\"{o}lder and the fact that $\psi^\frac{-m}{1-m}\Big||x|^\gamma\mathrm{div}\left(|x|^{-\beta}\nabla \psi\right)\Big|^\frac{1}{1-m}\le c\, R^\frac{-2-\beta+\gamma}{1-m}$, which can be easily derived following the lines of the proof of Lemma 5.2 of~\cite{Bonforte2019a}. Recall that the integration by parts done in the first line of~\eqref{proof.1} presents no difficulties, since $\psi=0$ in a neighborhood of the origin, where the weight $|x|^{-\gamma}$ can be singular or degenerate. Integrating the differential inequality~\eqref{proof.1}, as done in~\cite[Lemma 3.1]{Herrero1985}, one obtains~\eqref{Herrero.Pierre.inverse.2}. The proof is concluded. A rigorous proof starts from the integrated version of \eqref{proof.1} (which follows by definition of weak solutions) and then follows by a Grownwall-type argument. \qed

\medskip

\noindent\textbf{Proof of the necessary part of Theorem~\ref{equivalent.condition}}. The proof is based on Lemma~\ref{lem.herrero.pierre.inverse} proven in the Appendix, that we restate here for reader's convenience in the form that we need. Let $u$ be a solution to \eqref{cauchy.problem} with initial data $0 \leq u_0 \in \LL^1_\gamma(\RR^d)$. Then, there exist constants $C_1, C_2>0$ which depend on $m, d, \gamma, \beta$ such that for any $R>0$ and for any $t\ge0$
\begin{equation}\label{Herrero.Pierre.inverse}
\int_{B_{2R}^c(0)}u_0(x)|x|^{-\gamma}\dx  \le C_1\,\int_{B_{R}^c(0)}u(t,x)|x|^{-\gamma}\dx   + C_2\,t^\frac{1}{1-m}\, R^{(d-\gamma)-\frac{\sigma}{1-m}}\,.
\end{equation}

Let us proceed with the rest of the proof. Assume now that~\eqref{global.convergence.relative.error.equivalence} holds, hence there exists a time $\overline{t}>0$ such that for any $x\in \RR^d$
\[
\left|\frac{u(\overline{t},x)-\mathfrak{B}(\overline{t},x;M)}{\mathfrak{B}(\overline{t},x;M)}\right|<1\,,\qquad\mbox{hence}\qquad
u(\overline{t}, x)\le 2\,\mathfrak{B}(\overline{t}, x; M)\,.
\]
Integrating the latter inequality over $B_R^c$ we get that  there exists a constant $\kappa > 0$ which depend on $m, d, \gamma, \beta$ and on $\overline{t}$ such that
\[
\int_{B_{R}^c(0)}u(\overline{t}, x) |x|^{-\gamma}\dx   \le \, \kappa \, R^{d-\frac{\sigma}{1-m}}\qquad\mbox{for all $R>0$}\,.
\]
Combining~\eqref{Herrero.Pierre.inverse} at time $t=\overline{t}$ with the above estimate, we conclude that for any $R>0$
\[
(2R)^{\frac{\sigma}{1-m}-d}\, \int_{B_{2R}^c(0)}u_0(x)|x|^{-\gamma}\dx \le C\,\kappa^{\frac{1}{1-m}} + C_{m, d, \gamma, \beta}^\frac{1}{1-m}\,|\overline{t}|^\frac{1}{1-m}\,.
\]
As a consequence, the initial data $u_0$ satisfies the tail condition~\eqref{tail.condition} and the proof is concluded.\qed

\subsection{Harnack inequalities for quotients and sharp behaviour at infinity}\label{ssection:behaviour-infty}
In this Subsection we show a results which can be interpreted as a boundary Harnack inequality, since it extends to the whole space the Harnack inequalities on Parabolic Cones of Theorem \ref{thm.Harnack.Parabolic.Cones}. More precisely, there exist a constant $H>0$ such that for any $t>0$ (large enough) we have
\begin{equation}\label{Boundary-Harnack}
\frac{u(t,x)}{u(t,y)}\le H \,\frac{\mathfrak B(t, x; M)}{\mathfrak B(t, y; M)} \qquad\mbox{for all $x,y\in \RR^d$.}
\end{equation}
The above inequality is equivalent to \eqref{Boundary-Harnack.ineq.thm} and provides interesting information about the  behaviour at infinity of solutions to~\eqref{cauchy.problem} with data $u_0\in\mathcal{X}$. Somehow, the behaviour at infinity does not depend on the mass: indeed, we can show that for all $\overline{M},\underline{M}>0$, there exists $\overline{\tau}>0$ such that
\[
1\le \frac{\limsup\limits_{|x|\rightarrow \infty}\,u(t,x)\,\mathfrak{B}^{-1}(t,x; \overline{M})}{\liminf\limits_{|x|\rightarrow \infty}\,u(t,x)\,\,\mathfrak{B}^{-1}(t ,x; \underline{M})} \le   \left(1+\frac{\overline{\tau}}{t}\right)^\frac{1}{1-m}\qquad\mbox{if and only if }\qquad u_0\in \mathcal{X}\setminus\{0\}.
\]
The above inequalities are sharp, and are an equivalent statement of \eqref{liminf-limsup}. It is remarkable that equality is attained by Barenblatt profiles, possibly with different mass.
\begin{thm}
Let $u$ be a solution to \eqref{cauchy.problem} with $0\le u_0 \in \mathcal{X}$ such that $ \|u_0\|_{ m, \gamma, \beta}= A $ and $ \|u_0\|_{\LL^1_\gamma(\RR^d)}=M$ and let $R_0>0$ be such that $ \|u_0\|_{\LL^1_{\gamma}(B_{R_0}(0))} = M/2$. Then there exists a constant $\mathrm H>0$, which depends only on $m, d, \gamma, \beta$, such that
\begin{equation}\label{Boundary-Harnack.ineq.thm}
\sup_{x\in \RR^d}\frac{u(t,x)}{\mathfrak B(t, x; M)} \le \mathrm H \, \inf_{x \in \RR^d}\frac{u(t,x)}{\mathfrak B(t, x; M)}\quad\mbox{for any}\quad t\ge \overline{t}
\end{equation}
where
\[
\overline{t}= 3\, \max\left\{A^{1-m}\,\left(\frac{\ka_1}{\ka_2}\right)^\frac{1-m}{\sigma\vartheta}\,2^\frac{7}{\vartheta},  \kappa_*\,R_0^\frac{1}{\vartheta}\,\left(M/2\right)^\frac{1}{1-m}\right\}\,.
\]
The constants $\ka_1, \ka_2$ and $\kappa_*$ are as in~\eqref{local.upper.estimte} and as in~\eqref{Minimal.life.time} respectively.

Moreover, we obtain a characterization of the sharp behaviour at infinity, namely we have that
\begin{equation}\label{liminf-limsup}
1\le \frac{\limsup\limits_{|x|\rightarrow \infty}\,u(t,x)\,|x|^\frac{\sigma}{1-m}}{\liminf\limits_{|x|\rightarrow \infty}\,u(t,x)\,|x|^\frac{\sigma}{1-m}} \le   \left(1+\frac{\overline{\tau}}{t}\right)^\frac{1}{1-m}\qquad\mbox{if and only if }\qquad u_0\in \mathcal{X}\setminus\{0\}.
\end{equation}
Here, $\overline{\tau}$ depends on the initial data and is as in Theorem \ref{thm.upper}.
\end{thm}
\noindent\textbf{Proof.} We begin by proving inequality \eqref{Boundary-Harnack.ineq.thm}. In what follows we shall assume without loss of generality that
\begin{equation}\label{thm3.4.000}
\kappa_2^\frac{1}{1-m} \ge b_1^{-1}\,2^{4\sigma+m-2 }
\end{equation}
indeed, since $\kappa_2$ comes from the upper bound~\eqref{local.upper.estimte} we can choose it as large as needed. By applying Theorem~\ref{ghp.upper.bound} at time $t_0=A^{1-m}\,\left(\ka_1/\ka_2\right)^\frac{1-m}{\sigma\vartheta}\,2^\frac{7}{\vartheta}$ and Theorem~\ref{lower.estimate.theorem} at time $t_1= \kappa_\star\,R_0^\frac{1}{\vartheta}\,\left(M/2\right)^\frac{1}{1-m}$ we obtain that for any $t\ge \overline{t}$ the following inequality holds
\begin{equation}\label{ghp.applied}
\mathfrak B(t-\underline{\tau}, x; \underline{M})\le u(t,x)\le \mathfrak B(t+\overline{\tau}, x; \overline{M})\,,
\end{equation}
where
\[
\underline{\tau}=t_1/2=\left(\frac{\kappa_\star}{2}\right)\,R_0^\frac{1}{\vartheta}\left(\frac{M}{2}\right)^\frac{1}{1-m}\quad\quad  \underline{M}=b\,\frac{M}{2}
\]
and
\[
\overline{\tau}=(b_1\,2^{2-m-4\sigma}\left(\ka_2\right)^{1-m}-1)t_0 \quad\quad \overline{M}= \left(2b_0\ka_1\right)^{1-m}\left(b_1\,2^{2-m-4\sigma}\left(\ka_2\right)^{1-m}\right)^\frac{d-\gamma}{\sigma}M\,.
\]
Here is the point where the assumption \eqref{thm3.4.000} enters the game, since it implies that $\overline{\tau}\ge0$. By inequality~\eqref{ghp.applied} it is enough to show that there exist a constant $\mathrm H$ such that for any $t \ge \overline{t}$
\[
\sup_{x\in\RR^d} \frac{\mathfrak B(t+\overline{\tau}, x; \overline{M})}{\mathfrak B(t, x; M)} \le \mathrm H \inf_{x\in\RR^d} \frac{\mathfrak B(t-\underline{\tau}, x; \underline{M})}{\mathfrak B(t, x; M)}\,.
\]
A simple computation, which is left to the interested reader, shows that the previous inequality holds with a constant which depends only on $m, d, \gamma, \beta$ and not on the mass $M$ neither on the parameter $A$. The proof of \eqref{Boundary-Harnack.ineq.thm} is then concluded.

\noindent {\it Proof of \eqref{liminf-limsup}.~}We just  combine inequality~\eqref{Herrero.Pierre.Limit} of Corollary~\ref{herrero.pierre.limit.corollary} with inequality~\eqref{herrero.pierre.limit.upper.inq} of Remark~\ref{limit.herrero.pierre.upper.rem}.\qed

\subsection{Rates of convergence in $\mathcal{X}$}\label{SSec.ratesX}

As we have mentioned in Subsection~\ref{sec.dyn.sys.intro}, we know that solutions starting from $0\le u_0\in\mathcal{X}$, will eventually converge to a Barenblatt profile $\mathcal{B}_M$ (with the same mass as $u_0$), i.e. an element of the manifold $\mathcal{M}_{\mathcal{B}}$. The natural question that we address here is: are there ``universal rates'' of convergence towards $\mathcal{M}_{\mathcal{B}}$? More precisely:
\begin{center}
\textit{In self-similar variables can we find a speed of convergence to the\\ stationary profile valid for all solutions starting from data in $\mathcal{X}$? }
\end{center}

The answer to this question is delicate and can not be easily given for all $m\in (0,1)$, neither for all $m\in (m_c,1)$.
Some preliminary remarks are in order. In the case $\gamma=\beta=0$ the question has a long story (see \cite{Denzler2015}). When $\frac{d-2}{d}<m<1$ it has been proven in \cite{Carrillo2000,Carrillo2001,Carrillo2002,DelPino2002,Lederman2003,Otto2001,Denzler2005,McCann2006}, that, under suitable assumptions, there exist (sharp) rates of convergence in different topologies, the most common being the $\rd_1$ (see section~\ref{sec.dyn.sys.intro}). This results were proven by means of relative entropy functionals introduced in~\cite{Newman1984,Ralston1984},  or by means of the so called Bakry-\'Emery method,~\cite{Bakry1985}. The rate  $t^{-1}$ of uniform convergence in relative error in the whole range $\frac{d-2}{d}<m<1$,  has been computed first in \cite{Carrillo2003} for radial data, and later extended to a larger class of data in \cite{Kim2006}.  In a series of papers, similar results were obtained in the whole range $m< 1$, cf. \cite{Blanchet2007,Blanchet2009,Bonforte2010c,Bonforte2010}; notice that in the range $m <\frac{d-2}{d}$ there is a dramatic change in the behaviour of solutions since mass is not preserved and they can extinguish in finite time, see \cite{Vazquez2006,Vazquez2007}. In the general case $\gamma\neq0, \beta\neq0$, rates of convergence were  studied in~\cite{Bonforte2017,Bonforte2017a}.

In what follows we will show how we can combine the techniques of this paper with the ones used in \cite{Blanchet2009,Bonforte2010,Bonforte2010c}, to obtain rates of convergence to the Barenblatt profile with an (almost) uniform rate in the whole $\mathcal{X}$. For reasons that are not entirely clear up to now, we need to restrict ourselves to the range $\frac{d-1}{d}=m_1<m<1$ in the case $\gamma=\beta=0$, and to the range $\frac{2d-2-\beta-\gamma}{2(d-\gamma)}<m<1$ for the general case, see~\cite{Bonforte2017,Bonforte2017a} for further remarks. The latter restriction is somehow natural, since, at least when $\gamma=\beta=0$, we have that the FDE is a gradient flow of a displacement convex functional (the relative entropy) with respect to the so-called Wesserstein distance, see~\cite{McCann1997,Otto2001,McCann2006}. The displacement convexity is lost below $m_1$. The main result reads:
\begin{thm}[Almost Optimal Rates of Convergence in the non-weighted case]\label{almost.universal.rates.noweights}
Let $u$~be  the solution to \eqref{cauchy.problem} corresponding to the initial data $0\le u_0\in \mathcal{X}\setminus\{0\}$, $\int_{\RR^d}xu_0(x)\dx=0$ and assume that $\beta=\gamma=0$ and $m\in \big(\frac{d-1}{d},1\big)$. Then, for every $\delta\in (0,1)$  there exist $t_\delta, c_\delta>0$ (that may also depend on $u_0$) such that for all $t>t_\delta$
\begin{equation}\label{rates.noweight.infty}
\|u(t)-\mathfrak{B}(t;M)\|_{\LL^1(\RR^d)}\le \frac{c_\delta}{t^{1-\delta}}\qquad\mbox{and}\qquad
t^{d\vartheta}\|u(t)-\mathfrak{B}(t;M)\|_{\LL^\infty(\RR^d)}\le \frac{c_\delta}{t^{1-\delta}}\,,
\end{equation}
where $M=\|u_0\|_{\LL^1(\RR^d)}$.
\end{thm}

\noindent\textbf{Remark. }Notice that the above result new for the whole space $\mathcal{X}$ even if we are dealing with the case $\gamma=\beta=0$. Indeed, all the previous results deal with more restrictive assumption as radial data, a very precise control for $|x|\rightarrow\infty$ or being sandwiched between two Barenblatt profiles.

When dealing with CKN-weights, the result is a bit weaker, because of the possible lack of $C^k$ regularity at the origin and reads:
\begin{thm}[Minimal Rates of Convergence in the weighted case]\label{convergence.rates.weights}
 Assume $\gamma<0$ and let $u$ be  the solution to \eqref{cauchy.problem} corresponding to the initial data $0\le u_0\in \mathcal{X}\setminus\{0\}$ and let $\frac{2d-2-\beta-\gamma}{2(d-\gamma)}<m<1$. Then, there exists a $\delta_*\in (0,1)$ such that for every $\delta\in (0,\delta_*)$  there exist $t_\delta, c_\delta>0$ (that may also depend on $u_0$) such that for all $t>t_\delta$
\begin{equation}\label{rates.weight.infty}
\|u(t)-\mathfrak{B}(t;M)\|_{\LL^1_\gamma(\RR^d)}\le \frac{c_\delta}{t^{1-\delta}}\qquad\mbox{and}\qquad t^{(d-\gamma)\vartheta}\|u(t)-\mathfrak{B}(t;M)\|_{\LL^\infty(\RR^d)}\le \frac{c_\delta}{t^{1-\delta}}\,,
\end{equation}
where $M=\|u_0\|_{\LL^1_\gamma(\RR^d)}$.
\end{thm}

If we consider radial initial data in $\mathcal{X}$ we can provide a \emph{universal} rate of convergence, very much in the spirit of \cite{Carrillo2003} or \cite{Kim2006}.
\begin{thm}[Sharp Universal Rates for Radial Data]\label{rates.radial.data}
Let $u$ be  the solution to \eqref{cauchy.problem} corresponding to the radial initial data $0\le u_0\in \mathcal{X}\setminus\{0\}$ and assume that $\gamma=\beta=0$ and $m \in \left(\frac{d-2}{d}, 1\right)$. Then, there exist $t_0, c_0>0$ (that may also depend on $u_0$) such that for all $t>t_0$
\begin{equation}\label{rates.weight.infty.RE}
\left\|\frac{u(t)}{\mathfrak{B}(t;M)}-1\right\|_{\LL^\infty(\RR^d)}\le \frac{c_0}{t}\,,
\end{equation}
where $M=\|u_0\|_{\LL^1(\RR^d)}$.
\end{thm}
\begin{rem}\rm As an immediate consequence of \eqref{rates.weight.infty.RE} we obtain that for all $t\ge t_0$
\begin{equation}\label{rates.weight.infty.RAD}
\|u(t)-\mathfrak{B}(t;M)\|_{\LL^1(\RR^d)}\le \frac{c_0}{t}\qquad\mbox{and}\qquad t^{(d-\gamma)\vartheta}\|u(t)-\mathfrak{B}(t;M)\|_{\LL^\infty(\RR^d)}\le \frac{c_0}{t}.
\end{equation}
The above Theorem solves a problem left open in \cite{Carrillo2003}, i.e. identifying the largest class of nonnegative \emph{radial} $\LL^1$ data for which the above rate of convergence holds. Such rates are proven to be sharp, since they are fulfilled by two time-shifted Barenblatt, with the same mass, see \cite{Carrillo2003, Kim2006}.
Finally, we observe that, even if we restrict the analysis to radial data, the class $\mathcal{X}$ is much larger than those considered up to know in the literature: we refer to Section~\ref{sec:on.X} for examples of functions in $\mathcal{X}$ with a substantially different behaviour from  the Barenblatt profile.
\end{rem}

Finally, let us give the proof of the above statements.\\

\noindent\textbf{Proof of Theorem~\ref{almost.universal.rates.noweights}.} Here we exploit the techniques introduced in~\cite{Blanchet2009,Bonforte2010,Bonforte2010c}. Let us rescale $u(t,x)$ to $v(\tau, y)$ according to the change of variables~\eqref{change.fokker.planck} and define $w:=\frac{v(\tau, y)}{\mathcal{B}_M(y)}$ where $M=\|u_0\|_{\LL^1(\RR^d)}$. Let us define the \emph{Free Energy} or \emph{Relative Entropy} $\mathcal{F}[w]$ and the \emph{Fisher Information} $\mathcal{I}[w]$ as
\begin{equation}\label{entropy.fisher.def}\begin{split}
\mathcal F[w(\tau)] &:=\frac{m}{m-1}\, \int_{\RR^d}\left[\frac{w^m-1}{m}-(w-1)\right]\mathcal{B}_M^m \dy\,,\\
\mathcal I[w]&:=\frac{m}{1-m}\,\int_{\RR^d} w\,\mathcal{B}_M\Big|\nabla \left[\left(w^{m-1}-1\right)\mathcal{B}_M^{m-1}\right]\Big|^2\dy\,.
\end{split}\end{equation}
The Fisher information is related to the relative entropy by the time derivative along the flow
\begin{equation}\label{entropy.production}
\frac{\rm{d}}{\dtau}  \mathcal{F}[w] = - \mathcal{I}[w]\,.
\end{equation}
It is well known that the relative entropy controls the $\LL^1$ distance between the solution $v(\tau, y)$ and the Barenblatt profile $\mathcal{B}_M$, via the celebrated Csisz\'ar-Kullback inequality, see e.g. ~\cite{Kullback1967,Csiszar1967,Dolbeault2013,Carrillo2003, Otto2001},
more precisely
\begin{equation}\label{L1.entropy.inq}
\|v(\tau)-\mathcal{B}_M\|_{\LL^1(\RR^d)} \le \left(\frac{8}{m}\,\|\mathcal{B}_M^{2-m}\|_{\LL^1(\RR^d)}\right)^\frac{1}{2}\, \sqrt{\mathcal{F}[w]}\,.
\end{equation}
Therefore the decay of the relative entropy implies the same decay for of $\|u(t)-\mathcal{B}_M(t)\|_{\LL^1(\RR^d)}$.

For any $m \in \left(\frac{d-1}{d}, 1\right)$ the \emph{Entropy-Entropy Production} inequality reads
\begin{equation}\label{entropy.entropyproduction.in}
4\,\mathcal{F}[w]\le \mathcal{I}[w]\,,
\end{equation}
and it is well known to be equivalent to a member of a suitable family of (optimal) Gagliardo-Nirenberg inequalities, more details can be found in the pioneering work of Del Pino and Dolbeault \cite{DelPino2002}.
The best constant in~\eqref{entropy.entropyproduction.in} is $4$, and combining it with~\eqref{entropy.production}, we obtain the (sharp) exponential decay of the relative entropy along the flow, namely
\[
\mathcal{F}[w] \le \mathcal{F}[w_0]\,e^{-4\tau}\,.
\]
The strategy in~\cite{Blanchet2009,Bonforte2010,Bonforte2010c} consists in proving a faster decay of the entropy along the flow using an improved (with a larger constant) entropy-entropy production inequality  along the flow. Such improved inequality is obtained by means of Hardy-Poincar\'e type inequalities (with  improved constants), and by means of quantitative inequalities that compare the linear(ized) entropy and Fisher information with their nonlinear counterparts. Notice that we are in position to apply the results of  \cite{Blanchet2009, Bonforte2010}, since the running assumption guarantee that GHP (Theorem~\ref{thm.ghp}) holds, and implies the validity of assumption $(H1)''$ in those papers. Combining Lemma 3, Theorem 7 of \cite{Blanchet2009} with Lemma 1 of~\cite{Bonforte2010c}, we can prove the following claim.

\noindent\textit{Claim.} For any $0< \delta < 4\frac{1-\vartheta}{\vartheta}$ there exists a time $\tau_\delta>0$ such that
\begin{equation}\label{improved.entropy.fisher.inq}
\left(\frac{4}{\vartheta}-\delta\right)\,\mathcal{F}[w(\tau)]\le \mathcal{I}[w(\tau)]\,,\qquad\mbox{for any}\,\,\tau\ge \tau_\delta
\end{equation}
where $\vartheta$ is as in~\eqref{parameters.sigma.theta}.

\noindent {\sc Sketch of the proof of the claim.~}We shall not provide the lengthy details of the proof of the above  claim, we will just explain how to deduce it as a straightforward combination of already published results, adapting them to the current notations. The claim follows by formula (11) of~\cite{Bonforte2010c}, that in the current notations takes the form (at least for sufficiently large times)
\begin{equation}\label{PNAS-gap}
\frac{2\left[ \Lambda_{\alpha,d}-d(1-m)\left((1+\varepsilon)^{4(2-m)}-1\right)\right]}{(1+\varepsilon)^{7-3m}}\mathcal{F}[w(\tau)]\le \mathcal{I}[w(\tau)]\,,
\end{equation}
where $\varepsilon$ is (roughly speaking) the size of the relative error $|w-1|\sim\varepsilon$, which we need to be small in order to guarantee the validity of the result  (note that in formula (11) of~\cite{Bonforte2010c} $h=\max\{\sup\limits_{\RR^d}|w|, 1-\inf\limits_{\RR^d}|w|\}\sim 1+\varepsilon$).    Notice that everything is quantified explicitly in terms of $\varepsilon$ in the paper \cite{Bonforte2010c} which also relies on precise results of \cite{Blanchet2009, Bonforte2010}.  The smallness of $\varepsilon$ for sufficiently large times  follows by our Theorem~\ref{thm.ghp}, Global Harnack Principle, together with the uniform convergence in relative norm, Theorem~\ref{global.convergence.relative.error.thm}. Recalling now Lemma 1 of~\cite{Bonforte2010}, we get the expression for $\Lambda_{\alpha,d}= -4\alpha-2d$, which in our notations becomes $\Lambda_{\alpha,d}= \frac{2}{\vartheta(1-m)}$. Note that we need to assume that the first moment is fixed, but this is well-known to be true along the nonlinear flow as well, see \cite{Bonforte2010}. This concludes the proof of the claim.

\medskip

As a consequence of inequality~\eqref{improved.entropy.fisher.inq}, we obtain a faster decay of the relative entropy and conclude that
\begin{equation}\label{decay.l1}
\|v(\tau)-\mathcal{B}_M\|_{\LL^1(\RR^d)} \le C_\delta\, e^{-\left(\frac{2}{\vartheta}-\frac{\delta}{2}\right)\tau}\,,\qquad\mbox{for any}\,\,\tau\ge \tau_\delta.
\end{equation}
By re-scaling back to original variables and observing that $e^{2\tau}=R(t)\sim t^\vartheta$ one concludes that
\[
\|u(t, \cdot)-\mathfrak{B}(t, \cdot\,; M)\|_{\LL^1(\RR^d)} \le C_\delta\, t^{-1+\frac{\delta}{4}\vartheta}\,,
\]
Notice that $\delta>0$ is arbitrary, but $C_\delta$ may diverge as $\delta\to 0^+$.
This proves the left inequality in~\eqref{rates.noweight.infty}.

It only remains to prove the second inequality in~\eqref{rates.noweight.infty}, to do so we need to invoke the following interpolation Lemma which goes back to Gagliardo (see~\cite{Gagliardo1958}) and Nirenberg (see~\cite[Pag. 126]{Nirenberg1959}): let $f \in C^k(\RR^d)\cap\LL^1(\RR^d)$ for some $p\ge1$ and $k$ a positive integer, then
\begin{equation}\label{interpolation.inequality}
\|f\|_{\LL^\infty(\RR^d)}\le C_{p,k,d} \, \|f\|_{C^k(\RR^d)}^\frac{d}{d+k}\,\|f\|_{\LL^1(\RR^d)}^\frac{k}{d+k}\,,
\end{equation}
where $\|\cdot\|_{C^k(\RR^d)}$ is given by
\[
\|f\|_{C^k(\RR^d)}:=\max\limits_{|\eta|=k}\,\sup_{z\in\RR^d}\Big|\partial^\eta f(z)\Big|\,,
\]
where $|\eta|=\eta_1+\ldots+\eta_d$ is the length of the multi-index $\eta=\left(\eta_1, \ldots, \eta_d\right) \in \mathbb{Z}^d$. We recall that in the case $\gamma=\beta=0$ solution to~\eqref{cauchy.problem} are $C^\infty(\RR^d)$ and for any $k\ge1$ we have that
\[
\sup_{\tau\ge\tau_0}\|v(\tau)-\mathcal{B}_M\|_{C^j(\RR^d)}<\infty\,,
\]
for a proof of the above inequality see~\cite[Theorem 2 and Theorem 4]{Blanchet2009}.  Fix $k\ge1$ to be chosen later, combining the above interpolation inequality~\eqref{interpolation.inequality} with the decay of the $\LL^1$ norm given in~\eqref{decay.l1} one obtains
\[
\|v(\tau)-\mathcal{B}_M\|_{\LL^\infty} \lesssim e^{-\left(\frac{2}{\vartheta}-\frac{\delta}{2}\right)\left(\frac{k}{d+k}\right)}\,,
\]
and rescaling back to original variables we easily find that
\[
t^{d\vartheta}\|u(t,\cdot)-\mathfrak{B}(t, \cdot\,; M)\|_{\LL^\infty} \le C\, t^{-\left(1 -\frac{\delta\vartheta}{4}-\frac{d}{d+k}+\frac{d\,\delta\,\vartheta}{4(d+k)}\right)}\,,
\]
since both $k$ and $\delta$ we arbitrary we conclude that the second inequality in~\eqref{rates.noweight.infty} holds by choosing $k$ sufficiently large, for instance $k\ge \delta^{-4}$. This concludes the proof. \qed
\medskip
\noindent\textbf{Proof of Theorem~\ref{convergence.rates.weights}}. The proof is very similar to the one of Theorem~\ref{almost.universal.rates.noweights}, here we only explain the main differences.  We cannot reach the rate $t^{-1+\delta}$ for two reasons. The first: we can obtain an inequality as the one~\eqref{improved.entropy.fisher.inq}, however the constant is smaller that $\left(\frac{4}{\vartheta}-\delta\right)$, see~\cite{Bonforte2017,Bonforte2017a}. The second: we need to assume $\gamma<0$ to obtain an inequality similar to~\eqref{interpolation.inequality}, see Lemma~\ref{interpolation.lemma}. Finally, solutions to~\eqref{cauchy.problem} do not enjoy $C^\infty$ regularity, indeed they can be only $C^\alpha$ at the origin, cf.~\cite{Bonforte2019a}, hence inequality~\eqref{interpolation.inequality} does not apply as in the proof of Theorem~\ref{almost.universal.rates.noweights}: we can only interpolate with $C^\alpha$ norms, and we can not choose $k$ arbitrarily large. This concludes our considerations. \qed

\medskip
\noindent\textbf{Proof of Theorem~\ref{rates.radial.data}. }In \cite{Carrillo2003}, the authors proved Theorem~\ref{rates.radial.data} under the assumption that the initial data $u_0$ is bounded, radially symmetric and satisfies $u_0=O(|x|^{-\frac{2}{1-m}})$. It is only needed to show that radial data in $\mathcal{X}$  produce solutions that satisfy the decay assumption above for any time $t>t_0$ for some given $t_0$. This is exactly the statement of the GHP, Theorem~\ref{thm.upper}.  The proof is concluded.\qed

\section{Counterexamples and Generalized Global Harnack Principle}\label{Sec4}
In this section we carefully construct the family of sub/super solutions presented in the introduction, such a phenomena is possible only in $\mathcal{X}^c$. We show also examples of an anomalous ``fat-tail'' behaviour for both integrable and not-integrable solutions. At the end of this section we show that in $\mathcal{X}^c$ the convergence toward the Barenblatt is slower due to a different tail behaviour of solutions.

\subsection{Construction of a family of Subsolutions and anomalous tail behaviour}\label{ssec.subsols}
In the following Proposition we construct an explicit family of sub-solutions parameterized  by the powers of their  decay at infinity. Every subsolution decay in space slowly then the Barenblatt profile.
\begin{prop}[Family of $\LL^1_\gamma$-Subsolutions]\label{subsolution.prop}
Let $m \in (m_c, 1)$,  $\varepsilon \in (0, \frac{2}{1-m}-\frac{2}{\sigma}(d-\gamma))$,   $A, B >0$ and
\[
\alpha=\frac{1}{1-m} - \frac{\varepsilon}{2}>0.
\]
Define for some $t_0 \in \RR$ the function
\begin{equation}\label{definition.D}
D(t):= \left(\sigma\,A^{m-1}\,m\,B\,(d-\gamma)\, (1-\alpha(1-m)) \,t + t_0\right)^{\frac{1}{1-\alpha(1-m)}}.
\end{equation}
Then, for all $t>0$, the $\LL^1_\gamma(\RR^d)$ function
\begin{equation}\label{subsolution}
\underline{V}(t,x) = \frac{A}{(D(t)+B |x|^\sigma)^{\alpha}}
\end{equation}
is a subsolution to \eqref{cauchy.problem}.
If $m \in (\frac{d-\gamma}{(d-\gamma+\sigma)}, 1)$ and $\varepsilon \in (0, \frac{2}{1-m}-\frac{2(d-\gamma)}{\sigma}-2)$ then $|x|^\sigma\,\underline{V}(t, x) \in \LL^1_\gamma(\RR^d)$.
\end{prop}
\begin{rem}\rm
We notice that $\|V(t,\cdot)\|_{\LL^\infty(\RR^d)} \asymp  t^{-\frac{\alpha}{1-\alpha(1-m)}}$ as $t\rightarrow\infty$. This is not in contrast with the smoothing effect (inequality~\eqref{global.upper.estimate}) which implies that any solution $u(t,x)$ to~\eqref{cauchy.problem} decays in time less than $t^{-(d-\gamma)\vartheta}$: a simple computation shows that the condition~$\varepsilon \in \left(0,  \frac{2}{1-m}-2\frac{d-\gamma}{\sigma}\right)$ implies that $t^{-\frac{\alpha}{1-\alpha(1-m)}}~<~t^{-(d-\gamma)\vartheta}$.\\ However, as $|x|\rightarrow \infty$, $V(t, x)$ exhibits quite an interesting behaviour, namely $V(t, x)~\asymp~|x|^{-\sigma\alpha}$. The power $-\sigma\alpha$ do not match the one of the fundamental solution: indeed, we have $\mathfrak{B}(t,x; M)\asymp|x|^{-\frac{\sigma}{1-m}}$,   as $|x|\rightarrow\infty$. This proves that, for any choice of the parameters $A, B, t_0$ and for any choice of the mass $M$, the inequality $V(t,x) > \mathfrak{B}(t,x; M)$ holds for $|x|$ large enough. \\
As a final remark, we can define another family of subsolution. Indeed, for some choice of the parameters $ B', F'$ and $T$ the function $W$ defined as
\[
W(t, x)=\frac{\left(T-t\right)^{\frac{1}{1-m}}}{(B'+F' |x|^\sigma)^\alpha}\,,
\]
is a subsolution to~\eqref{cauchy.problem} which has the same qualitative behaviour as $|x|\rightarrow\infty$, the drawback is that this is meaningful only on a finite time interval, hence we prefer to use $\underline{V}$.
\end{rem}
\noindent{\bf Proof of Proposition \ref{subsolution.prop}: } We just need to verify that the function $\underline{V}(t,x)$ defined in~\eqref{subsolution} satisfies in a pointwise sense, the inequality
\begin{equation}\label{condition.subsolution}
\partial_t \underline{V}(t, x) \leq |x|^\gamma\mathrm{div}\left(|x|^{-\beta}\nabla \underline{V}^m\right)\,.
\end{equation}
Let $r=|x|$, and write -with a little abuse of notation- $\underline{V}(t,r)$ instead of $\underline{V}(t,x)$. We recall that the operator $\mathcal L_{\gamma, \beta} = |x|^{\gamma}\nabla \cdot \left(|x|^{-\beta}\nabla f \right)$ acts on a radial function $f(r)$ in the following way
\begin{equation}\label{operator.same.weights.radial}
\op(f)= r^{\gamma-\beta}\left(f''(r) + \frac{\left(d-1-\beta\right)}{r}\,f'(r)\right)\,.
\end{equation}
A straightforward computation shows the following identities:
\begin{equation*}\begin{split}
\partial_t \underline{V}(t, r) &= \frac{- A\,\alpha\,\partial_t D(t)}{\left(D(t)+F r^\sigma\right)^{\alpha+1}}\,,   \\
\op\left( \underline{V}(t, r)^m\right) &= \frac{-(\sigma\,\alpha m\, A^m\, F)}{(D(t)+F r^\sigma)^{\alpha m+2}}\left[(d-\gamma)D(t)+F\,r^\sigma\,(-\sigma\,\alpha\,m +d-2-\beta)\right]\,.
\end{split}\end{equation*}
As a consequence, the inequality $\partial_t V(t, r) \leq \op \left(V(t,r)^m\right)$ is satisfied if and only if
\begin{align}\label{inequality.subsolution}
\partial_t D(t) \geq \frac{\sigma\,m\,F\,A^{m-1}}{(D(t)+F r^\sigma)^{\alpha (m-1)+1}}  \left[(d-\gamma)D(t)+F\,r^\sigma\,(-\sigma\,\alpha\,m +d-2-\beta)\right]\,,
\end{align}
The reader may notice that if $\varepsilon < 2/(1-m)- 2\frac{(d-\gamma)}{\sigma}$ then in the right-hand-side of inequality~\eqref{inequality.subsolution} the term $F\,r^\sigma\,(-\sigma\,\alpha\,m +d-2-\beta)$ is negative. A simple computations then shows that the supremum of the right-hand-side of inequality~\eqref{inequality.subsolution} is achieved at $r=0$. Hence inequality $\partial_t V(t, r) \leq \op \left(V(t,r)^m\right)$ will follow by asking that
\begin{align*}
\partial_t D(t) &\geq \sigma\,m\,F\,A^{m-1}\,(d-\gamma)\, D(t)^{\alpha\,(1-m)}  \\
&= \sup_{r\ge0}\frac{\sigma\,m\,F\,A^{m-1}}{(D(t)+F r^\sigma)^{\alpha (m-1)+1}}  \left[(d-\gamma)D(t)+F\,r^\sigma\,(-\sigma\,\alpha\,m +d-2-\beta)\right]\,.
\end{align*}
We conclude the proof observing that, for any $t_0\in\RR^d$, such an inequality is satisfied by the function $D(t)$ defined in~\eqref{definition.D}. \qed

\noindent\textbf{Anomalous tail Behaviour. }As a corollary of Proposition~\ref{subsolution.prop} we have the following results about unexpected ``fat-tails'', both integrable and not-integrable.
\begin{cor}[Anomalous Integrable Tail Behaviour]\label{corollary.subsolution.1}
Under the assumptions of Proposition~\ref{subsolution.prop}:
\[
\mbox{If }\qquad u_0(x)\ge \frac{A}{\left(C+B|x|^\sigma\right)^\alpha},\qquad\mbox{then for any $t>0$ we have }\qquad\liminf_{|x|\rightarrow \infty} |x|^{\sigma\,\alpha}\,u(t,x) \ge \frac{A}{B}\,.
\]
Moreover, for any $t>0$ we have
\begin{equation}\label{relative.error.infinite}
\sup_{x\in\RR^d} \left| \frac{u(t,x)}{\mathfrak{B}(t,x; M)}-1 \right|=\infty \,.
\end{equation}
\end{cor}
\noindent\textbf{Proof.} The proof of the first statement of the above Corollary is an immediate application of Proposition~\ref{subsolution.prop}. As for formula \eqref{relative.error.infinite}, we see that if the initial data $u_0$ do not satisfy the tail condition~\eqref{tail.condition} there are no chances to conclude the convergence to the Barenblatt profile in \emph{uniform relative error}. Indeed, since
\[
u_0(x)\ge \frac{A}{\left(C+B|x|^\sigma\right)^\alpha}\,,\qquad\mbox{by Proposition~\ref{subsolution.prop} we get that}\qquad u(t,x) \ge \frac{A}{(D(t)+B |x|^\sigma)^{\alpha}}\,,
\]
for any $t>0$, where $D(t)=\left(\sigma\,A^{m-1}\,m\,B\,(d-\gamma)\, (1-\alpha(1-m)) \,t + C^{1-\alpha(1-m)}\right)^{\frac{1}{1-\alpha(1-m)}}$. A simple computation shows that the quotient
\[
\frac{1}{\mathfrak{B}(t,x; M)}\, \frac{A}{(D(t)+B |x|^\sigma)^{\alpha}}= \frac{A}{t^\frac{1}{1-m}}\,\frac{\left[b_0\frac{t^{\sigma\vartheta}}{M^{\sigma\vartheta(1-m)}}+b_1 |x|^{\sigma}\right]^\frac{1}{1-m}}{(D(t)+B |x|^\sigma)^{\frac{1}{1-m}-\frac{\varepsilon}{1}}} \sim \frac{A\,b_1^\frac{1}{1-m}}{B^{\frac{1}{1-m}-\frac{\varepsilon}{2}}} \, \frac{|x|^\frac{\sigma\varepsilon}{2}}{t^\frac{1}{1-m}}\quad\mbox{as}\quad |x|\rightarrow \infty\,,
\]
from which we deduce~\eqref{relative.error.infinite}. The proof is then complete. \qed

\noindent\textbf{A Family of Subsolutions not in $\LL^1_\gamma$. }A closer inspection of the proof of Proposition~\ref{subsolution.prop} reveals that the condition on $\varepsilon$ for $\underline{V}$ to be a subsolution is actually $m\,\varepsilon \le \frac{2}{1-m}-\frac{2}{\sigma}(d-\gamma)$, while for $\varepsilon > \frac{1}{m}\left(\frac{2}{1-m}-\frac{2}{\sigma}(d-\gamma)\right)$, $\underline{V}$ ceases to be a subsolution.
Indeed, when $\varepsilon > \frac{2}{1-m}-\frac{2}{\sigma}(d-\gamma)$, we can construct subsolutions which do not belong to $\LL^1_\gamma(\RR^d)$, so that, in general, initial data $u_0\not\in\LL^1_\gamma(\RR^d)$ do not produce $\LL^1_\gamma(\RR^d)$ solutions for any time $t>0$.  We resume this fact in the following corollary, considering the parameters $m, \varepsilon, A,B, t_0$ in the ``non integrability range'':
\begin{equation}\label{non-L1.hyp}
m \in (m_c, 1)\,\,, \varepsilon \in \left[\frac{2}{1-m}-\frac{2}{\sigma}(d-\gamma), \frac{1}{m}\left(\frac{2}{1-m}-\frac{2}{\sigma}(d-\gamma)\right)\right]\,\,,   A, B, t_0 >0\,,\quad  \alpha=\frac{1}{1-m} - \frac{\varepsilon}{2}>0\,.
\end{equation}
\begin{cor}[Family of Subsolutions not in $\LL^1_\gamma$]\label{not.integrable.subsolution}
Under assumptions \eqref{non-L1.hyp}, the function $\underline{V}(t,x)$ defined in~\eqref{subsolution} %
is a non-integrable sub-solution to~\eqref{cauchy.problem}. If $u(t,x)$ is a (super)solution to~\eqref{cauchy.problem}, we then have
\[
u_0\ge \underline{V}(0,\cdot)\not\in\LL^1_\gamma(\RR^d)\qquad\mbox{then }\qquad u(t, \cdot) \ge \underline{V}(t, \cdot) \not\in\LL^1_\gamma(\RR^d).
\]
\end{cor}

\subsection{Construction of a family of Supersolutions}
In this section we construct a family of supersolutions which share the same spatial tail behaviour with the subsolutions constructed in the previous Section in Proposition~\ref{subsolution.prop}.
\begin{prop}[Family of $\LL^1_\gamma$-Supersolutions]\label{supersolution.prop}
Let $m \in (m_c, 1)$,  $\varepsilon \in (0, \frac{2}{1-m}-\frac{2}{\sigma}(d-\gamma))$,  $E, F >0$ and $\alpha=\frac{1}{1-m} - \frac{\varepsilon}{2}>0$.
Define for some $t_0 \in \RR$ and $H>0$ the function
\begin{equation}\label{supersolution}
\overline{V}(t,x) = \frac{E\,G(t)^\alpha}{(G(t)+F |x|^\sigma)^{\alpha}} \in \LL^1_\gamma(\RR^d)\,, \qquad\mbox{where}\qquad G(t):=t_0 + H\,t\,.
\end{equation}
Then, $\overline{V}$ is a supersolution for all $t>0$ whenever $H$ is sufficiently big, more precisely for all
\begin{equation}\label{condition.B}
H \ge m\,\sigma\,F^2\,E^{m-1}\,\left(2+\beta-d+\sigma\,\alpha\,m\right)\,.
\end{equation}
\end{prop}
\noindent{\bf Proof of Proposition \ref{supersolution.prop}: }  We just need to verify that the function $\overline{V}(t,x)$ defined in~\eqref{supersolution} satisfies inequality
\begin{equation}\label{condition.supersolution}
\partial_t \overline{V}(t,x) \ge |x|^\gamma\mathrm{div}\left(|x|^{-\beta}\nabla \overline{V}^m\right)(t,x)\,.
\end{equation}
under the assumption~\eqref{condition.B}. Let $r=|x|$, and write -with a little abuse of notation- $\overline{V}(t,r)$ instead of $\overline{V}(t,x)$. We have the following identities (recalling the radial form of  $\mathcal L_{\gamma, \beta}$,  formula~\eqref{operator.same.weights.radial}.)
\begin{align*}
\partial_t \overline{V}(t, r) &=\frac{\alpha\, E\, G(t)^{\alpha-1}\, H}{\left(G(t)+ F\,r^\sigma\right)^{\alpha+1}} \,F\,r^\sigma\,. \\
\op\left( \overline{V}^m(t, r )\right) &= \frac{\sigma\,\alpha\,m\,F\,E^m G(t)^{\alpha\,m}}{(G(t)+F\,r^\sigma)^{\alpha m+2}} \left[F\,r^\sigma(2+\beta-d+\sigma\,\alpha\,m)-(d-\gamma)G(t)\right]
\end{align*}
It is straightforward to verify that~\eqref{condition.supersolution} holds at $r=0$ since for any $t>0$ the derivative in time $\partial_t \overline{V}(t, 0)=0$ and $\op\left( \overline{V}^m(t, 0)\right)$ is negative. When $r>0$ a direct computation shows that~\eqref{condition.supersolution} is equivalent to
\begin{equation}\label{condition.B.2}
 H\ge \left(\frac{G(t)}{G(t)+F\,r^\sigma}\right)^{1-\alpha(1-m)}\, m\,\sigma\,F\,E^{m-1}\,\left[F\,(2+\beta-d+\sigma\,\alpha\,m)-(d-\gamma)\frac{G(t)}{r^\sigma}\right]\,.
\end{equation}
Finally, we check that \eqref{condition.B} implies~\eqref{condition.B.2}, just by using that $\left(\frac{G(t)}{G(t)+F\,r^\sigma}\right)^{1-\alpha(1-m)}<1$ and $(d-\gamma)\frac{G(t)}{r^\sigma}>0$. Therefore, $\overline{V}(t,x)$ is a supersolution and the proof is concluded. \qed
\medskip
A closer look at the above proof, reveals that if we allow $\varepsilon < 0$, then $\overline{V}(t,x)$ ceases to be a supersolution. Indeed, if $\varepsilon <0$ we have that $1-\alpha(1-m)<0$ and so in~\eqref{condition.B.2} we would have that
\[
\left(\frac{G(t)}{G(t)+F\,r^\sigma}\right)^{1-\alpha(1-m)}= \left(1+\frac{F\,r^\sigma}{G(t)}\right)^{\alpha(1-m)-1}\rightarrow \infty\quad\mbox{as}\quad r\rightarrow \infty\,,
\]
and an inequality as~\eqref{condition.B.2} would be impossible.


\subsection{Slower convergence rates in $\mathcal{X}^c$}\label{ssec:no.rates}%
We have shown in Section \ref{SSec.ratesX}, that when $u_0\in\mathcal{X}$, then there are always power-type rates of convergence to a Barenblatt profile and in some cases the sharp decay rate $O(1/t)$ is obtained. In this paragraph we show, by means of an explicit counterexample, that power-like decay rates are simply not possible for general data outside $\mathcal{X}$. However, we are not able to exclude the possibility of slower decay rates (e.g. $\log, \log\log,$ etc.). The latter question is really delicate and deserves a thorough study that goes beyond the scope of this paper.
\begin{thm}\label{no.rates}
For any $\delta>0$ there exists initial data $u_{0,\delta} \in \LL^1_\gamma(\RR^d)$ such that the corresponding solution $u_\delta(t,x)$ to~\eqref{cauchy.problem} satisfies
\begin{equation}\label{linfinity.limit}
\lim\limits_{t\rightarrow\infty}t^{\delta+(d-\gamma)\vartheta}\|u_\delta(t)-\mathfrak{B}(t;M)\|_{\LL^\infty(\RR^d)}=\infty\,,
\qquad\mbox{and}\qquad
\lim\limits_{t\rightarrow\infty}t^{\delta}\|u_\delta(t)-\mathfrak{B}(t;M)\|_{\LL^1_\gamma(\RR^d)}=\infty.
\end{equation}
\end{thm}

\noindent\textbf{Proof. }We only give a detailed proof of the left limit in~\eqref{linfinity.limit} in the non-weighted case $\gamma=\beta=0$, the weighted case being completely analogous. For the same reasons, we will just sketch the proof of the right-limit in the non-weighted case.

\noindent\textit{Proof of the left-limit. }Fix any $\delta>0$ and let $\varepsilon \in \left(0, \frac{2}{1-m}-d\right)$ be such that
\begin{equation}\label{choice.delta}
\delta> \frac{2\vartheta}{\varepsilon(1-m)}\left(\frac{2}{1-m}-d-\varepsilon\right)\,.
\end{equation}
A simple computation shows that this choice of $\varepsilon$ is always possible. Let $u_\delta(t,x)$ be the solution to \eqref{cauchy.problem} corresponding to the initial data
\[
u_{0, \delta}(x)=\frac{A}{\left(1+B|x|^2\right)^\alpha}\,,\qquad\mbox{where $\alpha=\frac{1}{1-m}-\frac{\varepsilon}{2}$  and  $A, B >0$ are chosen such that $\|u_{0, \delta}\|_{\LL^1(\RR^d)}=1$.}
\]
A simple rescaling shows that for any choice of $B$, we can always choose $A$ such that $\|u_{0, \delta}\|_{\LL^1(\RR^d)}=1$, namely
\begin{equation}\label{condition.A.B}
A\,\int_{\RR^d}\frac{\dy}{\left(1+|y|^2\right)^\alpha}= B^\frac{d}{2}\,.
\end{equation}
Let us  consider the subsolution $\underline{V}$ given in Proposition~\eqref{subsolution.prop}
\[
\underline{V}(t,x) = \frac{A}{(D(t)+B |x|^2)^{\alpha}}\,,\quad\mbox{where}\quad\,D(t):= \left(2\,A^{m-1}\,m\,B\,d\, (1-\alpha(1-m)) \,t + 1\right)^{\frac{1}{1-\alpha(1-m)}}\,.
\]
We recall that $u_{0, \delta}(x)=\underline{V}(0, x)$ so that by comparison $u_\delta(t, x)\ge\underline{V}(t,x)$ for all $t\ge0$. The following claim allows to conclude the proof:

\noindent\textit{Claim}. For any $B$ sufficiently large, there exists $t_0=t_0(B)>0$ and $\underline{c}>0$ such that for any  for any $t\ge t_0$ and for any $|x|^2 \in [D(t), 2 D(t)]$, we have that
\begin{equation}\label{subsolution.Barenblatt}
\underline{V}(t,x)> \mathfrak{B}(t, x;1)\,\quad\mbox{and}\quad (\underline{V}(t, x)-\mathfrak{B}(t, x;1))\,\,>\,\,\underline{c}\,t^{- \frac{\alpha}{1-\alpha(1-m)}}\,.
\end{equation}
Let us assume momentarily the validity of the claim and conclude the proof. The above inequality~\eqref{subsolution.Barenblatt}  immediately implies that $u_\delta(t,x)\ge \mathfrak{B}(t, x;1)$ for any $|x|^2 \in [D(t), 2 D(t)]$. Let now $|x|^2=D(t)$, we then have
\begin{equation}\label{computation.linfinity}\begin{split}
t^{d\vartheta+\delta}\,\|u_\delta(t)-\mathfrak{B}(t; 1)\|_{\LL^\infty(\RR^d)} &\ge t^{d\vartheta+\delta}\, |u_\delta(t, x)-\mathfrak{B}(t, x; 1)|\\
&\ge t^{d\vartheta+\delta}\, (u_\delta(t, x)-\mathfrak{B}(t, x; 1)) \\
&\ge t^{d\vartheta+\delta}(\underline{V}(t, x)-\mathfrak{B}(t, x;1))\,\,\ge\,\,\underline{c}\,t^{d\vartheta+\delta - \frac{\alpha}{1-\alpha(1-m)}}
\end{split}\end{equation}
where in the last line we have used again~\eqref{subsolution.Barenblatt}. Notice that $d\vartheta+\delta - \frac{\alpha}{1-\alpha(1-m)}>0$, since $\delta$ is as in~\eqref{choice.delta}, therefore~\eqref{linfinity.limit} follows as a consequence of~\eqref{computation.linfinity}. The proof of the left-limit is complete.

\noindent\textit{Proof of the right-limit. }The proof follows by integrating on the region $|x|^2 \in [D(t), 2 D(t)]$, the last line of inequality \eqref{computation.linfinity} and recalling that $u_\delta(t, x)\ge\underline{V}(t,x)$ for all $t\ge0$.

It only remains to prove the Claim.

\noindent {\it Proof of the Claim. } Let us define $R(t)=\{x\in\RR^d: D(t)\le|x|^2\le 2D(t)\}$, for any $t>0$. To prove the first inequality in~\eqref{subsolution.Barenblatt} we need to check that for $t$ large enough we have
\[
\inf_{x\in R(t)}\underline{V}(t, x)> \sup_{x\in R(t)}\mathfrak{B}(t, x;1)\,,\qquad\mbox{which amounts to prove that}\qquad\underline{V}(t, 2D(t))> \mathfrak{B}(t, D(t); 1)\,,
\]
where in the last inequality we made a small abuse of language: we write $\mathfrak{B}$ as a radial function $\mathfrak{B}(t, r; 1)$.  We rewrite $\underline{V}(t, 2D(t))> \mathfrak{B}(t, D(t); 1)$ in the equivalent form
\begin{equation}\label{inequality.limit.1}
\frac{A}{(1+2B)^\alpha}\, \frac{D(t)^\frac{\varepsilon}{2}}{t^\frac{1}{1-m}} - \left[b_1 + b_0\,\frac{t^{2\vartheta}}{D(t)}\right]^\frac{-1}{1-m}\ge a >0\,,
\end{equation}
for some $a>0$. The proof of~\eqref{inequality.limit.1} follows by observing that in the limit $t\rightarrow\infty$, we have that for all $B$ large enough
\[
(1+2B)^\frac{\varepsilon}{2}\,\left(\frac{B}{1+2B}\right)^\frac{1}{1-m}\left(m\,d\,\varepsilon\,(1-m)\right)> b_1^{\frac{1}{m-1}}\,.
\]
Hence, inequality~\eqref{inequality.limit.1} holds for $B$ and $t$ large enough, since it is true in the limit $t\rightarrow \infty$ and all the quantities that appear in~\eqref{inequality.limit.1} are continuous with respect to $t,B>0$.

\noindent It only remains to prove the last inequality in~\eqref{subsolution.Barenblatt}:  for any $x \in R(t)$ we have that
\begin{equation*}\begin{split}
\underline{V}(t,x)-\mathfrak{B}(t, x; 1)& \ge \underline{V}(t, 2D(t)) - \mathfrak{B}(t, D(t); 1) \\
& = \left(\frac{t}{D(t)}\right)^\frac{1}{1-m}\, \left[\frac{A}{(1+2B)^\alpha}\,\frac{D(t)^\frac{\varepsilon}{2}}{t^\frac{1}{1-m}}-\frac{1}{\left(b_1 + b_0\,\frac{t^{2\vartheta}}{D(t)}\right)^\frac{1}{1-m}}\right]
  \ge \underline{c}\, t^\frac{-\alpha}{1-\alpha(1-m)}\,,
\end{split}\end{equation*}
where we have used that $t/D(t)\asymp t^\frac{-\alpha(1-m)}{1-\alpha(1-m)}$ and~\eqref{inequality.limit.1}. The proof of the claim and of the Theorem is now concluded.\qed

\section{On the Fast Diffusion Flow in $\mathcal{X}$}\label{sec:on.X}

In this section we analyze some properties of the tail space $\mathcal{X}$ that plays a key role in the proof of the upper estimates of Theorem~\ref{thm.upper} and it is the optimal space for \emph{uniform convergence in relative error}.
\subsection{An equivalent tail condition}\label{ssect5.1}
Here we analyze the tail condition~\eqref{tail.condition.1} introduced by Vazquez in~\cite{Vazquez2003}: we will prove that it is equivalent to~\eqref{tail.condition}, but this fact is not trivial: its proof needs the GHP of Theorem~\ref{thm.ghp}, as we shall see below.
\begin{prop}
Let $d\geq 3$, $\gamma, \beta < d$ real numbers such that  $\gamma -2 < \beta \leq \gamma(d-2)/d$ and $m \in \left(m_c, 1\right)$. Assume $f\in \LL^1_\gamma(\RR^d)$. Then,
\[
\mbox{$f$ satisfies~\eqref{tail.condition} if and only if it satisfies~\eqref{tail.condition.1}}
\]
\end{prop}
\noindent\textbf{Proof.} We will first prove that~\eqref{tail.condition} implies~\eqref{tail.condition.1}. Assume that $|f|_{\mathcal{X}}<\infty$ and let $x\in\RR^d$, $x\neq0$. We have the following chain of inequalities
\[
\int_{B_{\frac{|x|}{2}(x)}}|f(y)|\frac{\dy}{|y|^\gamma}\le \int_{B_{\frac{|x|}{2}(0)}^c}|f(y)|\frac{\dy}{|y|^\gamma} \le 2^{\frac{2+\beta-\gamma}{1-m}-(d-\gamma)}\,\,|f|_{\mathcal{X}} \,\, |x|^{d-\gamma-\frac{2+\beta-\gamma}{1-m}}\,=\mathrm{O}(|x|^{d-\gamma-\frac{2+\beta-\gamma}{1-m}})\,,
\]
which is exactly~\eqref{tail.condition.1}. In the above line we have used that $B_{\frac{|x|}{2}(x)}\subset B_{\frac{|x|}{2}(0)}^c$. Assume now that $f$ satisfies~\eqref{tail.condition.1}, without loss of generality we can assume that $f\neq0$. Let $u(t,x)$ be the solution to~\eqref{cauchy.problem} with initial data $u(0,x)=|f(x)|$. As proven by Vazquez in~\cite{Vazquez2003} and also by Vazquez and one of the authors in~\cite{Bonforte2006}, $u(t,x)$ satisfies inequality~\eqref{ghp.inequality}, i.e. the GHP. Therefore, by Theorem~\ref{thm.ghp} we have that $u_0=|f|\in\mathcal{X}\setminus\{0\}$, which means that $f$ satisfies~\eqref{tail.condition}. The proof is concluded. \qed
\subsection{A non-equivalent tail condition: an example of a ``bad'' functions in $\mathcal{X}$}\label{ssec.BadGuy}
In \cite{Bonforte2006} a (non sharp) sufficient condition for the validity of the GHP takes the form: there exists $R>0$ and $A>0 $ such that
\begin{equation}\label{point.wise.decay}
|f(x)|\le \frac{A}{|x|^\frac{2+\beta-\gamma}{1-m}}\quad\mbox{for any $|x|\ge R$}\,.
\end{equation}
It is easy to check that the above condition is sufficient to guarantee that $f\in\mathcal{X}$, but not necessary, as we shall explain by means of the following example: we construct a function $f\in \mathcal{X}$ which does not satisfy~\eqref{point.wise.decay}.

Let $d\ge3$ and, to fix ideas, let us assume that $\gamma=\beta=0$.  Define the function $f$ to be
\begin{equation*}
f(y)=\sum_{N=2}^\infty \frac{\chi_{B_{N^{-2}}(x_N)}(y)}{N^{\frac{2}{1-m}-1}}
\end{equation*}
where  $N\in [2,\infty)\cap\NN$, $x_N=(N,\mathbf{0})$, where $\mathbf{0} \in \RR^{d-1}$ is the zero vector.
The function $f$ is well defined, since by construction $B_{N^{-2}}(x_N)\cap B_{M^{-2}}(x_M)= \emptyset$ unless $N=M$, and $\frac{2}{1-m}-1>1$ in since $\frac{d-2}{d}<m<1$.   We then have that
\begin{equation*}\begin{split}
\int_{\RR^d} |x|^{\frac{2}{1-m}}f(x)\dx &=  \sum_{N=2}^\infty \int_{B_{\frac{1}{N^2}}(x_N)}|x|^\frac{2}{1-m}f(x)\dx \le c\,\sum_{N=2}^\infty \frac{N^{\frac{2}{1-m}-2d}}{N^{\frac{2}{1-m}-1}}\le c\,  \sum_{N=2}^\infty \frac{1}{N^{2d-1}} < \infty\,,
\end{split}\end{equation*}
where $0<c=c(m,d)<\infty$ and we have used the fact that if  $x \in B_{\frac{1}{N^2}}(x_N)$ then $|x|\le 2N$ and we recall that the last series converges since $d\ge 3$. As a consequence, we have that $f\in \mathcal{X}$: indeed for any $R>1$
\[
R^{\frac{2}{1-m}-d}\int_{B_R^c(0)}f\dx \le R^{\frac{2}{1-m}}\int_{B_R^c(0)}f\dx \le \int_{B_R^c(0)}|x|^\frac{2}{1-m}f\dx<\infty\,.
\]
On the other hand, a straightforward computation shows that $f$ does not satisfy the pointwise decay condition in~\eqref{point.wise.decay}. As   expected, the pointwise condition is more restrictive than the integral one.

To conclude, we give also an example of a \textit{radial function} $h\in \mathcal{X}$ which does not satisfy~\eqref{point.wise.decay}. Let
\begin{equation*}
h(y)=\sum_{N=2}^\infty \frac{\chi_{A_N(y)}}{|N-|y||^{\eta}}\,,
\end{equation*}
where $N\ge2$ is an integer,  $A_N:=\{x\in\RR^d: N\le|x|\le N+N^{-\alpha}\}$ with $0<\eta<1$ and $(1-\eta)\alpha>2/(1-m)$.

\subsection{The Fast Diffusion flow as a curve in $\mathcal{X}$}\label{X-flow} In this section we will consider solutions to~\eqref{cauchy.problem} as continuous curves in $\mathcal{X}$. To this end, we provide some details about the natural topology of $\mathcal{X}$. Indeed, the metric associated to the  ``natural'' norm $|\cdot|_\mathcal{X}$ would provide $\mathcal{X}$ with a non-complete topology and this is a bit unpleasant. To see this, just consider $0<\varepsilon<\frac{\sigma}{1-m} -(d-\gamma)$ and define the function
\[
f(x)= |x|^{-\frac{\sigma}{1-m}}\,(1-\chi_{B_1(0)}) + |x|^{-(d-\gamma)-\varepsilon}\,\chi_{B_1(0)}\,.
\]
The above function $f$ does not belong to $\LL^1_\gamma(\RR^d)$ nevertheless $|f|_{\mathcal{X}}<\infty$, hence $f\in \mathcal{X}$. Unfortunately, $f$ can be approximated in the topology induced by $|\cdot|_\mathcal{X}$, by the family $\{ f_r(x)\}_{r\in(0,1]}\subset\LL^1_\gamma(\RR^d)$, where:
\[
f_r(x) = |x|^{-\frac{\sigma}{1-m}}\,(1-\chi_{B_1(0)}) + |x|^{-(d-\gamma)-\varepsilon}\,\chi_{B_1(0)\setminus B_r(0)}\,,\quad0<r<1\,.
\]
We have that $f_r\in\mathcal{X}$ for any $0<r \le 1$ and a simple (but lengthy) computation shows that $|f_r-f|_\mathcal{X}\rightarrow0$ as $r \rightarrow 0$. Hence, we prefer to introduce the following norm on $\mathcal{X}$
\begin{equation}\label{norm.X}
\|f\|_{ \mathcal{X}}:=\sup_{R > 0} \left(1 \vee R\right)^{\frac{2+\beta-\gamma}{1-m}-(d-\gamma)}\int_{B^{c}_R(0)}|f(x)||x|^{-\gamma}\dx <\infty\,.
\end{equation}
The main difference between $|\cdot|_{\mathcal{X}}$ and $\|\cdot\|_{\mathcal{X}}$ is that the latter takes into account the influence of  the $\LL^1_\gamma(\RR^d)$-norm, and provides a complete topology, as in the following:
\begin{prop}\label{info.space.X}
Let $d\geq 3$, $\gamma, \beta < d$ real numbers such that  $\gamma -2 < \beta \leq \gamma(d-2)/d$ and $m \in \left(m_c, 1\right)$.  Then,
\begin{itemize}
\item[i)] For any $f\in\mathcal{X} $ we have that
\begin{equation}\label{norm.X.controls.L.1}
\|f\|_\mathcal{X}= \max\{\|f\|_{\LL^1_\gamma(\RR^d)}, |f|_\mathcal{X}\}\,;
\end{equation}
\item[ii)] $\mathcal{X}$ equipped with the norm $\|\cdot\|_\mathcal{X}$, defined in~\eqref{norm.X}, is a Banach space;
\item[iii)] Compactly supported functions are \textsl{not dense} in $\mathcal{X}$ equipped with the norm $\|\cdot\|_\mathcal{X}$.
\item[iv)] The Barenblatt profile has finite $\mathcal{X}$ norm, indeed
\begin{equation}\label{limit.infinity.X}
|\mathcal{B}_M|_\mathcal{X} = \lim\limits_{R\rightarrow\infty} R^{\frac{\sigma}{1-m}-(d-\gamma)}\,\int_{B_R^c(0)}\mathcal{B}_M |x|^{-\gamma}\dx=(1-m)\,\vartheta\,\omega_d\,.
\end{equation}
\end{itemize}
\end{prop}
The proof of the above Proposition is long but straightforward, hence we refrain from giving it here; however, it can be found in~\cite[Chapter 4]{Simonov2020}.

As we already explained before, the space $\LL^1_{\gamma, +}(\RR^d)$ can be split into two disjoint sets $\mathcal{X}$ and $\mathcal{X}^c$. A remarkable fact is that $\mathcal{X}$ and $\mathcal{X}^c$ are two invariant sets of $\LL^1_{\gamma, +}(\RR^d)$ for the fast diffusion flow, in a sense  made precise in the following Proposition.
\begin{prop}[Invariance of $\mathcal{X}$ and $\mathcal{X}^c$]\label{curve.X.prop}
Let $u(t)$ be a solution to~\eqref{cauchy.problem} with $ u_0 \in\LL^1_{\gamma, +}(\RR^d)$. Then,
\begin{itemize}
\item[i)]  $\mathcal{X}$ is invariant under the flow, namely $u_0 \in \mathcal{X}$ if and only if $u(t, \cdot)\in \mathcal{X}$ for all $t>0$\,.
\item[ii)]$\mathcal{X}^c$ is invariant under the flow, namely $u_0 \in \mathcal{X}^c$ if and only if $u(t, \cdot)\in \mathcal{X}^c$ for all $t>0$\,.
\item[iii)]   If $u_0 \in \mathcal{X}^c$, then
\[
\Big\|\frac{u(t)}{\mathcal{B}(t; M)}-1\Big\|_{\LL^\infty(\RR^d)}=\infty\qquad\mbox{for all $t>0$.}
\]

\item[iv)]  If $u_0 \in \mathcal{X}\setminus\{0\}$, then the following limit holds
\begin{equation}\label{limit.X.norms}
\lim_{t\rightarrow\infty}  \frac{|u(t)|_\mathcal{X}}{t^\frac{1}{1-m}} = \lim_{t\rightarrow\infty} \frac{\|u(t)\|_\mathcal{X}}{t^\frac{1}{1-m}} = \left(\sigma\,m\right)^\frac{1}{1-m}\, \left(\frac{\vartheta}{1-m}\right)^\frac{m}{1-m}\,\omega_d\,.
\end{equation}
\item[v)] If $u_0 \in \mathcal{X}\setminus\{0\}$, then for any $t>0$ the function $t\mapsto t^\frac{1}{m-1} |u(t, \cdot)|_\mathcal{X} $ is non increasing, and
\begin{equation}\label{lower.bound.X.norm}
\frac{|u(t-h)|_\mathcal{X} }{(t-h)^{\frac{1}{1-m}}}\ge\frac{ |u(t)|_\mathcal{X}}{t^{\frac{1}{1-m}}}\ge \,\left(\sigma\,m\right)^\frac{1}{1-m}\, \left(\frac{\vartheta}{1-m}\right)^\frac{m}{1-m}\,\omega_d\,,\qquad\mbox{for all $0\le h\le  t$.}
\end{equation}
\end{itemize}
\end{prop}
\begin{rem} \rm
Inequality~\eqref{lower.bound.X.norm} is sharp, because equality is achieved by the Barenblatt profile $\mathfrak{B}(t\,,\cdot\,; M)$ for any $M>0$.
\end{rem}

\noindent\textbf{Proof. }Recall first Lemma~\ref{lem.herrero.pierre.inverse}, whose proof is contained in the Appendix: for any $R, t, s \ge 0$ we have that
\begin{equation}\label{herrero.pierre.for.X}\begin{split}
\left(\int_{B_R^c(0)}u(t,x)\frac{\dx}{|x|^\gamma}\right)^{1-m}\le \left(\int_{B_{2R}^c(0)}u(s,x)\frac{\dx}{|x|^\gamma}\right)^{1-m} + C\,|t-s|\,R^{(d-\gamma)(1-m)-\sigma}\,,
\end{split}\end{equation}
where $C$ is a positive constant which only depends on $d, m, \gamma, \beta$. By taking the supremum in $R>0$ in inequality~\eqref{herrero.pierre.for.X} we can deduce that there exist $c>0$ depending only on $d,m, \gamma, \beta$, such that for any $t, s \ge0$
\begin{equation}\label{inequality.norms.X.nocontrol}
|u(t, \cdot)|_\mathcal{X}\le c\left(|u(s, \cdot)|_\mathcal{X} + |t-s|^\frac{1}{1-m}\right)\,.
\end{equation}
Let us prove first $i)$. From~\eqref{inequality.norms.X.nocontrol} we deduce that if $u_0 \in\mathcal{X}$ then $u(t) \in\mathcal{X}$ for all $t>0$ just by letting $s=0$. The opposite choice lead to the converse implication.
Part  $ii)$ follows analogously from~\eqref{inequality.norms.X.nocontrol}.
 To prove $iii)$ we proceed by contradiction. Suppose that there exists $\overline{t}$ such that
\[
\Big\|\frac{u(\overline{t})}{\mathcal{B}(\overline{t}; M)}-1\Big\|_{\LL^\infty(\RR^d)}\le C < \infty\,.
\]
Then, reasoning as in the proof of Theorem~\ref{equivalent.condition} (the necessary part, Subsection~\ref{sec:proof.thm.equivalent.condition}) we conclude that $u(\overline{t})\in \mathcal{X}$, and by point $i)$, we would have that $u_0 \in \mathcal{X}$, a contradiction.
\noindent As for the proof of $iv)$,  we will use the GHP of Theorem \ref{thm.ghp}: fix $t_0>0$, then there exist $\underline{\tau} ,\overline{\tau}>0$ and $\underline{M}, \overline{M}>0$ such that for any $t\ge t_0>0$  we have
\[
\mathfrak{B}(t-\underline{\tau},x; \underline{M}) \leq u(t, x) \leq \mathfrak{B}(t+\overline{\tau},x; \overline{M})\,,
\qquad\mbox{for any $x \in \RR^d$}\,.
\]
Then, a lengthy still not difficult computation, gives the following inequality
\begin{equation}\label{inequality.limit}
\left(\frac{R_\star(t-\underline{\tau})}{\zeta}\right)^\frac{1}{(1-m)\vartheta} \le \frac{|u(t, \cdot)|_\mathcal{X}}{\omega_d\,(1-m)\,\vartheta\,} \le\left(\frac{R_\star(t+\overline{\tau})}{\zeta}\right)^\frac{1}{(1-m)\vartheta}\,.
\end{equation}
Recall that $R_\star$ is defined in~\eqref{parameters.sigma.theta}, so that inequality~\eqref{limit.X.norms} follows from~\eqref{inequality.limit}, together with the following
\[
\lim_{t\rightarrow\infty} R_\star(t+T)^\frac{1}{(1-m)\vartheta}\,t^{-\frac{1}{1-m}} = \vartheta^{-\frac{1}{1-m}}\,, \qquad\mbox{for any $T\in \RR$.}
\]
We prove $v)$ through a smooth approximation by means of auxiliary norms: Let $k>0$ be a positive integer and  $\phi_k(x)$ be such that
\[
\phi_k(x)=1\,\mbox{on}\, |x|\ge 1 +\frac{1}{k}\,,\,\,\phi_k(x)=0\,\mbox{on}\,|x|\le 1\,,\quad\mbox{and}\,\,\phi_k(x)>0\,\mbox{on}\,1<|x|<1+\frac{1}{k}\,.
\]
Let us define
\begin{equation}\label{third.norm.X}
||| f |||_{k,\mathcal{X}} = \sup_{R>0} R^{\frac{\sigma}{1-m}-(d-\gamma)}\, \int_{\RR^d} f(x)\,\phi_k\left(\frac{x}{R}\right)|x|^{-\gamma}\dx\,.
\end{equation}
For any $k\ge1$ and for any $f\in \mathcal{X}$ we have that
\[
\left(\frac{k}{k+1}\right)^{\frac{\sigma}{1-m}-(d-\gamma)} \, |f|_\mathcal{X} \le |||f|||_{k, \mathcal{X}} \le |f|_\mathcal{X}\,,
\]
as a consequence of the above inequality for any $f\in\mathcal{X}$ the following limit holds
\begin{equation}\label{equivalence.tail.condition.third.norm}
\lim_{k\rightarrow \infty} |||f|||_{k, \mathcal{X}} =|f|_\mathcal{X}\,.
\end{equation}
We take advantage of the auxiliary norms~\eqref{third.norm.X}. Let $k >0$ be a positive integer and $R>0$ and define $Y_{k}(t)=\int_{\RR^d}\phi_k\left(\frac{x}{R}\right)\,u(t,x)\frac{\dx}{|x|^\gamma}$. Using now time monotonicity, the so-called Benilan-Crandall estimates \cite{Benilan1981}, $u_t\le \frac{u}{(1-m)t}$  valid in the distributional sense, we find that
\[
Y_k'(t)\le \frac{1}{(1-m)t} \,Y_k(t)\,,\qquad\mbox{so that for all $ \tau > s >0$}\qquad
\frac{Y_k(s)}{s^{\frac{1}{1-m}}} \ge \frac{Y_k(\tau)}{\tau^{\frac{1}{1-m}}}\,.
\]
Multiplying by $R^{\sigma/(1-m)-(d-\gamma)}$ and taking the supremum in $R>0$ in the above inequality we get
\[
\frac{|||u(s,\cdot)|||_{k, \mathcal{X}}}{s^\frac{1}{1-m}} \ge \frac{|||u(t,\cdot)|||_{k, \mathcal{X}}}{t^\frac{1}{1-m}} \,,
\]
taking the limit as $k\rightarrow \infty$ in the above inequality one gets the monotonicity of $t^{-1/(1-m)}|u(t, \cdot)|_\mathcal{X}$. Let $\tau > s >0$ as before (we can take any $\tau>s$, hence also $\tau\to\infty$) to obtain inequality~\eqref{lower.bound.X.norm}, namely
\[
\frac{|u(s, \cdot)|}{s^{\frac{1}{1-m}}} \ge \frac{|u(\tau, \cdot)|}{\tau^{\frac{1}{1-m}}}\ge \lim_{\tau\to \infty}\frac{|u(\tau, \cdot)|}{\tau^{\frac{1}{1-m}}}
=\left(\sigma\,m\right)^\frac{1}{1-m}\, \left(\frac{\vartheta}{1-m}\right)^\frac{m}{1-m}\,\omega_d\,,
\]
where in the last step we have used~\eqref{limit.X.norms}.  The proof is then concluded.\qed

\subsection{Convergence to the Barenblatt in $\mathcal{X}$}
Finally we address here the question of convergence to the Barenblatt profile of solutions to~\eqref{cauchy.problem} in $\mathcal{X}$ with the topology induced by $\|\cdot\|_\mathcal{X}$. We find that it is false in general that
\[
\|u(t, \cdot)-\mathfrak{B}(t, \cdot\,; M)\|_\mathcal{X} \rightarrow 0 \quad \mbox{as} \quad t\rightarrow \infty\,.
\]
To see this fact we provide an explicit counterexample. Consider the Barenblatt solution $\mathfrak{B}(t, x; M)$ and its translation in time $\mathfrak{B}(t+\tau, x; M)$, for $R>0$ large enough we have that
\[
\Big|\mathfrak{B}(t+\tau, x; M)-\mathfrak{B}(t, x; M) \Big|\ge \frac{1}{2}\left[\left(1+\frac{\tau}{t}\right)^{\frac{1}{1-m}} -1\right] \, \mathfrak{B}(t, x; M)\quad \mbox{for any }\quad |x|\ge R\,,
\]
we therefore conclude, thanks to \eqref{limit.infinity.X} and~\eqref{lower.bound.X.norm}, that
\[
\Big|\mathfrak{B}(t+\tau, x; M)-\mathfrak{B}(t, x; M) \Big|_\mathcal{X} \gtrsim t^\frac{m}{1-m}\,.
\]
However, if we suitably renormalize  $\|\cdot\|_\mathcal{X}$ by the factor $t^\frac{1}{1-m}$ we find the following result.
\begin{prop}\label{convergence.barenblatt.Xnorm}
Under the assumption of Theorem~\ref{equivalent.condition} we have that
\begin{equation}\label{convergence.barenblatt.Xnorm.limit}
\lim\limits_{t\rightarrow\infty}\frac{\|u(t, \cdot)-\mathfrak{B}(t, \cdot; M)\|_\mathcal{X}}{t^\frac{1}{1-m}}=0\,.
\end{equation}
\end{prop}
\begin{rem} \rm
It is interesting to stress that the above limit, rewritten  in \emph{self-similar} variables~\eqref{change.fokker.planck},  does not need the renormalization factor $t^\frac{1}{1-m}$, namely
\[
\lim\limits_{\tau\rightarrow\infty}\|v(\tau, \cdot)-\mathcal{B}_M(\cdot)\|_\mathcal{X}=0\,.
\]
\end{rem}
\noindent\textbf{Proof.} By Proposition~\ref{info.space.X} we know that $\|\cdot\|_\mathcal{X}=\max\{|\cdot|_\mathcal{X},\, \|\cdot\|_{\LL^1_\gamma(\RR^d)}\}$. To prove~\eqref{convergence.barenblatt.Xnorm} we need to show
 \[
 \lim_{t\rightarrow\infty}\frac{|u(t, \cdot)-\mathfrak{B}(t, \cdot; M)|_\mathcal{X}}{t^\frac{1}{1-m}}=0\qquad\mbox{and}\qquad\lim_{t\rightarrow\infty}\frac{\|u(t, \cdot)-\mathfrak{B}(t, \cdot; M)\|_{\LL^1_\gamma(\RR^d)}}{t^\frac{1}{1-m}}=0\,.
 \]
 By conservation of mass the second limit is straightforward, hence we only need to prove the first. Under the running assumption we know by Theorem \ref{equivalent.condition} that $u(t, x)$ converge to $\mathfrak{B}(t, x; M)$ uniformly in relative error. We restate this result in the following way: there exists a positive function  $g(t)\rightarrow $0 as $t\rightarrow \infty$ such that for any $x\in \RR^d$  and for any $t$ large enough we have
\[
|u(t, x)-\mathfrak{B}(t, x; M)| \le g(t)\, \mathfrak{B}(t, x; M)\,.
\]
By the above inequality we deduce that
\[
\limsup_{t\rightarrow\infty} \frac{|u(t, \cdot)-\mathfrak{B}(t, \cdot; M)|_\mathcal{X}}{t^\frac{1}{1-m}} \le \lim_{t\rightarrow \infty} g(t)\,\frac{|\mathfrak{B}(t, \cdot; M)|_\mathcal{X}}{t^\frac{1}{1-m}} =0\,,
\]
where we have used identity~\eqref{limit.X.norms} and the fact that $g(t)\rightarrow 0$ as $t\rightarrow \infty$. The proof is therefore concluded.\qed

\section{Appendix}
\subsection{How to recover the Mass of the Barenblatt profile $\mathcal{B}_M$}\label{appendix:mass}
The following identity is a consequence of the integral representation formula of the Euler Beta function, see~\cite[6.2.1, pag. 258]{Abramowitz1964}:
\begin{equation*}
\overline{M}=\int_{\RR^d}\left(1+|x|^{2+\beta-\gamma}\right)^\frac{1}{m-1} |x|^{-\gamma} \dx = |\mathbb S^{d-1}|\,\frac{\Gamma \left(\frac{d-\gamma }{2+\beta-\gamma}\right) \Gamma \left(\frac{1}{1-m}-\frac{d-\gamma }{2+\beta-\gamma }\right)}{(2+\beta-\gamma)\;\Gamma
   \left(\frac{1}{1-m}\right)}\,.
\end{equation*}
By scaling we obtain
\begin{equation*}\begin{split}
M  &= \int_{\RR^d}\left(C(M)+|x|^{2+\beta-\gamma}\right)^\frac{1}{m-1}|x|^{-\gamma}\dx = C(M)^{\frac{1}{m-1}-\frac{\gamma}{\sigma}} \int_{\RR^d}\left(1+|x/C(M)^{1/\sigma}|^{\sigma}\right)^\frac{1}{m-1}|x/C(M)^{1/\sigma}|^{-\gamma}\dx \\
& = C(M)^{\frac{1}{m-1}+\frac{d-\gamma}{\sigma}}  \int_{\RR^d}\left(1+|x|^{2+\beta-\gamma}\right)^\frac{1}{m-1} |x|^{-\gamma} \dx  = C(M)^{\frac{1}{m-1}+\frac{d-\gamma}{\sigma}} \overline{M},
\end{split}\end{equation*}
therefore we have that $ C(M)= \left(\frac{\overline{M}}{M}\right)^{\frac{\sigma(1-m)}{\sigma-(d-\gamma)(1-m)}}$.

\subsection{Interpolation Inequality}
Let $\Omega\subset \RR^d$ be a bounded domain and  $u: \Omega\rightarrow \RR$ be a function and define for any $\nu \in \left(0,1\right)$
\begin{equation}\label{C.NU.norm}
\lfloor u\rfloor_{C^\nu\left(\Omega\right)} := \sup_{\substack{x,y \in \Omega \\ x\neq y}} \frac{|u(x)-u(y)|}{|x-y|^\nu}\,.
\end{equation}
We say that $u \in C^\nu(\Omega)$ whenever $\lfloor u\rfloor_{C^\nu\left(\Omega\right)}< \infty$. Notice that  $\lfloor u\rfloor_{C^\nu\left(\Omega\right)}=0$ if and only if $u$ is constant, since in what follows we need to use strictly positive quantities we shall use the following inequality which hold for $u\in C^\nu\left(\Omega\right)$
\begin{equation}\label{C.nup1.norm}
|u(x)-u(y)| \le \left(1+\lfloor u\rfloor_{C^\nu\left(\RR^d\right)}\right) \, |x-y|^\nu\,.
\end{equation}
Let $\Omega' \subset \Omega$ be a subdomain, we define the distance between $\Omega$ and $\Omega'$ as
\[
\mathrm{dist}(\Omega, \Omega') = \inf_{\substack{x \in \partial\Omega, \\ y \in \partial\Omega'}} |x-y|\,,
\]
where $\partial\Omega$ is the boundary of $\Omega$ and $\partial\Omega'$ is the boundary of $\Omega'$. The purpose of this appendix is to prove the following lemma.
\begin{lem}\label{interpolation.lemma}
Let $p \geq 1$,  $\nu \in \left(0,1\right)$ and $u: \Omega  \rightarrow \RR$ be a function such that $u \in \LL^p_\gamma(\Omega)\cap C^\nu\left(\Omega\right)$. Assume $\gamma\le0$ and let $\Omega'\subset\Omega$ be  such that $\mathrm{dist}(\Omega, \Omega')>0$, then there exists a positive constant $C_{d,\gamma, \nu, p}$, which depends on $d, \gamma, p$ and $\nu$, such that
\begin{equation}\label{interpolation.inequality.gamma.negative}
\|u\|_{\LL^\infty(\Omega')} \leq C_{d,\gamma, \nu, p}\, \left(1+\frac{\|u\|_{\LL^p_\gamma(\Omega)}}{\left(1+\lfloor u\rfloor_{C^\nu\left(\Omega\right)}\right)\,\mathrm{dist}(\Omega, \Omega')^\frac{1}{p}}\right)^\frac{d-\gamma}{d-\gamma+p\nu} \left(1+ \lfloor u\rfloor_{C^\nu\left(\Omega\right)}\right)^{\frac{d-\gamma}{d-\gamma+p\nu}} \|u\|_{\LL^p_\gamma\left(\Omega\right)}^{\frac{p\nu}{d-\gamma+p\nu}}\,.
\end{equation}

Assume $0<\gamma<d$, and let in addiction $\Omega'\subset\Omega$ be two bounded domains, then there exists a positive constant $C_{d,\gamma, \nu, p}$, which depends on $d, \gamma, p$ and $\nu$, such that
\begin{equation}\label{interpolation.inequality.gamma.positive}\begin{split}
\|u\|_{\LL^\infty(\Omega')} &\le
C_{d,\gamma, \nu, p}\, \left(\mathrm{dist}(\Omega, \Omega')+\sup_{x\in\Omega' }|x|\right)^\frac{\gamma}{p}\,
\left(1+\frac{\|u\|_{\LL^p_\gamma(\Omega)}}{\left(1+\lfloor u\rfloor_{C^\nu\left(\Omega\right)}\right)\,\mathrm{dist}(\Omega, \Omega')^\frac{1}{p}}\right)^\frac{d }{d +p\nu}\\
&\times\left(1+\lfloor u\rfloor_{C^\nu\left(\Omega\right)}\right)^\frac{d}{d+p\nu}\,\|u\|_{\LL^p_\gamma(\Omega)}^\frac{p\nu}{d+p\nu}\,,
\end{split}\end{equation}
\end{lem}
\noindent\textbf{Proof}. For any  $x, y \in \Omega'$ we have, by the triangle inequality, that
\begin{equation*}
|u(x)|^p \leq \left(|u(x)-u(y)|+|u(y)|\right)^p \leq 2^p \left(|u(x)-u(y)|^p +|u(y)|^p\right)\,.
\end{equation*}
Let $x\in \Omega'$, $0\le R < \mathrm{dist}(\Omega, \Omega')$, averaging on a ball $B_R(x)\subset \Omega$  we have
\begin{equation}\label{to.optimize}\begin{split}
|u(x)|^p &\leq  \frac{2^p }{\mu_\gamma(B_R(x))} \int_{B_R(x)}|u(x)-u(y)|^p|y|^{-\gamma}\dy + \frac{2^p }{\mu_\gamma(B_R(x))} \int_{B_R(x)}|u(y)|^p|y|^{-\gamma}\dy   \\
&\leq 2^{p}  R^{p\nu} \left(1+\lfloor u\rfloor_{C^\nu\left(\Omega\right)}\right)^p
+ 2^{p}\,\frac{\|u\|^p_{\LL^p_\gamma(\Omega)}}{\mu_\gamma(B_R(x))}\,,
\end{split}\end{equation}
where in the last step we have used~\eqref{C.nup1.norm} and that $\int_{B_R(x)}|u(y)|^p|y|^{-\gamma}\dy \leq\|u\|^p_{\LL^p_\gamma(\Omega)}$. We claim that for any $x_0 \in \RR^d$ and for any $R\ge0$ there exist positive constants $c_{\gamma, d}, C_{\gamma, d}$ such that
\begin{equation}\label{measure.inequality}
c_{\gamma, d}\, R^{d} \, \left(|x_0|\vee \frac{R}{2}\right)^{-\gamma} \le \mu_\gamma(B_R(x_0)) \le C_{\gamma, d}\, R^{d} \, \left(|x_0|\vee \frac{R}{2}\right)^{-\gamma}\,.
\end{equation}
The above inequality can be proven using the techniques developed in \cite[Lemma 5.2, Appendix B]{Bonforte2019a}, we will not include the proof here. Now we consider two cases, namely $\gamma\le 0$ and $\gamma>0$.

\noindent Assume that $\gamma\le0$, plugging the lower bound of~\eqref{measure.inequality} in~\eqref{to.optimize} we deduce
\begin{equation}\label{evaluate.gamma.negative}
|u(x)|^p \leq C\left( R^{p\nu} \, \left(1+\lfloor u\rfloor_{C^\nu\left(\Omega\right)}\right)^p
+ R^{\gamma-d}\,\|u\|^p_{\LL^p_\gamma(\Omega)}\right)\,.
\end{equation}
Then, inequality~\eqref{interpolation.inequality.gamma.negative} follows by letting
\[
2 R = \left(\frac{(d-\gamma)\,\|u\|^p_{\LL^p_\gamma(\Omega)}}{p\nu\,\left(1+\lfloor u\rfloor_{C^\nu\left(\Omega\right)}\right)^p}\right)^\frac{1}{d-\gamma+p\nu}\wedge \mathrm{dist}(\Omega, \Omega')\,.
\]

Assume that $0<\gamma<d$, then using~\eqref{measure.inequality} in~\eqref{to.optimize} we have for any $0\le R< \mathrm{dist}(\Omega, \Omega'):=D$
\begin{equation}\label{evaluate.gamma.positive}\begin{split}
|u(x)|^p &\le C\left( R^{p\nu} \, \left(1+\lfloor u\rfloor_{C^\nu\left(\Omega\right)}\right)^p
+ \frac{\left(|x_0|\vee \frac{R}{2}\right)^\gamma}{R^d}\,\|u\|^p_{\LL^p_\gamma(\Omega)}\right) \nonumber \\
&\le C \left(D+|x_0|\right)^\gamma\, \left( R^{p\nu} \, \left(1+\lfloor u\rfloor_{C^\nu\left(\Omega\right)}\right)^p
+ R^{-d}\,\|u\|^p_{\LL^p_\gamma(\Omega)}\right)
\end{split}\end{equation}
where in the last step we have used that $\left(|x_0|\vee \frac{R}{2}\right)^\gamma \le \left(D+|x_0|\right)^\gamma$. Letting
\[
2R  = \left(\frac{d\,\|u\|^p_{\LL^p_\gamma(\Omega)}}{p\nu\, \left(1+\lfloor u\rfloor_{C^\nu\left(\Omega\right)}\right)^p}\right)^\frac{1}{d+p\nu}\wedge \mathrm{dist}(\Omega, \Omega')\,,
\]
taking the supremum in $x_0 \in \Omega'$ we get~\eqref{interpolation.inequality.gamma.positive}. The proof is then complete. \qed

\subsection{Holder Continuity of solution to weighted equations}\label{app.holder.cont}
We present here some regularity results for nonnegative local weak solutions to both linear and nonlinear parabolic equations with weights. The results contained in this section are mainly contained in \cite[Part III]{Bonforte2019a} and references therein. We provide similar results here, adapted to the present setting and  assumptions, for convenience of the reader. Consider local weak solutions in the cylinder $Q:=(0,T)\times\Omega$  to the equation
\begin{equation}\label{WHE.DW.GEN}
 v_t=w_\gamma\sum_{i,j=1}^N\partial_i  \left( A_{i,j}(t,x)\,  \partial_j v  \right),
\end{equation}
where $A_{i,j}=A_{j,i}$ and there exist constants $0 < \lambda_0 \le \lambda_1 < +\infty$ such that for some $\gamma,\beta < N$, satisfying $\gamma -2 < \beta \leq \left(\frac{N-2}{N}\right) \gamma$, we have for any $\xi \in \RR^d$ and any $x \in \RR^d$
\begin{equation}\label{WHE.DW.GEN.paramt}
w_\gamma\asymp |x|^{\gamma} \qquad\mbox{and}\qquad 0<\lambda_0|x|^{-\beta} |\xi|^2\le \sum_{i,j=1}^N A_{i,j}(t,x)\xi_i \xi_j\le \lambda_1 |x|^{-\beta}|\xi|^2\,.
\end{equation}
We shall restrict ourselves to the class of  bounded, nonnegative, local weak solutions to equation~\eqref{WHE.DW.GEN}, precisely defined in  \cite{Bonforte2019a,Chiarenza1984b,Gutierrez1991}. Notice that this class of solutions is large enough for our purposes.

It is convenient to introduce the notion of distance between nested cylinders of the form $Q=(0, T)\times \Omega$. Let $Q'=(T_1, T_2) \times \Omega' \subset Q$, we define
\begin{equation}\label{parabolic.distance.sets}
d_{\gamma, \beta}(Q, Q'):= \inf\limits_{\substack{(t,x)\in \{[0,T]\times\partial\Omega\}\cup\{\{0\}\times\Omega\},\\ (s,y) \in Q'}}\, |x-y|\vee \left( \rho^{\gamma, \beta}_y\right)^{-1}(|t-s|)\,.
\end{equation}
where  $\gamma, \beta$ are as above and $\left( \rho^{\gamma, \beta}_y\right)^{-1}$ is the inverse of  $\rho_{y}^{\gamma,\beta}$ (well defined for any $y \in \RR^d$) defined as
\[
 \rho_{y}^{\gamma,\beta}(R) := \left( \int_{B_R(y)} |x|^{(\beta - \gamma) \frac{N}{2}} \dx \right)^{\frac{2}{N}}\,.
\]
Finally, we introduce the a notion of $C^\alpha$ norm which takes into account the presence of the weights. With the above notation we define
\begin{equation}\label{C.NU.weithed.nomr}
\lfloor u\rfloor_{C^\alpha_{\gamma, \beta}\left(Q\right)}:=\sup_{\substack{(t,x), (\tau,y)\in Q'\\ (t,x), \neq(\tau,y)}}\frac{|v(t,x)-v(\tau,y)|}{(|x-y|+|t-\tau|^{\frac{1}{2\vee \sigma}})^\alpha}
\end{equation}
The proof of the following result can be found in \cite[Proposition 4.2, Corollary 4.3]{Bonforte2019a}.

Notice that  the following results involve both the ``parabolic'' H\"older norm $\lfloor u\rfloor_{C^\alpha_{\gamma, \beta}\left(Q'\right)}$ defined in \eqref{Prop.Lin.Harn.HoCont.ineq.111}, and the ``elliptic'' one, $\lfloor u\rfloor_{C^\nu\left(\Omega\right)}$, defined in \eqref{C.NU.norm}.
\begin{prop}[H\"older Continuity for linear equations with weights]\label{Prop.Lin.Harn.HoCont}
Let $v$ be a nonnegative bounded local weak solution to equation \eqref{WHE.DW.GEN} on $Q:=(0,T)\times\Omega$, under the assumption \eqref{WHE.DW.GEN.paramt}. Let $Q':=(T_1, T_2) \times \Omega' \subset Q$. Then there exist $\alpha\in (0,1) $  and $\ka_\alpha >0$\,, such that for all $(t,x), (s,y) \in Q' $
\begin{equation}\label{Prop.Lin.Harn.HoCont.ineq.111}
\lfloor u\rfloor_{C^\alpha_{\gamma, \beta}\left(Q'\right)}\le\frac{\ka_\alpha}{d_{\gamma, \beta}(Q, Q')^\alpha} \|v\|_{\LL^\infty(Q)},
\end{equation}
where $\ka_\alpha>0$ is given by
\begin{equation}\label{ka.alpha}
\ka_\alpha= \ka_\alpha' \left\{
\begin{array}{lll}
1\,,&\,\mbox{if }  \sigma \ge 2 ,\\
\left(T^{\frac{1}{\sigma}}\vee \sup\limits_{x_0\in\Omega}|x_0|\right)^{\frac{\gamma-\beta}{2}}\,,&\,\mbox{if }  0<\sigma < 2\,.\\
\end{array}\right.
\end{equation}
The constants $\alpha, \ka_\alpha' $depend only on $N, \gamma, \beta, \lambda_0, \lambda_1$.
\end{prop}
Proposition~\ref{Prop.Lin.Harn.HoCont.ineq.111} can be fruitfully used to deduce regularity results also for nonlinear parabolic equation: for example we can consider nonnegative bounded solutions to $u_t = |x|^\gamma \nabla\left(|x|^{-\beta}\nabla u^m\right)$ as solutions to the linear equation $u_t = |x|^\gamma \nabla\left(|x|^{-\beta}\, a(t,x)\,\nabla u\right)$ where $a(t,x)=m u(t,x)^{m-1}$. Indeed the same can be done for solutions to the Fokker-Planck type equation~\eqref{fokker.planck.equation} as follows.
\begin{lem}\label{Technical.Lemma}
Let $\rho, \tau_0 > 0$, $0 <\lambda_0 \le \lambda_1 < \infty$, $m \in (0,1)$ and let $v(\tau,y):(0, \infty) \times \RR^d \rightarrow \RR$ be a nonnegative bounded solution to~\eqref{fokker.planck.equation}, assume that
\begin{equation*}
\lambda_0\leq m\,v^{m-1}(\tau,y) \leq \lambda_1 \,\,\,\,\mbox{for any}\,\,\, \tau \geq \tau_0 \,\,\,\mbox{and}\,\,\,\ |y|\leq \rho\,.
\end{equation*}
Then there exist $\nu>0$ and $\ka >0$ such that if $ \tau_1 > \tau_0$ and $\tau \in \left[\tau_1+\frac{1}{\sigma}\log{R_\star(2)},\tau_1+\frac{1}{\sigma}\log{R_\star(3)}\right]$ then
\begin{equation}\label{c.nu.selfsimilar}
\lfloor v(\tau, \cdot) \rfloor_{C^\nu\left(B_{\rho/2}(0)\right)}\leq \ka\,\,\|v\|_{\LL^{\infty}\left(\left[\tau_1,\tau_1+\frac{1}{\sigma}\log{R_\star(4)}\right] \times B_{\rho}(0)\right)}\,.
\end{equation}
The constants $\nu, \ka$ depend  on $m, d, \gamma, \beta, \lambda_0, \lambda_1$; $\ka$ depends also on $\rho$.
\end{lem}
\noindent {\bf Proof.~}The proof is divided in several steps. It is convenient to consider a time-shifted solution:
\begin{equation*}
\overline{v}\left(\tau,y\right):= v\left(\tau+\tau_1\right) \,\,\,\,\mbox{for any}\,\,\,\tau\geq0\,.
\end{equation*}
\noindent\textit{Rescaling to originals variables}.
The rescaled function $\overline{u}(t,x)$ defined by
\begin{equation}
\overline{u}(t,x):=\frac{\zeta^{d-\gamma}}{R_\star(t+1)^{d-\gamma}} \overline{v}\left(\frac{1}{\sigma}\log \frac{R_\star(t+1)}{R_\star(1)}, \frac{\zeta x}{R_\star(t)}\right)= \frac{\zeta^{d-\gamma}}{R_\star(t+1)^{d-\gamma}} \overline{v}\left(\tau, y\right),
\end{equation}
satisfies~\eqref{cauchy.problem}. Define the following domains
\begin{equation*}
\overline{Q}_1:=\left\{(t, x) : 0\leq t \leq 3, |x|\leq \frac{\rho\,R_\star(t+1)}{\zeta}\right\}\,,\,\,\,\,
\overline{Q}_2:=\left\{(t, x) : 1\leq t \leq 2, |x|\leq \frac{\rho\,R_\star(t+1)}{2\,\zeta}\right\}\,\,\,\,
\end{equation*}
On both $\overline{Q}_1$ and $\overline{Q}_2$ the following estimate holds true
\begin{equation*}
\frac{R_\star(1)^{(d-\gamma)(1-m)}}{\zeta^{(d-\gamma)(1-m)}}\, \lambda_0  \leq m \tilde{u}^{m-1}(t,x) \leq \lambda_1 \, \frac{R_\star(4)^{(d-\gamma)(1-m)}}{\zeta^{(d-\gamma)(1-m)}}
\end{equation*}
\noindent\textit{Application of the linear result}. We can consider $\overline{u}$ as a bounded solution to the linear equation
\[
\overline{u}_t = |x|^{\gamma}\,\nabla\left( |x|^{-\beta}\,a(t,x)\, \nabla \overline{u}\right)\,\,\,\mbox{where}\,\,\,\, a(t,x)=m\overline{u}(t,x)\,,
\]
on the domain $\overline{Q}_1$. From Proposition~\ref{Prop.Lin.Harn.HoCont} we deduce that there exists $\nu>0$ and $ \ka_\nu>0$ such that
\begin{equation}\label{Holder.Norm.Estimates}
\|\overline{u}\|_{C^{\nu}_{\gamma, \beta}(\overline{Q}_2)}\leq \ka_\nu \frac{\|\overline{u}\|_{\LL^{\infty}(\overline{Q}_1)}}{d_{\gamma,\beta}(\overline{Q}_1, \overline{Q}_2)^{\nu}}\,.
\end{equation}
The constant $\nu$ shall depend only on $d, m, \gamma, \beta$ and $\lambda_0, \lambda_1$, since $R_\star(1), R_\star(4), \zeta$ are numerical constants which only depend  on $d, m, \gamma, \beta$. However, the constant $\ka_\nu$ will depend on $\rho$ when $0< \sigma<2$, see the expression of the constant $\ka_\alpha$ in Proposition~\ref{Prop.Lin.Harn.HoCont}. Finally, we notice that $d_{\gamma,\beta}(\overline{Q}_1, \overline{Q}_2)^{\nu}$ depend as well on $\rho$.
We shall now freeze the time variable and consider $\overline{u}(t,x)$ as a function in space only. From~\eqref{Holder.Norm.Estimates} we deduce that for any $ t \in \left[1,2\right] $ we have that
\begin{equation}\label{c.nu.original.variables}
\lfloor u(t, \cdot) \rfloor_{C^\nu\big(B_{\widetilde{R}(t)}(0)\big)} \le \ka_\nu \frac{\|\overline{u}\|_{\LL^{\infty}(\overline{Q}_1)}}{d_{\gamma,\beta}(\overline{Q}_1, \overline{Q}_2)^{\nu}}
\qquad\mbox{where}\qquad \widetilde{R}(t)=\frac{\rho\,R_\star(t+1)}{2\,\zeta}.
\end{equation}
\noindent\textit{Rescaling back to self-similar variables}. The domains $\overline{Q}_1$ and $\overline{Q}_2$ will be  back  to $\left[\tau_1,\tau_1+\frac{1}{\sigma}\log{R_\star(4)}\right] \times B_{\rho}(0)$ and  $\left[\tau_1+\frac{1}{\sigma}\log{R_\star(2)},\tau_1+\frac{1}{\sigma}\log{R_\star(3)}\right] \times B_{\rho/2}(0)$ respectively. While~\eqref{c.nu.original.variables} become~\eqref{c.nu.selfsimilar} where
\[
\ka = \left(\frac{R_\star(4)}{R_\star(2)}\right)^{d-\gamma}\,\left(\frac{\zeta}{R_\star(1)}\right)^\nu\, \frac{\ka_\nu}{d_{\gamma,\beta}(\overline{Q}_1, \overline{Q}_2)^{\nu}}\,.
\]
The proof is then concluded. \qed

\noindent {\large \sc Acknowledgments. }This work was partially funded by Projects MTM2017-85757-P (Spain) and by the E.U. H2020 MSCA programme, grant agreement 777822. N.~S. was partially funded by the FPI-grant BES-2015-072962, associated to the project MTM2014-52240-P (Spain). This work has been partially supported by the Project EFI ANR-17-CE40-0030 of the French National Research Agency.

\smallskip\noindent {\sl\small\copyright~2020 by the authors. This paper may be reproduced, in its entirety, for non-commercial purposes.}

\addcontentsline{toc}{section}{~~~References}
\bibliographystyle{siam}

\end{document}